\RequirePackage{fix-cm}
\documentclass[smallextended]{svjour3}       
\smartqed  
\usepackage{graphicx}
\usepackage{amssymb}
\usepackage{isomath}
\usepackage{mathtools}
\usepackage{tikz}
\usepackage{pgfplots}
\usetikzlibrary{arrows.meta}
  \pgfplotsset{compat=newest}
  \usetikzlibrary{plotmarks}
  \usepackage{grffile}

\pgfplotsset{every axis/.append style={
        scaled ticks = false, 
        tick label style={/pgf/number format/fixed}
    }
}

\usetikzlibrary{external}
\newlength\figureheight 
\newlength\figurewidth 
\usetikzlibrary{patterns}
\usetikzlibrary{plotmarks}
\usetikzlibrary{patterns}
\usepackage{tikz-3dplot}
\usepackage{multirow}
\usepackage{subfigure}
\usepackage{esint} 
\usepackage{multirow}
\usepackage{booktabs}
\usepackage{ulem}
\usepackage{cancel}
\usepackage[autostyle]{csquotes} 
\usepackage{stackengine} 
\usepackage{float}
\usepackage{setspace}
\usepackage{placeins}
\usepackage{lipsum}
\usepackage{amsfonts}
\usepackage{stmaryrd} 
\usepackage{algorithm2e}
\usepackage{adjustbox}
\usepackage{hyperref}

\begin{document}

\title{The Hermite-Taylor Correction Function Method for Maxwell's Equations
}

\titlerunning{The Hermite-Taylor Correction Function Method for Maxwell's Equations}        

\author{Yann-Meing Law \and Daniel Appel\"{o}
}


\institute{Y.-M. Law  \and  D. Appel\"{o} \at
              Department of CMSE, Michigan State University, East Lansing, USA \\
              \email{lawkamci@msu.edu, appeloda@msu.edu}           
}

\date{Received: date / Accepted: date}

\maketitle

\begin{abstract}
The Hermite-Taylor method, introduced in 2005 by Goodrich, $\phantom{Hagstrom}$ Hagstrom and Lorenz, is highly efficient and accurate when applied to linear hyperbolic systems on periodic domains. Unfortunately its widespread use has been prevented by the lack of a systematic approach to implementing boundary conditions. In this paper we present the Hermite-Taylor Correction Function method, which provides exactly such a systematic approach for handing boundary conditions. Here we focus on Maxwell's equations but note that the method is easily extended to other hyperbolic problems.   
\end{abstract}

\section{Introduction}
The property of waves to travel over large distances and long time without changing their shape is an important feature used in current technologies, such as communication devices and other electromagnetic products. The governing equations for electromagnetic problems are the Maxwell's equations and it is to these we seek approximate solutions in this paper. To make the numerical approximation to the solution accurate either low order methods on fine meshes, which can be computationally costly, or high-order methods on coarser meshes can be used. The latter approach is usually preferable for large scale problems.

Several high-order methods in computational electromagnetics have been proposed, such as high-order finite-difference time-domain (FDTD) methods \cite{Yee1966,Xie2002}, {\color{black} discontinuous Galerkin (DG) methods \cite{Hesthaven2002,cockburn2001runge,Balsara2019,Hazra2019}} and pseudo-spectral methods \cite{Fan2002,Galagusz2016,Yang1997}, to name a few. High-order explicit FDTD methods require a restrictive stability condition and wide stencils, which complicate the enforcement of boundary conditions. Unconditionally stable alternating-direction-implicit (ADI) FDTD methods have been developed to circumvent the time step constraints \cite{Namiki1999,Zheng1999,Tan2008,Chen2010,Liang2013}, however, methods are difficult to generalize to high order and treating complex geometry is not straightforward. 
 
Discontinuous Galerkin methods achieve high-order convergence rates by approximating the function using local high order polynomials and are an excellent choice for problems where a high quality mesh can be generated. The main drawbacks of DG methods is their restrictive time step at high order of accuracy and the duplication of degrees of freedom on the edges of elements. 	
	 
Another avenue to handle time dependent wave problems is the Hermite-Taylor method, which consists of a Hermite interpolation procedure in space and a Taylor method in time \cite{Goodrich2005} {\color{black}(see also \cite{hagstrom2015solving} for a review of Hermite methods)}. The key idea is to evolve, in time, the numerical solution as well as {\color{black} its space derivatives through order $m$} to achieve a $(2\,m+1)$ order accurate method using only $(m+1)^d$ degrees of freedom per element in $d$-dimensions. 

As was shown for linear symmetric hyperbolic problems in \cite{Goodrich2005}, this method provides a stability condition that only depends on the largest wave-speed, independent of the order. Hence, large time-step sizes can be used for these high-order methods and therefore ease the computational burden for large-scale problems. As the $(m+1)^d$ degrees of freedom in a Hermite method are collocated at a single node the imposition of general boundary conditions can be challenging. Typically, in addition to the physical boundary conditions the method needs to be augmented with a relatively large number of numerical boundary conditions {\color{black}(sometimes called compatibility boundary conditions or, more recently, inverse Lax-Wendroff conditions)}. While this has been successfully done for the wave equation on both Cartesian and curvilinear meshes in \cite{compat_wave_hermite_AAL_DEAA_WDH}, it has proven difficult to use this technique for first order hyperbolic systems.  

A possible solution to this is to use a hybrid DG-Hermite method \cite{Chen2014} for Maxwell's equations. The method in \cite{Chen2014}  takes advantage of the flexibility of DG solvers to handle complex geometries and boundary conditions by considering two non-overlapping meshes, an unstructured mesh for the DG method and a staircased Cartesian mesh where the Hermite method is used. This approach requires a hybrid structured-unstructured mesh and the use of local time-stepping to maintain large time-step sizes in the Hermite method. In \cite{OversetHermiteDG} an overset grid method that combines a Hermite method (on Cartesian meshes) and a DG method (on structured curvilinear meshes) for the wave equation is proposed. This method does not require a hybrid non-overlapping mesh and as such it is somewhat more geometrically flexible but again, it is not easy to extend to first order hyperbolic systems.  	

In this work, we propose an alternative solution for imposing boundary conditions for Maxwell's equations within the framework of Hermite methods.  Our new method is based on the correction function method (CFM). The CFM was first proposed in \cite{Marques2011} to handle Poisson's equation with interface conditions and continuous coefficients in a finite-difference context. Given a numerical solution (for example from a finite difference method) that has been updated near but not on the boundary from the CFM seeks a polynomial approximation to the solution in the vicinity of a boundary or interface using a minimization procedure. A functional that is based on a square measure of the residual of the original PDE problem and that also contains terms from the finite difference solver is minimized over a suitable space of polynomials. Once this polynomial approximation, also called the correction function, is found, the numerical solution can be corrected so that it satisfies the boundary conditions to high order of accuracy. The CFM method has been used for Poisson's equation \cite{Marques2017,Marques2019}, the wave equation \cite{Abraham2018} and for electromagnetic problems with both interface and boundary \cite{LawMarquesNave2020,LawNave2021,LawNave2022}. 

In this paper we introduce a CFM - Hermite-Taylor method. An advantage with using Hermite based methods for the base scheme is that the Hermite stencil remains the same regardless of its order. This is not the case for FDTD methods. Additionally, the Hermite-Taylor method directly provides a space-time polynomial approximating the solution that is required in the CFM functional. In this paper we focus exclusively on the case when the geometry of the problem can be represented on a Cartesian mesh or on a logically Cartesian curvilinear mesh. Already in this setting the Hermite stencil provides a good advantage but we expect that in future work where we treat interfaces and non-grid aligned boundaries the advantage will be even greater. 

{\color{black} 
We are focusing exclusively on the enforcement of boundary conditions. 
Other important concerns, 
	such as the preservation of the divergence-free constraints and the energy, 
	will not be addressed here.} 
	
The paper is organized as follows. We introduce Maxwell's equations with the considered boundary conditions in Section~\ref{sec:def_problem}. In Section~\ref{sec:HermiteTaylor}, the 1-D Hermite-Taylor method is described in detail and some remarks are provided for higher dimensional cases. The correction function method is introduced and described in detail in the Hermite-Taylor setting in Section~\ref{sec:CFM}. Finally, numerical examples in 1-D and 2-D that verify the properties of the Hermite-Taylor correction function method are presented in Section~\ref{sec:num_examples}. 

\section{Problem Definition} \label{sec:def_problem}
In this work, we seek approximate solutions to Maxwell's equations
\begin{equation} \label{eq:problemMaxwell}
\begin{aligned}
	\mu\,\partial_t \mathbold{H} + \nabla\times \mathbold{E} =&\,\, 0, \\
	\epsilon\,\partial_t \mathbold{E} - \nabla\times\mathbold{H} =&\,\, 0,\\
	\nabla\cdot(\epsilon\,\mathbold{E}) =&\,\, 0,\\
	\nabla\cdot(\mu\,\mathbold{H})=&\,\, 0,
\end{aligned}
\end{equation}
in the domain $\Omega \subset \mathbb{R}^d$ with $d=1,2$ and the time interval $I=[t_0,t_f]$.
Here $\mathbold{H}$ is the magnetic field, 
	$\mathbold{E}$ is the electric field, 
	$\mu$ is the magnetic permeability and $\epsilon$ is the electric permittivity.
To complete the system \eqref{eq:problemMaxwell}, 
	we consider the initial conditions 
	\begin{equation*}
		\begin{aligned}
			\mathbold{H}(\mathbold{x},t_0) =&\,\, \mathbold{H}_0 \quad \text{in } \Omega,\\ 
			\mathbold{E}(\mathbold{x},t_0) =&\,\, \mathbold{E}_0  \quad \text{in } \Omega,
		\end{aligned}
	\end{equation*}
	{\color{black} and the boundary conditions on the electromagnetic fields}.
	
{\color{black} In this work, 
	we focus on the following boundary conditions:
\begin{itemize}
	\item[1.] Perfect electric conductor (PEC):
		\begin{equation} \label{eq:nxE}
			\mathbold{n}\times\mathbold{E} = 0  \quad \text{on } \Gamma \times I,
		\end{equation}
	\item[2.] Perfect magnetic conductor (PMC):
		\begin{equation} \label{eq:nxH}
			\mathbold{n}\times\mathbold{H} = 0 \quad \text{on } \Gamma \times I,
		\end{equation}
	\item[3.] Impedance boundary condition:
		\begin{equation} \label{eq:impedance}
			\mathbold{E}\times\mathbold{n} +Z\,\mathbold{n}\times(\mathbold{H}\times\mathbold{n})= 0 \quad \text{on } \Gamma \times I.
		\end{equation}
\end{itemize}
Here $Z=\sqrt{\tfrac{\mu}{\epsilon}}$ is the impedance, 
	$\Gamma$ is the boundary of the domain $\Omega$ and $\mathbold{n}$ is the outward unit normal to $\Gamma$. 
For further discussions on Maxwell's equations with these boundary conditions and results on their well-posedness, 
	we refer the reader to \cite{Assous2018,Lindell2019}.
Note that we consider the non-homogeneous case of these boundary conditions to facilitate the verification of the Hermite-Taylor correction function method. We denote the given right-hand side function by $\mathbold{g}(x,y,t)$.}
	
\section{Hermite-Taylor Method} \label{sec:HermiteTaylor}
In the following, a brief review of the Hermite-Taylor method, introduced by Goodrich et al. \cite{Goodrich2005}, is provided. For simplicity, we consider the 1-D case and include some comments regarding higher dimensions. 
	
\begin{figure}
 	\centering
	\includegraphics[width=3.5in]{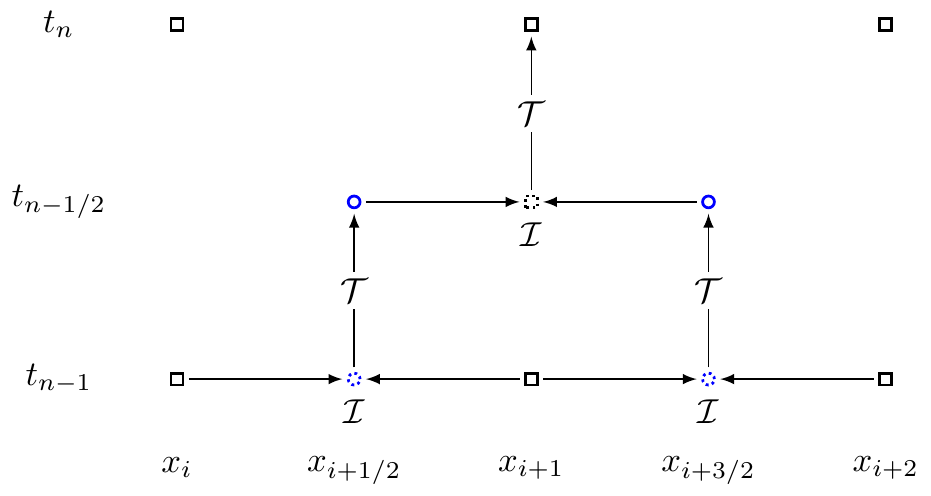}
       \caption{Illustration of the Hermite-Taylor procedure to evolve the data from {\color{black} $(x_{i+1},t_{n-1})$ to $(x_{i+1},t_{n})$}. The Hermite interpolation procedure and the Taylor method are 
       denoted respectively by $\mathcal{I}$ and $\mathcal{T}$. The primal and dual nodes are respectively 
       represented by black squares and blue circles.}
       \label{fig:Hermite_method}
\end{figure}

The Hermite method uses a mesh staggered in both space and time as illustrated in Fig.~\ref{fig:Hermite_method}.
Consider the domain $\Omega = [x_\ell, x_r]$ and a time interval $I = [t_0,t_f]$. 
We then define the primal mesh to be 
	$$x_i = x_\ell + i\,\Delta x, \quad i=0,\dots,N_x, \quad \Delta x = \frac{x_r-x_\ell}{N_x}.$$
Here $N_x$ is the number of cells on the primal mesh.
The dual mesh is then defined as the cell centers of the primal mesh 
	$$x_{i+1/2} = x_\ell + (i+1/2)\,\Delta x, \quad i=0,\dots,N_x-1.$$
The approximate solution on the primal mesh is centered at times 
	$$ t_n = t_0 + n\,\Delta t, \quad n=0,\dots,N_t, \quad \Delta t = \frac{t_f-t_0}{N_t},$$
	while the approximation on the dual mesh is centered at times     
	$$t_{n+1/2} = t_0 + (n+1/2)\,\Delta t, \quad n=0,\dots,N_t-1.$$
Here $N_t$ is the number of time steps. 

The Hermite-Taylor method requires three processes:
\begin{itemize}
	\item[1.] \underline{Hermite interpolation:}

	Assume that the values of the electromagnetic fields and {\color{black} their derivatives through order $m$} 
		(or sufficiently accurate approximation of these) are available on the primal mesh at {\color{black} $t_{n-1}$}. 
	Then, 
		for each cell in the primal mesh, 
		for each electromagnetic field,
		we construct the unique polynomial of degree $2\,m+1$ coinciding with the electromagnetic field and 
		its {\color{black} derivatives through order $m$} at the endpoints of the cell,
		that is the Hermite interpolant of the electromagnetic field. 
	In Fig.~\ref{fig:Hermite_method}, 
		this step is represented by $\mathcal{I}$.
	
	\item[2.] \underline{Recursion relation:}
	
	The recursion relation constructs a space-time polynomial,
		referred as a Hermite-Taylor polynomial in this work, 
		approximating each electromagnetic field.
	Considering a cell and a given Hermite interpolant of each electromagnetic field on this cell, 
		we identify the derivatives of the electromagnetic field as scaled coefficients of the polynomial at the cell center.
	By expanding,
		in time, 
		each scaled coefficient in a Taylor polynomial and enforcing the PDE at the cell center, 
		we obtain a recursion relation for the coefficients of the Hermite-Taylor polynomials. 
	This step is represented in Fig.~\ref{fig:Hermite_method} by either blue dashed circles or black dashed squares.
	
	\item[3.] \underline{Time evolution:}
	
	Finally, 
		we update the electromagnetic fields and their {\color{black} derivatives through order $m$} at the dual mesh points by simply evaluating the Hermite-Taylor polynomials. 
	This step is represented by $\mathcal{T}$ in Fig.~\ref{fig:Hermite_method}.
		
\end{itemize}
Let us now detail each time step of the method.

\subsection{Hermite Interpolation} 
Assuming that the {\color{black} space derivatives through order $m$} of the electromagnetic fields at the initial time $t_0$ are available on the primal mesh, 
	we compute the $(2\,m+1)$ degree Hermite interpolant $p_{i+1/2}^f(x)$ on each cell $[x_i, x_{i+1}]$ satisfying 
	\begin{equation*}
		\frac{d^\ell p_{i+1/2}^f(x_i,t_0)}{dx^\ell} = \frac{d^\ell f(x_i,t_0)}{dx^\ell}, \quad
		\frac{d^\ell p_{i+1/2}^f(x_{i+1},t_0)}{dx^\ell} = \frac{d^\ell f(x_{i+1},t_0)}{dx^\ell}, \quad \ell=0,\dots,m.
	\end{equation*}
Here $f$ is either the magnetic field $H$ or the electric field $E$. 
We then obtain a polynomial approximating each electromagnetic field on the cell $[x_i,x_{i+1}]$ and centered at the cell center $x_{i+1/2}$,
\begin{equation*}
	\begin{aligned}
        H(x,t)|_{t=t_0} \approx&\,\, p^H_{i+1/2}(x) = \sum_{\ell=0}^{2\,m+1} c^H_\ell(t)|_{t=t_0}\,\bigg(\frac{x-x_{i+1/2}}{\Delta x}\bigg)^\ell, \\
        E(x,t)|_{t=t_0} \approx&\,\, p^E_{i+1/2}(x) = \sum_{\ell=0}^{2\,m+1} c^E_\ell(t)|_{t=t_0}\,\bigg(\frac{x-x_{i+1/2}}{\Delta x}\bigg)^\ell,
       	\end{aligned}
\end{equation*}
	where $c^H(t)$ and $c^E(t)$ are time-dependent coefficients.

\subsection{Recursion Relation} 

Let us now compute a Hermite-Taylor polynomial approximating each electromagnetic field. 	
To do so, 
	we expand the coefficients in a Taylor polynomial of degree $q$ centered at $t_0$, 
	which leads to 
\begin{equation} \label{eq:Hermite_Taylor_polynomials}
	\begin{aligned}
        H(x,t) \approx&\,\, p^H_{i+1/2}(x,t) = \sum_{\ell=0}^{2\,m+1} \sum_{s=0}^q c^H_{\ell,s}\,\bigg(\frac{x-x_{i+1/2}}{\Delta x}\bigg)^\ell\,\bigg(\frac{t-t_{0}}{\Delta t}\bigg)^s, \\
        E(x,t) \approx&\,\, p^E_{i+1/2}(x,t) = \sum_{\ell=0}^{2\,m+1} \sum_{s=0}^q c^E_{\ell,s}\,\bigg(\frac{x-x_{i+1/2}}{\Delta x}\bigg)^\ell\,\bigg(\frac{t-t_{0}}{\Delta t}\bigg)^s.
       	\end{aligned}
\end{equation}
Here $c^H_{\ell,0}$ and $c^E_{\ell,0}$ are known from the initial data and the interpolation step. 
Consider Maxwell's equations in 1-D with constant coefficients, 
	\begin{equation*}
		\begin{aligned}
    			\frac{\partial H}{\partial t} =&\,\, -\frac{1}{\mu}\frac{\partial E}{\partial x}, \\
    			\frac{\partial E}{\partial t}  =&\,\, -\frac{1}{\epsilon}\frac{\partial H}{\partial x}.
		\end{aligned}
	\end{equation*}
For smooth solutions, 
	we then have 
	\begin{equation} \label{eq:recursion_Maxwell_1D}
		\begin{aligned} 
		\frac{\partial^{\ell+s+1}H}{\partial t^{s+1}\partial x^\ell}=&\,\, -\frac{1}{\mu}\frac{\partial^{\ell+s+1}E}{\partial t^{s}\partial x^{\ell+1}},\\
		\frac{\partial^{\ell+s+1}E}{\partial t^{s+1}\partial x^\ell}=&\,\, -\frac{1}{\epsilon}\frac{\partial^{\ell+s+1}H}{\partial t^{s}\partial x^{\ell+1}}.
		\end{aligned}
	\end{equation}
Substituting $H$ and $E$ by their Hermite-Taylor approximations $p_{i+1/2}^H(x,t)$ and $p_{i+1/2}^E(x,t)$, 
	in the system \eqref{eq:recursion_Maxwell_1D} and evaluating them at $(x_{i+1/2},t_0)$,
	we obtain the following recursion relations for the coefficients
\begin{equation*} 
	c^H_{\ell,s} = -\frac{(\ell+1)\,\Delta t}{\mu\,s\,\Delta x}\,c^E_{\ell+1,s-1}, \quad
	c^E_{\ell,s} = -\frac{(\ell+1)\,\Delta t}{\epsilon\,s\,\Delta x}\,c^H_{\ell+1,s-1}, \quad
	\ell = 0,\dots,2\,m+1, \quad s=1,\dots,q.
\end{equation*}
Knowing $c^H_{\ell,0}$ and $c^E_{\ell,0}$, 
	these recursion relations allow the computation of the Hermite-Taylor polynomials approximating the electromagnetic fields.

{\color{black}
\subsubsection{Variable Coefficients Problems}

For spatially variable coefficients, 
	the recursion relations for the polynomial coefficients involve high-order derivatives of the coefficients.
As in \cite{Goodrich2005}, 
	we consider high-order derivatives of the coefficients $u(x) = \frac{1}{\mu(x)}$ and 
	$e(x) = \frac{1}{\epsilon(x)}$ and take advantage of the Leibnitz' rule.  
	
For sufficiently smooth solutions and coefficients $u$ and $e$, 
	we have 
\begin{equation} \label{eq:recursion_Maxwell_1D_var_coeff}
	\begin{aligned} 
		\frac{\partial^{\ell+s+1}H}{\partial t^{s+1}\partial x^\ell}=&\,\, -\frac{\partial^\ell}{\partial x^\ell} \bigg(u\,\frac{\partial^{s+1}E}{\partial t^{s}\partial x}\bigg) = \sum_{i=0}^\ell \binom{\ell}{i}\frac{\partial^{\ell-i}u}{\partial x^{\ell-i}}\frac{\partial^{i+s+1}E}{\partial t^s\partial x^{i+1}} ,\\
		\frac{\partial^{\ell+s+1}E}{\partial t^{s+1}\partial x^\ell}=&\,\, -\frac{\partial^\ell}{\partial x^\ell} \bigg(e\,\frac{\partial^{s+1}H}{\partial t^{s}\partial x}\bigg) = \sum_{i=0}^\ell \binom{\ell}{i}\frac{\partial^{\ell-i}e}{\partial x^{\ell-i}}\frac{\partial^{i+s+1}H}{\partial t^s\partial x^{i+1}}.
	\end{aligned}
\end{equation}
Identifying the coefficients of the Hermite-Taylor polynomials as scaled derivatives and enforcing the system \eqref{eq:recursion_Maxwell_1D_var_coeff} at $(x_{i+1/2},t_0)$,
	we obtain
\begin{equation*} 
	\begin{aligned}
	c^H_{\ell,s} =&\,\, -\sum_{i=0}^{\ell}  \frac{(i+1)\,\Delta t\, \Delta x^{\ell-i-1}}{(\ell-i)!\,s}\,\frac{\partial^{\ell-i} u}{\partial x^{\ell-i}}\,c^E_{i+1,s-1}, \\
	c^E_{\ell,s} =&\,\, -\sum_{i=0}^{\ell}  \frac{(i+1)\,\Delta t\, \Delta x^{\ell-i-1}}{(\ell-i)!\,s}\,\frac{\partial^{\ell-i} e}{\partial x^{\ell-i}}\,c^H_{i+1,s-1}, 
	\end{aligned}
\end{equation*}
	for $\ell = 0,\dots,2\,m+1$ and  $s=1,\dots,q$.
}
	
{\color{black} 
We note that this step can be generalized for other problems including linear, 
	non-linear and variable coefficient problems. 
We refer the reader to \cite{hagstrom2015solving} for more details.}

\subsection{Time Evolution}
Finally, 
	we evolve	the electromagnetic fields and their {\color{black} space derivatives through order $m$} on the dual mesh nodes,
	located at $(x_{i+1/2},t_{1/2})$ for the cell $[x_i, x_{i+1}]$, 
	by evaluating \eqref{eq:Hermite_Taylor_polynomials}
$$\frac{\partial^{\ell} p^H_{i+1/2}(x_{i+1/2},t_{1/2})}{\partial x^\ell}, \quad \frac{\partial^{\ell} p^E_{i+1/2}(x_{i+1/2},t_{1/2})}{\partial x^\ell}, \quad \ell=0,\dots,m.$$
	
A similar process is repeated to evolve the data from the dual mesh at $t_{1/2}$ to the primal mesh at $t_1$ and therefore to complete the time step. 
The overall procedure is repeated until the final time is reached. 
Fig.~\ref{fig:Hermite_method} illustrates the Hermite-Taylor method at a given primal node.
	 
\begin{remark}
For linear constant coefficients hyperbolic problems, 
	the Taylor expansion in time of the coefficients of the Hermite polynomials is computed exactly for $q$ sufficiently large \cite{Goodrich2005},
	for example $q=2\,m+1$ in \eqref{eq:Hermite_Taylor_polynomials} for the 1-D case.
In general, 
	 we set $q=\nu\,(2\,m+1)$ in $\mathbb{R}^{\nu}$ to obtain an exact time expansion of the coefficients.
\end{remark}
		
\begin{remark}
In higher dimensions, 
	the primal mesh is defined as the classical Cartesian mesh while the dual nodes are defined 
	at the cell center. 
Hence, 
	this differs from the mesh used in FDTD methods. 
As for the Hermite interpolation procedure, 
	approximations are computed using a tensor product of 1-D Hermite polynomials.
We refer the interested reader to \cite{Goodrich2005} for more details 
	on the Hermite-Taylor setting for higher dimensions. 
\end{remark}
	
As mentioned before,
	a challenge for the Hermite-Taylor method is to enforce general boundary conditions. 
Indeed, 
	this method requires to know all information on the boundary,
	including the {\color{black} space derivatives through order $m$}, 
	which are usually not available. 
In the next section,
	we present a way to obtain the needed information using the correction function method.

\section{Correction Function Method} \label{sec:CFM}

In this section, we describe the correction function method that computes approximations to the electromagnetic fields and their {\color{black} derivatives through order $m$} at the nodes located on the boundary of the domain. There are two key ingredients to the CFM: the minimization of functionals describing the electromagnetic fields near the boundary, and careful definition of the space-time domains of the functionals along the boundary. We refer to a space-time domain of a functional as a local patch. Once the minimization procedure is completed, we obtain space-time polynomials, called correction functions, approximating each electromagnetic field in the vicinity of the boundary. The correction functions are used to update the solution at the boundary nodes. In the following, we first describe the method in detail in 1-D and then generalize it in higher dimensions.

\subsection{The Hermite CFM Method in One Dimension} \label{sec:cfm_1D}
On the mesh in Fig.~\ref{fig:Hermite_method}, {\color{black} the first step has allowed for the update of the Hermite solution on the dual mesh at time level $t_{n-1/2}$ and the second step has allowed for the update of the numerical solution on the primal mesh at $t_{n}$, 	except near the boundary. At $(x_0,t_n)$ and $(x_{N_x},t_n)$ for $n=1,\dots,N_t$  the solution will be updated using the CFM.} 
	
We define a node where the numerical solution is updated using the Hermite-Taylor method as a Hermite node and a node where the numerical solution is computed using the CFM we denote as a CF node. In the following, the subscript $i$ refers to the $i^{\text{th}}$ CF node in the mesh and the superscript $n$ refers to the time level $t_n$. In the 1-D case, $i=0$ and $i=1$ refer respectively to the boundary nodes $x_0$ and $x_{N_x}$. 

We further note that although the functional just to be defined can depend on time, as manifested by the $n$ superscript, (for example to account for a moving geometry) but for all the problems considered here it will not. When there is no time dependence all the small linear system of equations (one at each CF node) resulting from the quadratic optimization problem, will not change in time and can thus be formulated, factored and stored once and for all before the time stepping loop. Consequently the complexity of the Hermite-CFM method will approach that of the Hermite method in the limit $h \rightarrow 0$.    

The CFM minimizes a functional unique to each CF node composed of three parts 
\begin{equation} \label{eq:functional_J}
	J_i^n = \mathcal{G}_i^n + \mathcal{B}_i^n + \mathcal{H}_i^n.
\end{equation}
Here, $\mathcal{G}_i^n$ weakly enforces the governing equations, $\mathcal{B}_i^n$ weakly enforces the boundary conditions and $\mathcal{H}_i^n$ weakly enforces that the correction functions match the Hermite solution near the $i^\text{th}$ CF node.

The domains over which the different terms in the functional are computed are not the same.
The domain of $\mathcal{B}_i^n$ should include the part of the boundary in the vicinity of the $i^{\text{th}}$ CF node to weakly enforce the boundary conditions. The domain of $\mathcal{H}_i^n$ should be the same as the space-time domains of the Hermite nodes closest to the $i^{\text{th}}$ CF node. We then weakly enforce the correction functions to match the Hermite solution in the domain of $\mathcal{H}_i^n$ while avoiding extrapolation procedures of the Hermite solution. 
Finally, the domain of integration for $\mathcal{G}^n_i$ should enclose the $i^{\text{th}}$ CF node, the domain of integration for $\mathcal{B}_i^n$ and the domain of integration for $\mathcal{H}_i^n$ to enforce Maxwell's equations over the whole local patch of the functional $J_i^n$.

\begin{figure}
 	\centering
	\includegraphics[width=3.0in]{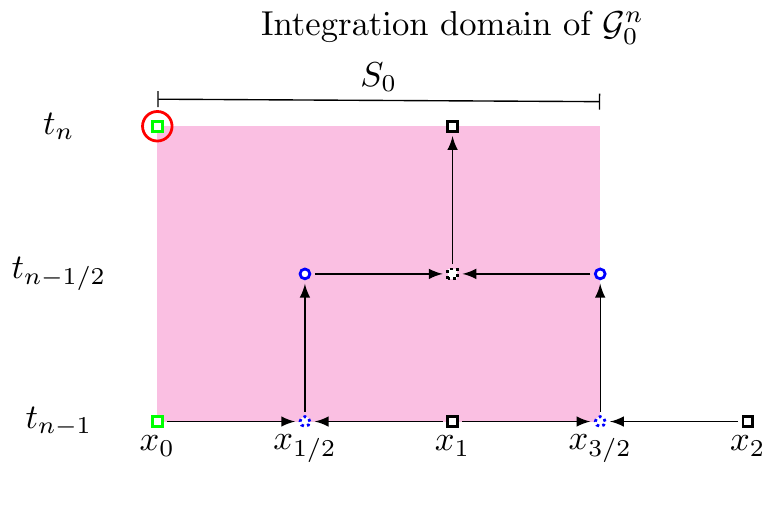}
	\caption{Illustration of the domain of integration $S_0\times I_n$ of $\mathcal{G}_0^n$. 
		   The primal CF and Hermite nodes are respectively represented by green squares and black squares 
		   while the dual Hermite nodes are represented by blue circles.
		   The CFM seeks the information located at $(x_0,t_{n})$ which is enclosed by the red circle.
		   The space-time local patch $S_0\times I_n$ is denoted by a dashed magenta box.}
	\label{fig:local_patch_1D_x0}
\end{figure}

As an example for the CF node $x_0$ at time level $t_n$, $\mathcal{G}_0^n$ contains the residual of the PDE and it is integrated over the rectangular space-time region (the local patch) consisting of the direct product of the space interval $S_0=[x_0,x_{3/2}]$ with the time interval {\color{black} $I_n=[t_{n-1},t_n]$} as illustrated in Fig.~\ref{fig:local_patch_1D_x0}. 
We then have 
\begin{equation*}
	\mathcal{G}_0^n(H^n_{h,0},E^n_{h,0}) =  \frac{\ell_0}{2} \,\int\limits_{I_n}\!\int\limits_{S_0} (\mu\,\partial_t H^n_{h,0}+\partial_x E^n_{h,0})^2 + (\epsilon\,\partial_t E^n_{h,0} + \partial_x H^n_{h,0})^2\,\mathrm{d}x\,\mathrm{d}t,
\end{equation*}
where $\ell_0 = x_{3/2}-x_0 = 1.5\,\Delta x$ is the characteristic length of the space interval $S_0$. Here $H^n_{h,0}$ and $E^n_{h,0}$ are the sought correction functions approximating the electromagnetic fields and are used to update the numerical solution at $(x_0,t_n)$.

\begin{figure}
 	\centering
	\includegraphics[width=3.0in]{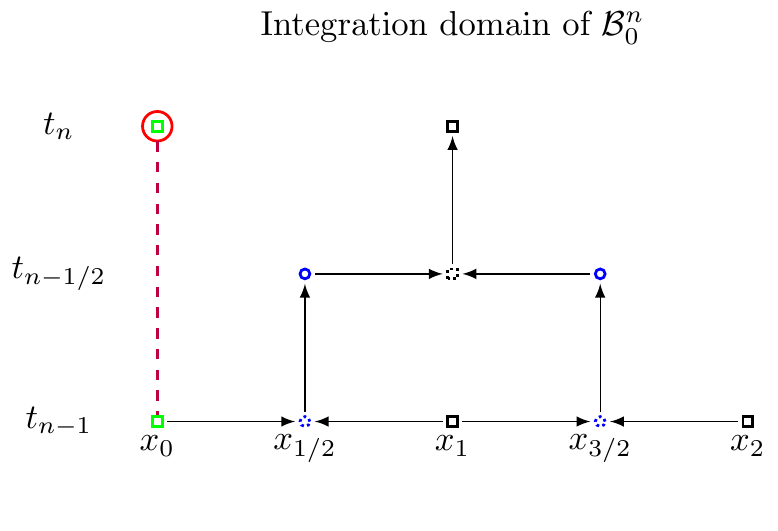}
	\caption{Illustration of the domain of integration $I_n$ at $x_0$ of $\mathcal{B}_0^n$. 
		   The primal CF and Hermite nodes are respectively represented by green squares and black squares 
		   while the dual Hermite nodes are represented by blue circles.
		   The CFM seeks the information located at $(x_0,t_{n})$ which is enclosed by the red circle.
		   The intersection between the boundary and the local patch,
			that is the line connecting $(x_0,t_{n-1})$ to $(x_0,t_n)$,  
			is denoted by a dashed purple line.}	
	\label{fig:boundary_1D_x0}
\end{figure}

\begin{figure}
 	\centering
	\includegraphics[width=3.0in]{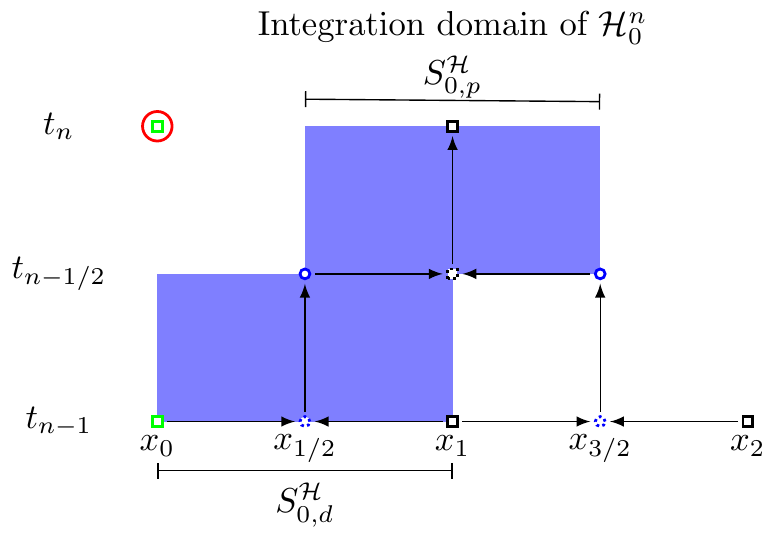}
       \caption{Illustration of {\color{black} the domains of integration $S_{0,d}^\mathcal{H} \times [t_{n-1},t_{n-1/2}]$ and $S_{0,p}^\mathcal{H} \times [t_{n-1/2},t_n]$ of $\mathcal{H}_0^n$.} 
		   The primal CF and Hermite nodes are respectively represented by green squares and black squares 
		   while the dual Hermite nodes are represented by blue circles.
		   The CFM seeks the information located at $(x_0,t_{n})$ which is enclosed by the red circle.
		   {\color{black} The domains $S_{0,d}^\mathcal{H} \times [t_{n-1},t_{n-1/2}]$ and $S_{0,p}^\mathcal{H} \times [t_{n-1/2},t_n]$},
			where we enforce the correction functions to match the Hermite-Taylor polynomials,
			is denoted by a dashed blue 
		   box.}
	\label{fig:HermiteTaylor_patch_1D_x0}
\end{figure}

The term $\mathcal{B}_0^n$ contains the residual of the boundary condition at $x_0$ and it is integrated over the time interval $I_n$ as shown in Fig.~\ref{fig:boundary_1D_x0}. 
{\color{black} As an example, we have  
\begin{equation*}
	\mathcal{B}_0^n(E^n_{h,0}) = \frac{1}{2} \, \int\limits_{I_n} (E^n_{h,0}(x_0,t)-g_E(t))^2\,\mathrm{d}t,
\end{equation*}
for the boundary condition \eqref{eq:nxE}.} 
We now require the correction functions to weakly match the Hermite solution over the  {\color{black} space-time domains of the primal Hermite node $x_1$ and the dual Hermite node $x_{1/2}$}. This is what connects the two methods and is needed for the minimization problem to be well-posed. 
{\color{black}
The first part of $\mathcal{H}_0^n$ contains the Hermite-Taylor polynomials $H^*(x,t) = p_{1/2}^H(x,t)$ and $E^*(x,t) = p_{1/2}^E(x,t)$, 
	which are associated with the cell of the dual Hermite node $x_{1/2}$,
	and it is integrated over the rectangular region consisting of the direct product of the space interval $S^{\mathcal{H}}_{0,d} = [x_{0},x_{1}]$ with the time interval $[t_{n-1},t_{n-1/2}]$.
The second part of the term $\mathcal{H}_0^n$ contains the Hermite-Taylor polynomials $H^*(x,t)=p_{1}^H(x,t)$ and $E^*(x,t)=p_{1}^E(x,t)$, 
	and it is integrated over the rectangular space-time region consisting of the direct product of the space interval $S^{\mathcal{H}}_{0,p} = [x_{1/2},x_{3/2}]$, 
	which is the cell associated with the primal Hermite node $x_1$, 
	with the time interval $[t_{n-1/2},t_n]$. 
The space-time regions $S^{\mathcal{H}}_{0,d}\times [t_{n-1},t_{n-1/2}]$ and $S^{\mathcal{H}}_{0,p}\times [t_{n-1/2},t_{n}]$ are illustrated in Fig.~\ref{fig:HermiteTaylor_patch_1D_x0}. We then have 
\begin{equation} \label{eq:infoHermite}
	\begin{aligned}
	\mathcal{H}_0^n(H^n_{h,0},E^n_{h,0}) =&\,\,  \frac{c_H}{2}\,\int\limits_{t_{n-1}}^{t_{n-1/2}} \int\limits_{S_{0,d}^\mathcal{H}} (H^n_{h,0}-H^*)^2 + (E^n_{h,0}-E^*)^2\,\mathrm{d}x\,\mathrm{d}t\\
	+&\,\, \frac{c_H}{2}\,\int\limits_{t_{n-1/2}}^{t_n} \int\limits_{S_{0,p}^\mathcal{H}} (H^n_{h,0}-H^*)^2 + (E^n_{h,0}-E^*)^2\,\mathrm{d}x\,\mathrm{d}t,
	\end{aligned}
\end{equation}
where $c_H$ is a given penalization function that is such that $0< c_H(\Delta x)\leq 1$.}

A similar procedure is used to define the local patch and the functional associated with the second CF node $x_{N_x}$ at the time level $t_n$.

\subsubsection{The Linear System of Equations that Solves the Optimization Problem}
At each CF node we must solve the following problem. 
	\begin{equation} \label{eq:minPblm_1D}
		\begin{aligned}
			&\text{Find } ({H}^n_{h,i},{E}^n_{h,i}) \in V \times V \text{ such that }\\
 			&\qquad ({H}^n_{h,i},{E}^n_{h,i}) =  \underset{{v},{w}\in V}{\arg\min}\, J_i^n({v},{w}).
		\end{aligned}
	\end{equation}
Here $V = \mathbb{Q}^k\big(S_i\times I_n \big)$ is the space of polynomials of degree $k$. In this work, we use space-time Legendre polynomials. In our one dimensional example $i=0,1$. Note that although $n=1,\dots,N_t$, since the boundary does not change in time, there is in fact only one optimization problem for each CF node.
 
We formally compute the gradient of $J_i^n$ with respect to the coefficients of the polynomial approximations $H^n_{h,i}$ and $E^n_{h,i}$, and use that it vanishes at a minimum to find a solution to the minimization problem \eqref{eq:minPblm_1D}.
This leads to a linear system 
\begin{equation*}
	M_i^n\, \mathbold{c}_i^n = \mathbold{b}_i^n,
\end{equation*}
where $\mathbold{c}_i^n$ contains the coefficients of $H^n_{h,i}$ and $E^n_{h,i}$. 

Again, since the boundary of the domain does not move, we have $M_i = M_i^n$, so the matrices $M_i$, their scaling and LU factorization are found in a pre-computation step. Consequently, the only computations needed at each time step is the computation of the right-hand side $\mathbold{b}_i^n$, followed by forward and backward substitutions to find $\mathbold{c}_i^n$.

\subsubsection{Summary of the Hermite-CFM Method in One Dimension}
Given the numerical solution on the primal mesh at $t_{n-1}$, the algorithm of the Hermite-Taylor correction function method to evolve the numerical solution at $t_{n}$ is:
\begin{itemize}
\item[1.] Update the numerical solution on the dual mesh at $t_{n-1/2}$ using the Hermite-Taylor method {\color{black} and store the Hermite-Taylor polynomials needed for the CFM};
\item[2.] Update the numerical solution on the primal Hermite node at $t_{n}$ using the Hermite-Taylor method and store the Hermite-Taylor polynomials needed for the CFM;
\item[3.] Update the numerical solution at the CF nodes using the CFM by computing the right hand sides $\mathbold{b}_i^n$ and solve for $\mathbold{c}_i^n$. This is done independently for each $i$ and can thus be done in parallel without any communication step.
\end{itemize}

\subsection{The Hermite-CFM Method in Two Dimensions}
\begin{figure} 
	\centering
	\includegraphics[width=2.25in]{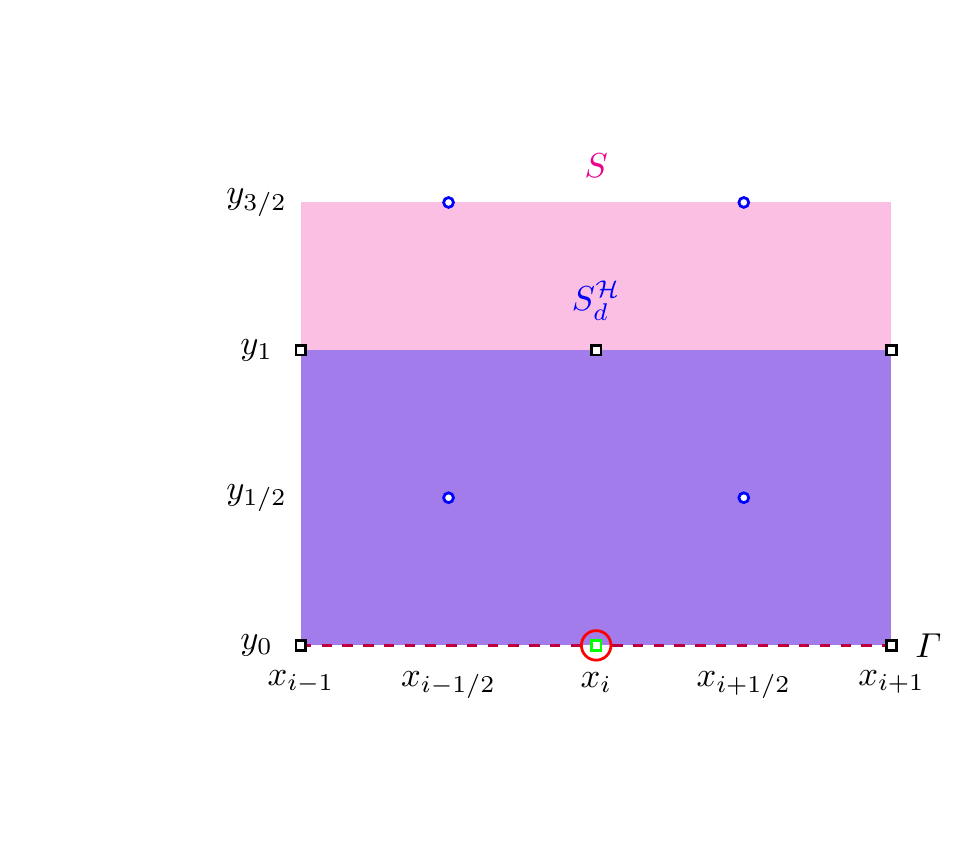} 
	\includegraphics[width=2.25in]{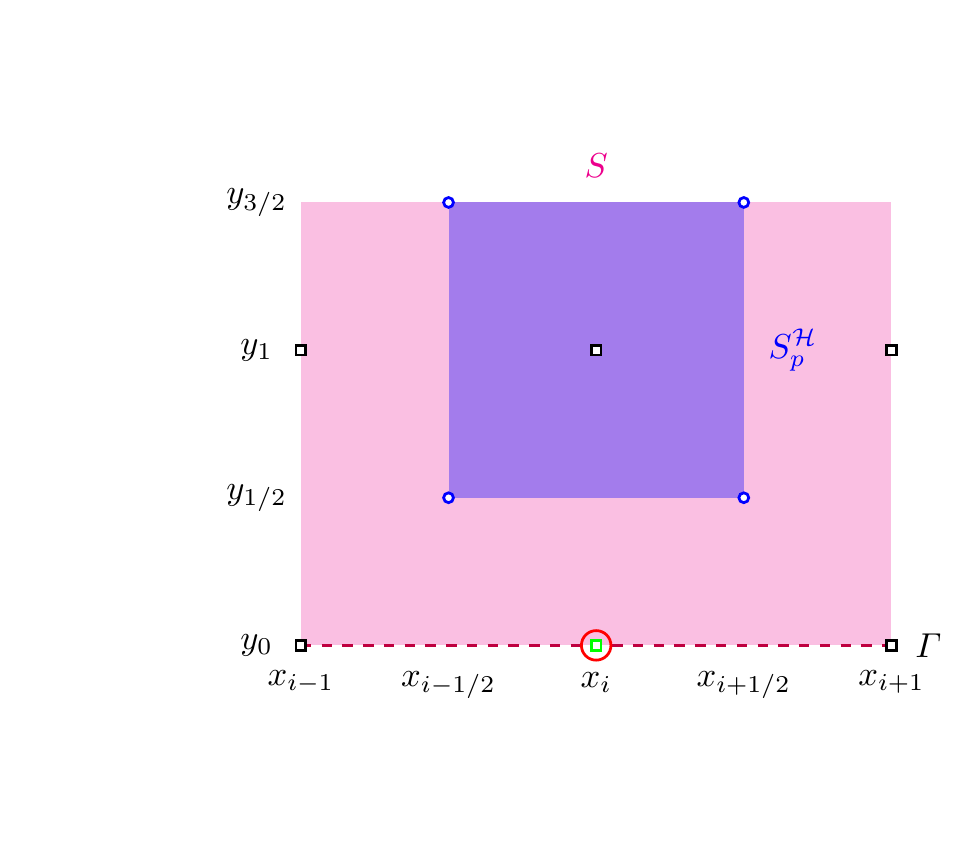}
       \caption{Illustration of a 2-D local patch for a bottom edge CF node. 
       		   {\color{black} The left and right plots show respectively the spatial component of the local patch over the time intervals $[t_{n-1},t_{n-1/2}]$ and $[t_{n-1/2},t_n]$.
       		   The primal CF and Hermite nodes are respectively represented by green squares and black squares 
		   while the dual Hermite nodes are represented by blue circles.
		   The CFM seeks the information located at $(x_i,y_0)$,
			which is enclosed by the red circle.
		   The spatial domain $S$ of local patches is denoted by a dashed magenta box.
		   The part of the boundary $\Gamma$ include in the local patch is represented by a dashed purple line.
		   The spatial domains $S_d^{\mathcal{H}}$ and $S_p^{\mathcal{H}}$ where we enforce the correction functions to match the Hermite-Taylor polynomials
		   	are denoted by a dashed blue box.}
		   }
	\label{fig:edgePatch}
\end{figure}

\begin{figure} 
 	\centering
	\includegraphics[width=2.25in]{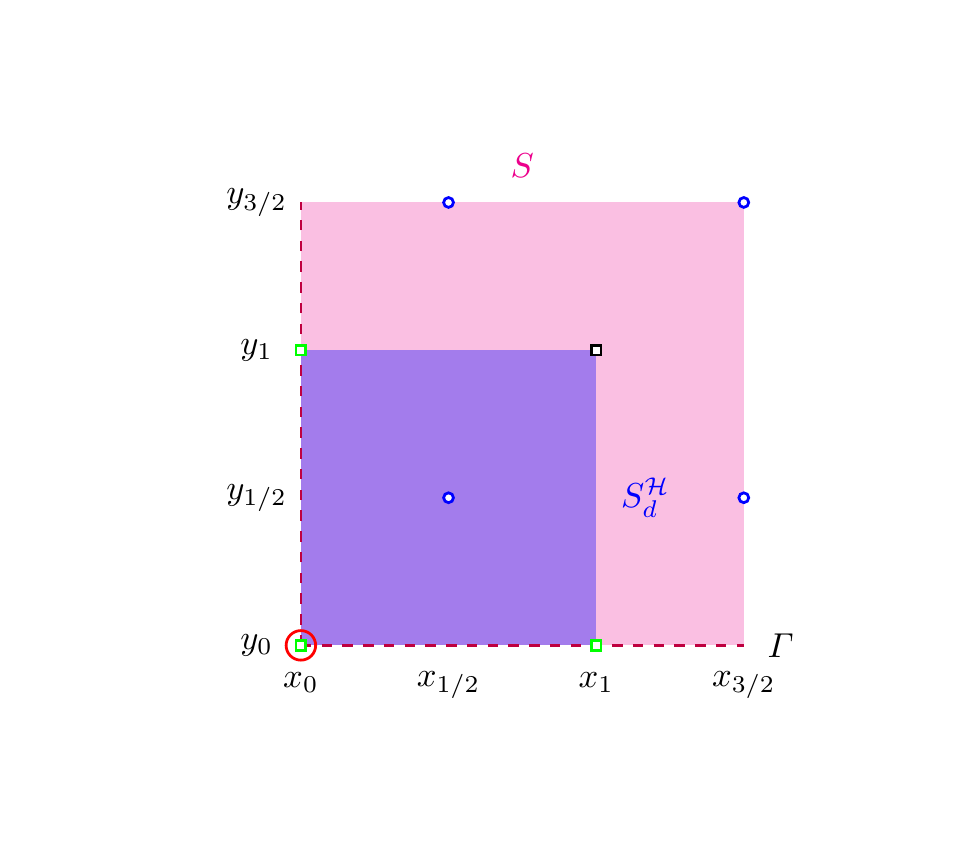}
	\includegraphics[width=2.25in]{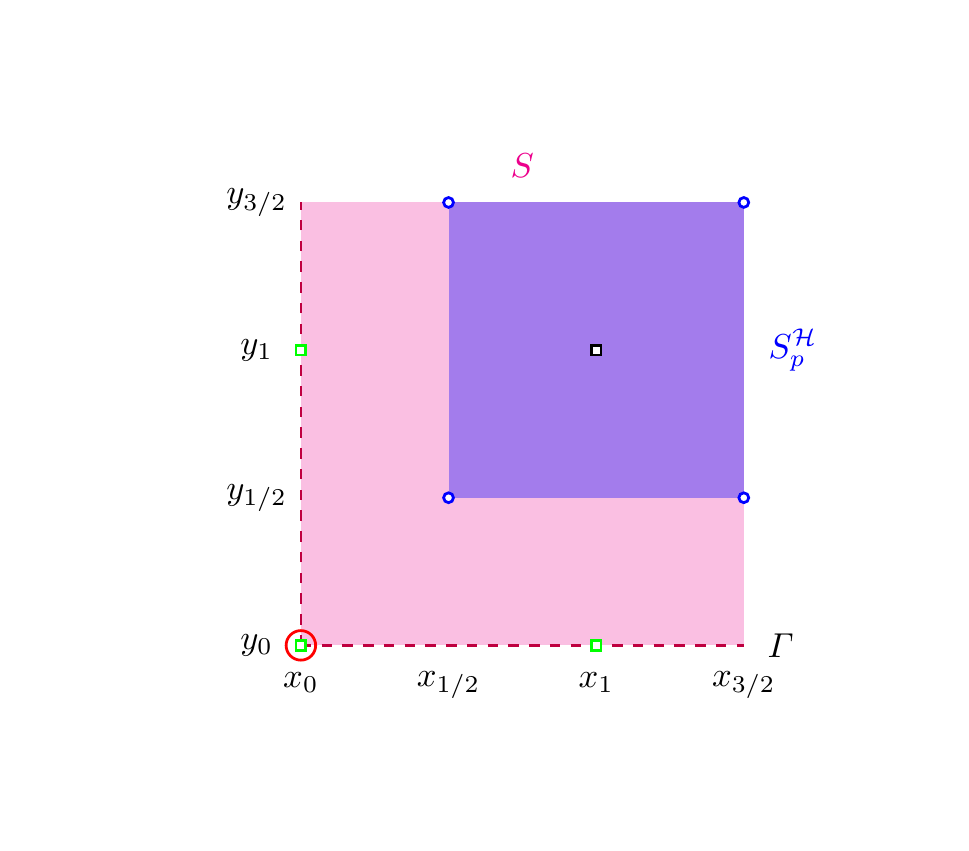}
	\caption{Illustration of a 2-D local patch for a bottom-left corner CF node. 
       		   {\color{black} The left and right plots show respectively the spatial component of the local patch over the time intervals $[t_{n-1},t_{n-1/2}]$ and $[t_{n-1/2},t_n]$.
		   The primal CF and Hermite nodes are respectively represented by green squares and black squares 
		   while the dual Hermite nodes are represented by blue circles.
		   The CFM seeks the information located at $(x_0,y_0)$,
			which is enclosed by the red circle.
		   The spatial domain $S$ of local patches is denoted by a dashed magenta box.
		   The part of the boundary $\Gamma$ include in the local patch is represented by a dashed purple line.
		   The spatial domains $S_d^{\mathcal{H}}$ and $S_p^{\mathcal{H}}$ where we enforce the correction functions to match the Hermite-Taylor polynomials
		   	are denoted by a dashed blue box.}
		   }
		   \label{fig:cornerPatch}
\end{figure}

\begin{figure} 
 	\centering
	\includegraphics[width=2.25in]{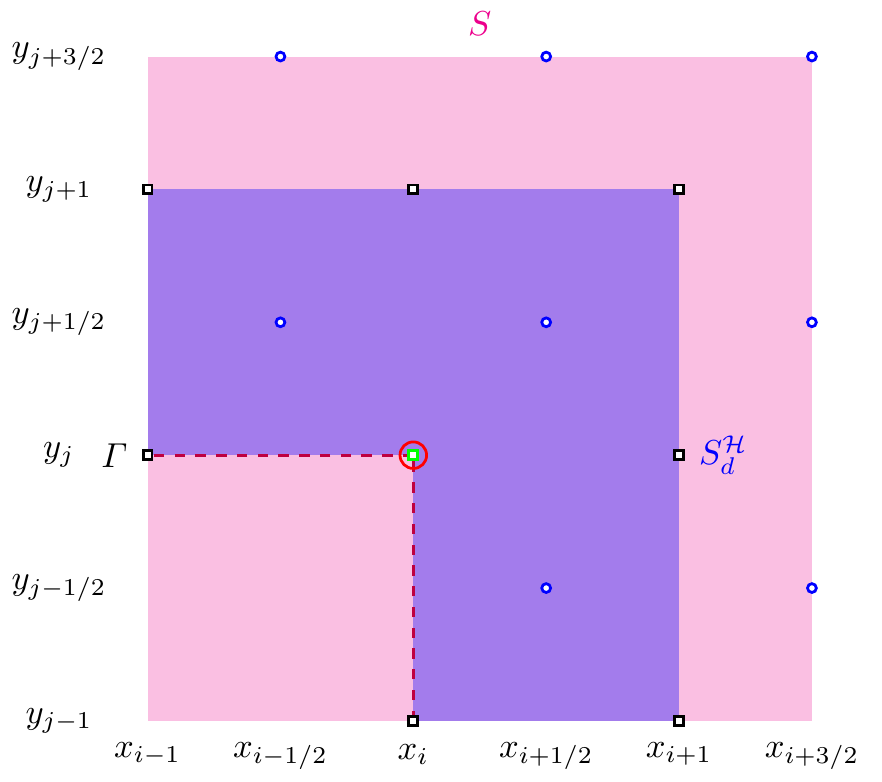} \hspace*{0.05in}
	\includegraphics[width=2.25in]{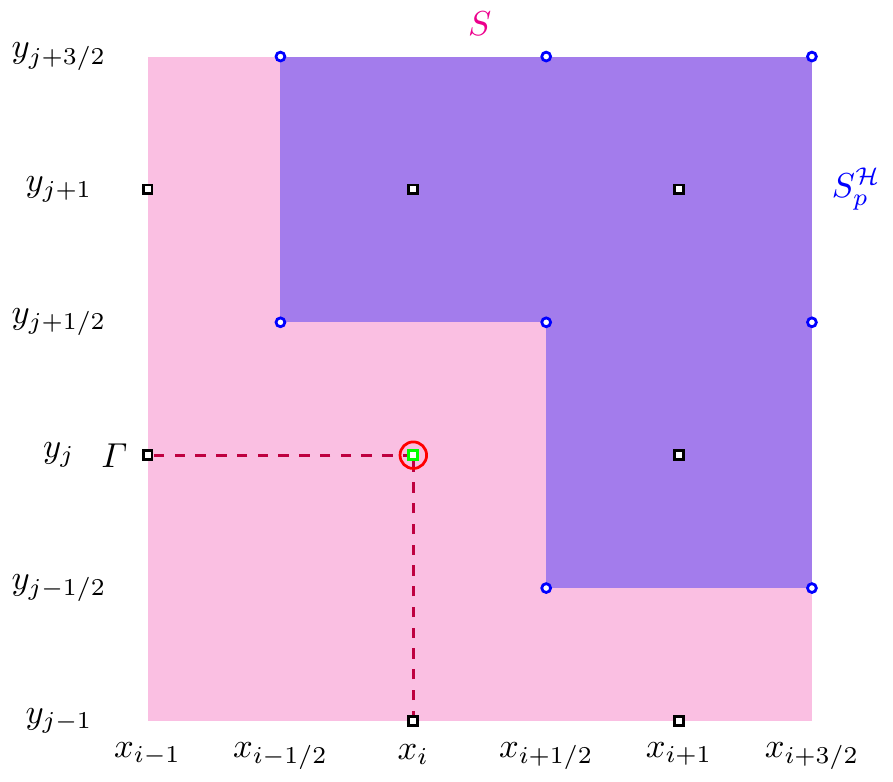}
       \caption{Illustration of a 2-D local patch for a reentrant corner CF node. 
       	          {\color{black} The left and right plots show respectively the spatial component of the local patch over the time intervals $[t_{n-1},t_{n-1/2}]$ and $[t_{n-1/2},t_n]$.
       		   The primal CF and Hermite nodes are respectively represented by green squares and black squares 
		   while the dual Hermite nodes are represented by blue circles.
		   The CFM seeks the information located at $(x_i,y_i)$,
			which is enclosed by the red circle.
		   The spatial domain $S$ of local patches is denoted by a dashed magenta box.
		   The part of the boundary $\Gamma$ include in the local patch is represented by a dashed purple line.
		   The spatial domains $S_d^{\mathcal{H}}$ and $S_p^{\mathcal{H}}$ where we enforce the correction functions to match the Hermite-Taylor polynomials
		   	are denoted by a dashed blue box.}
		   }
		   \label{fig:reentrant_corner_patch}
\end{figure}

We only consider piecewise rectangular domains composed of straight lines between primal nodes. For higher dimensions, the spatial domain of a local patch is adapted depending on the geometry of the boundary and where the Hermite solution is available 
	in the vicinity of its CF node while the time domain $I_n$ remains the same. 
The spatial domain $S_i$ of a local patch needs to satisfy three constraints:
	\begin{itemize}
	\item[1.] The $i^{\text{th}}$ CF node must be inside;
	\item[2.] Part of the boundary of the domain close to the $i^{\text{th}}$ CF node must be contained in it;
	\item[3.] It must contain the cells of the Hermite nodes closest to the $i^{\text{th}}$ CF node.
	\end{itemize}
Examples of the spatial domains of local patches in 2-D that satisfy these constraints are shown in Fig.~\ref{fig:edgePatch}, 
	Fig.~\ref{fig:cornerPatch} and Fig.~\ref{fig:reentrant_corner_patch}.
For simplicity, 
	we omit the subscript associated with the CF node in the description of the local patches.
	
Let us first consider a CF node $(x_i,y_0)$ along an edge as depicted in Fig.~\ref{fig:edgePatch}.
{\color{black} In this case,
	the spatial domain of the local patch is $S = [x_{i-1},x_{i+1}]\times[y_0,y_{3/2}]$ while its intersection with 
	the boundary of the domain, 
	$S\cap\Gamma$, 
	is the line connecting the points $(x_{i-1},y_0)$ and $(x_{i+1},y_0)$.
The spatial domains where we weakly enforce the Hermite solution are $S_d^{\mathcal{H}} = [x_{i-1},x_{i+1}]\times[y_{0},y_{1}]$ over the time interval $[t_{n-1},t_{n-1/2}]$
	and $S_p^{\mathcal{H}} = [x_{i-1/2},x_{i+1/2}]\times[y_{1/2},y_{3/2}]$ over the time interval $[t_{n-1/2},t_n]$.}

{\color{black} For a CF node located at a corner $(x_0,y_0)$ as illustrated in Fig.~\ref{fig:cornerPatch}, 
	we have $S = [x_0,x_{3/2}]\times[y_0,y_{3/2}]$,
	$S_d^{\mathcal{H}} = [x_{0},x_{1}]\times[y_{0},y_{1}]$ and $S_p^{\mathcal{H}} = [x_{1/2},x_{3/2}]\times[y_{1/2},y_{3/2}]$.
The intersection of $S$ with the boundary is composed of the line connecting $(x_0,y_0)$ to $(x_0,y_{3/2})$ and that connecting $(x_0,y_0)$ to $(x_{3/2},y_0)$.}

As a last example, 
	we consider the situation in Fig.~\ref{fig:reentrant_corner_patch} where a CF node is located at a reentrant corner $(x_i,y_i)$.
{\color{black} We then have $S=[x_{i-1},x_{i+3/2}]\times[y_{j-1},y_{j+3/2}]$. 
The spatial domain where the Hermite solution is enforced $S_d^\mathcal{H}$ over the time interval $[t_{n-1},t_{n-1/2}]$ is the union of $[x_{i-1},x_{i+1}]\times[y_{j},y_{j+1}]$
 	and $[x_{i},x_{i+1}]\times[y_{j-1},y_{j}]$.
The spatial domain where the Hermite solution is enforced $S_p^\mathcal{H}$ over the time interval $[t_{n-1/2},t_n]$ is the union of $[x_{i-1/2},x_{i+3/2}]\times[y_{j+1/2},y_{j+3/2}]$
 	and $[x_{i+1/2},x_{i+3/2}]\times[y_{j-1/2},y_{j+1/2}]$.
The intersection between the spatial domain $S$ of the local patch and the boundary is composed of the line connecting $(x_i,y_{j-1})$ to $(x_i,y_j)$ and 
	that connecting $(x_{i-1},y_j)$ to $(x_i,y_j)$.}

Let us now consider Maxwell's equations in 3-D and seek polynomial approximations of the magnetic field and the electric field 
	in each local patch, 
	that is $\mathbold{H}^n_{h,i}$ and $\mathbold{E}^n_{h,i}$ for $i=0,\dots,N_\Gamma$ and $n=1,\dots,N_t$.
Here $N_\Gamma$ is the total number of CF nodes.
The first part of the functional \eqref{eq:functional_J} becomes
	\begin{equation*}
		\begin{aligned}
			\mathcal{G}_i^n(\mathbold{H}^n_{h,i},\mathbold{E}^n_{h,i}) =&\,\,  \frac{\ell_i}{2} \,\int\limits_{I_n}\!\int\limits_{S_i} (\mu\,\partial_t \mathbold{H}^n_{h,i} + 
			\nabla\times\mathbold{E}^n_{h,i})\cdot(\mu\,\partial_t \mathbold{H}^n_{h,i} + \nabla\times\mathbold{E}^n_{h,i}) \\
			+&\,\, ( \epsilon\,\partial_t \mathbold{E}^n_{h,i} - \nabla\times\mathbold{H}^n_{h,i})\cdot( \epsilon\,\partial_t \mathbold{E}^n_{h,i} - 
			\nabla\times\mathbold{H}^n_{h,i})\\
			+&\,\, (\nabla\cdot (\mu\, \mathbold{H}^n_{h,i}))^2 + (\nabla\cdot (\epsilon\, \mathbold{E}^n_{h,i}))^2\,\mathrm{d}\mathbold{x}\,\mathrm{d}t,
		\end{aligned}
	\end{equation*}
where $\ell_i = \beta\,h$ is the characteristic length of the spatial domain $S_i$ that depends on the mesh size $h$ and $\beta >0$.
The second part of the functional $J_i^n$ that weakly enforces the boundary conditions is either
	\begin{equation*}
		\mathcal{B}_i^n(\mathbold{E}^n_{h,i}) = \frac{1}{2} \,\int\limits_{I_n}\!\int\limits_{\Gamma\cap S_i} (\mathbold{n}\times\mathbold{E}^n_{h,i} - \mathbold{g}_E)\cdot(\mathbold{n}\times\mathbold{E}^n_{h,i} - \mathbold{g}_E)\, \mathrm{d}s\,\mathrm{d}t,
	\end{equation*}
	for the boundary condition \eqref{eq:nxE},  
	\begin{equation*}
		\mathcal{B}_i^n(\mathbold{H}^n_{h,i}) = \frac{1}{2} \, \int\limits_{I_n}\!\int\limits_{\Gamma\cap S_i} (\mathbold{n}\times\mathbold{H}^n_{h,i} - \mathbold{g}_H) \cdot (\mathbold{n}\times\mathbold{H}^n_{h,i} - \mathbold{g}_H)\, \mathrm{d}s\,\mathrm{d}t,
	\end{equation*}
	for the boundary condition \eqref{eq:nxH} {\color{black} or 
	\begin{equation*}
		\mathcal{B}_i^n(\mathbold{H}^n_{h,i},\mathbold{E}^n_{h,i}) = \frac{1}{2} \, \int\limits_{I_n}\!\int\limits_{\Gamma\cap S_i} (\mathbold{E}_{h,i}^n\times\mathbold{n} +Z\,\mathbold{n}\times(\mathbold{H}_{h,i}^n\times\mathbold{n}) - \mathbold{g}) \cdot (\mathbold{E}_{h,i}^n\times\mathbold{n} +Z\,\mathbold{n}\times(\mathbold{H}_{h,i}^n\times\mathbold{n}) - \mathbold{g})\, \mathrm{d}s\,\mathrm{d}t,
	\end{equation*}
	for the boundary condition \eqref{eq:impedance}.}
	
{\color{black} The final part of $J_i^n$ that weakly enforces the correction functions to match the Hermite solution is given by
	\begin{equation*} 
	\begin{aligned}
		\mathcal{H}_i^n(\mathbold{H}^n_{h,i},\mathbold{E}^n_{h,i}) =&\,\, \frac{c_H}{2} \int\limits_{t_{n-1}}^{t_{n-1/2}}\int\limits_{S^\mathcal{H}_{i,d}} (\mathbold{H}^n_{h,i}-\mathbold{H}^*)\cdot(\mathbold{H}^n_{h,i}-\mathbold{H}^*) 
		+ (\mathbold{E}^n_{h,i}-\mathbold{E}^*)\cdot(\mathbold{E}^n_{h,i}-\mathbold{E}^*)\,\mathrm{d}\mathbold{x}\,\mathrm{d}t \\
		+&\,\,\frac{c_H}{2} \int\limits_{t_{n-1/2}}^{t_{n}}\int\limits_{S^\mathcal{H}_{i,p}} (\mathbold{H}^n_{h,i}-\mathbold{H}^*)\cdot(\mathbold{H}^n_{h,i}-\mathbold{H}^*) 
		+ (\mathbold{E}^n_{h,i}-\mathbold{E}^*)\cdot(\mathbold{E}^n_{h,i}-\mathbold{E}^*)\,\mathrm{d}\mathbold{x}\,\mathrm{d}t.
	\end{aligned}
	\end{equation*}}
 
We then have the following problem statement: 
	\begin{equation} \label{eq:minPblm}
		\begin{aligned}
			&\text{Find } (\mathbold{H}^n_{h,i},\mathbold{E}^n_{h,i}) \in V \times V \text{ such that }\\
 			&\qquad (\mathbold{H}^n_{h,i},\mathbold{E}^n_{h,i}) =  \underset{\mathbold{v},\mathbold{w}\in V}{\arg\min}\, J_i^n(\mathbold{v},\mathbold{w}),
		\end{aligned}
	\end{equation}
	for $i=0,\dots,N_\Gamma$ and $n=1,\dots,N_t$.
Here 
\begin{equation*}
	V = \big\{ \mathbold{v} \in \big[\mathbb{Q}^k(S_i\times I_n )\big]^3 \big\}.
\end{equation*}

As in 1-D, 
	we use that the gradient of the functional $J_i^n$ with respect to the coefficients of the polynomial approximations $\mathbold{H}_{h,i}^n$ and $\mathbold{E}_{h,i}^n$ 
	vanishes at a minimum to obtain a linear system of equations to solve.
The dimension of the minimization problems is independent of the mesh size and the time step size, 
	and is $3\,(k+1)^{3}$ in 2-D and $6\,(k+1)^{4}$ in 3-D. 
However, 
	the number of minimization problems $(N_\Gamma+1)\,N_t$ increases as the mesh size and the time step size diminish.
Once the minimization problem is solved on a local patch, 
	the electromagnetic fields and their {\color{black} space derivatives through order $m$} are estimated at its CF node using 
	$\mathbold{H}^n_{h,i}$ and $\mathbold{E}^n_{h,i}$.

\begin{remark} 
The terms in $\mathcal{G}_i^n$ enforcing the residual of Maxwell's equations \eqref{eq:problemMaxwell} 
		are scaled by $\ell_i$ to guarantee that all the terms in $\mathcal{G}_i^n$ and $\mathcal{B}_i^n$ behave in a 
		similar way as the mesh size diminishes \cite{Marques2011}. 
Let us assume that the correction functions are polynomials of degree $k$ that leads to an accuracy 
	of $\mathcal{O}(\ell_i^{k+1})$ and that $k=2\,m$.
Using $\mathbold{H}_i^n = \mathbold{H} + \mathcal{O}(\ell_i^{k+1})$ and $\mathbold{E}_i^n = \mathbold{E} + \mathcal{O}(\ell_i^{k+1})$ in the functional $J_i^n$, 
	we have that the terms in $\mathcal{G}_i^n$ and $\mathcal{B}_i^n$ behave 
	as $\mathcal{O}(\ell_i^{2\,k+5})$ 
	while the term in $\mathcal{H}_i^n$ scales as $\mathcal{O}(\ell_i^{2\,k+6})$. 
Hence, 
	the functional $J_i^n$ is dominated by the boundary conditions and Maxwell's equations as $\ell_i$ diminishes.
\end{remark}

\begin{remark} 
The number of matrices to construct can be further reduced depending on the geometry of the domain and the physical properties of the material $\mu$ and $\epsilon$.
As an example, 
	let us consider a 2-D geometry discretized with a Cartesian mesh with $\Delta x = \Delta y$.
We also assume the boundary $\Gamma$ of the domain 
	to coincide only with primal nodes.
For problems with constant coefficients on a rectangular domain, 
	the number of matrices is reduced to eight because the spatial domain $S_i$ of local patches on an edge translates along it.
If reentrant corners are also considered, 
	there is a maximum of twelve matrices to compute. 
\end{remark}

\begin{remark} \label{rem:propertyHTCFM}
	Assuming that the correction functions are polynomials of degree $k$ that lead to an accuracy of $\mathcal{O}(\ell_i^{k+1})$, 
		we then have $k\geq2\,m$ to preserve the accuracy of a $(2\,m+1)$ order Hermite-Taylor method.
	As was remarked for FDTD methods in \cite{LawNave2021}, 
		the CFM impacts the stability of the original method because of the Hermite-Taylor polynomials $\mathbold{H}^*$ and $\mathbold{E}^*$.
	Since a rigorous proof of the stability of the proposed method is out of reach for the moment, 
		we investigate numerically its stability properties in Section~\ref{sec:num_examples}. 
\end{remark}

\section{Numerical Examples} \label{sec:num_examples}
In this section, we numerically investigate the stability of the proposed method and perform convergence studies in 1-D and 2-D. 

\subsection{Examples in One Dimension}

Let us seek approximate solutions to Maxwell's equations 
\begin{equation*} \label{eq:1DMaxwellEq}
\begin{aligned}
    \mu\,\partial_t H + \partial_x E =&\,\, 0, \\
    \epsilon\,\partial_t E + \partial_x H =&\,\, 0,
\end{aligned}
\end{equation*}
 	in the domain $\Omega = [x_\ell, x_r]$ and the time interval $I=[t_0,t_f]$.
The initial conditions are $H(x,t_0) = a(x)$ and $E(x,t_0) = b(x)$, 
	and we focus on the boundary conditions $E(x_\ell,t) = g_\ell(t)$ and $E(x_r,t) = g_r(t)$.
 Here $a(x)$,
    $b(x)$, 
    $g_\ell(t)$ and $g_r(t)$ are known functions.

In this subsection, 
	we use the Hermite-Taylor correction function method with $1\leq m \leq 5$. 
We set the degree of the correction functions to be $2\,m$.
The CFM should not therefore impact the convergence rate of the Hermite-Taylor method.

\subsubsection{Stability}
Let us first investigate the stability of the Hermite-Taylor correction function method.  
We  consider $\Omega = [0,1]$,  
	and set $\mu=1$ and $\epsilon =1$.
The stability condition of the Hermite-Taylor method depends only on the largest wave speed and is given here by 
	$\Delta t < h$,
	where $h$ is the mesh size. 
As mentioned in Remark~\ref{rem:propertyHTCFM}, 
	the stability of the Hermite-Taylor method is impacted by the CFM because we use Hermite-Taylor polynomials $\mathbold{H}^*$ and $\mathbold{E}^*$ 
	in the minimization problem \eqref{eq:minPblm}.
Although we do not have a rigorous proof of the stability of the Hermite-Taylor correction function method, 
	we provide numerical evidences of it by investigating the eigenvalues of the global matrix 
	associated with the method. 
	
Since Maxwell's equations is a linear system of PDEs and assuming $g_\ell = g_r = 0$,
	the proposed numerical method can be written as 
\begin{equation*}
	\mathbold{W}_p^{n+1} = A\,\mathbold{W}_p^n,
\end{equation*}
	where $A$ is a square matrix of dimension $2\,(N_x+1)\,(m+1)$ and 
	$\mathbold{W}_p^n$ is a vector containing all the degrees of freedom on the primal mesh at time $t_n$.
A stable method should have all the eigenvalues of $A$ inside the unit circle of the complex plane. 
In the following, 
	we compute numerically the eigenvalues of $A$ and consider that the scheme is stable if the spectral radius $\rho(A)$ of the matrix $A$ 
	is at most one with an error of $\mathcal{O}(10^{-10})$.
 \begin{figure} 
\begin{adjustbox}{max width=1.1\textwidth,center}
 \centering
		\includegraphics[width=2.5in,trim={1.5cm 0cm 1.75cm 0cm},clip]{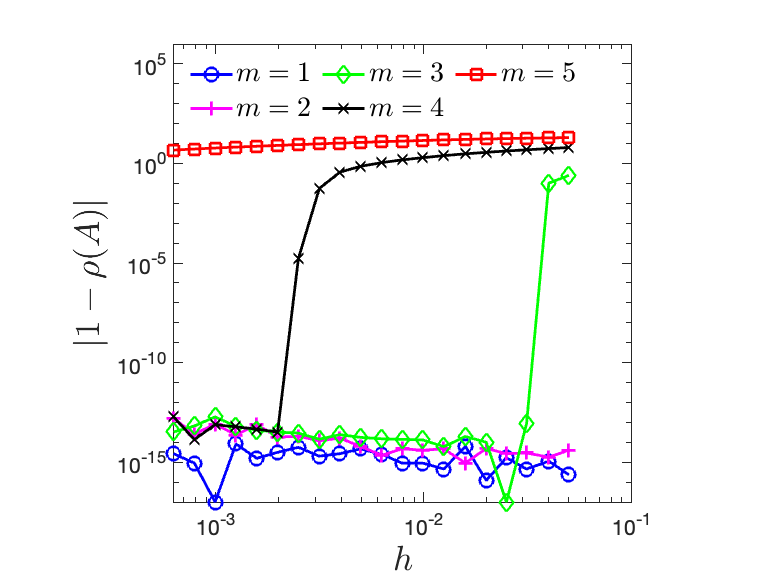}
		\includegraphics[width=2.5in,trim={1.5cm 0cm 1.75cm 0cm},clip]{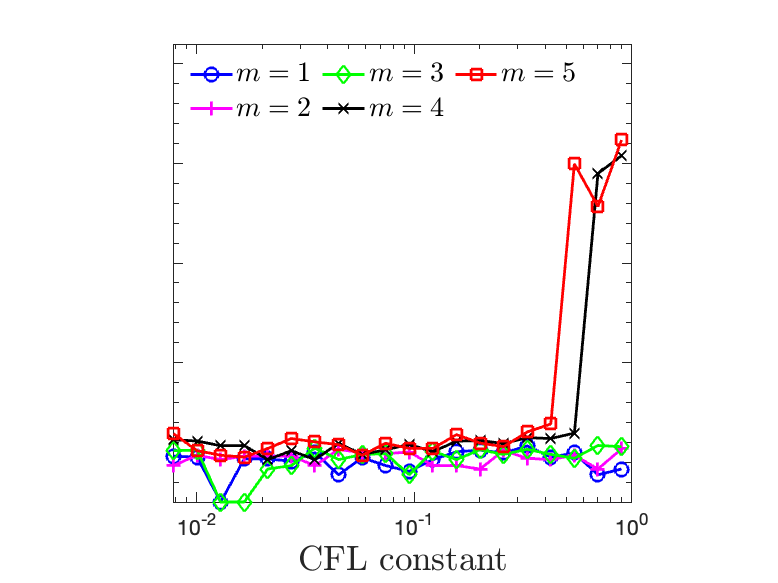}
	   	\includegraphics[width=2.5in,trim={1.5cm 0cm 1.75cm 0cm},clip]{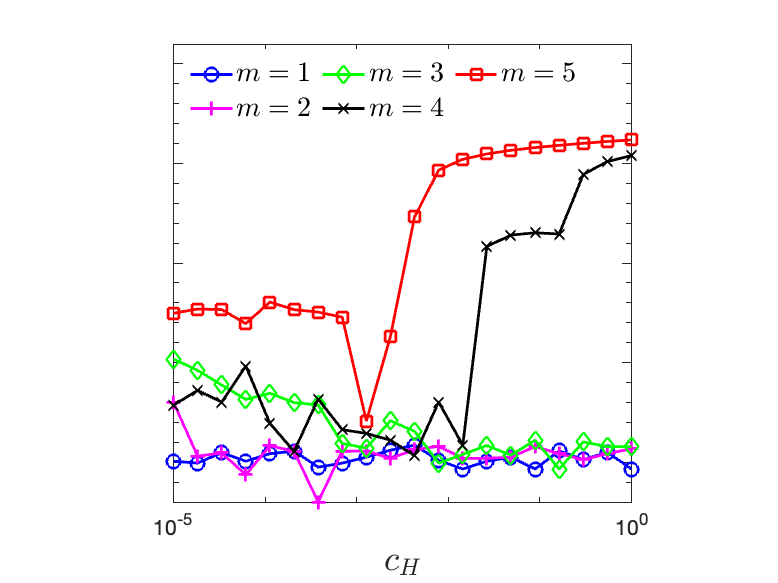}
\end{adjustbox}
  \caption{Absolute difference between one and the spectral radius of the matrix $A$ as a function of the mesh size, the CFL constant and the penalization parameter $c_H$ for various values of $m$. For the left plot, the CFL constant is set to $0.9$ and $c_H = 1$. For the middle plot, the mesh size is $h=\tfrac{1}{80}$ and $c_H=1$. For the right plot, the CFL constant is set to $0.9$ and $h=\tfrac{1}{80}$.} 
   \label{fig:spect_1D}
\end{figure}

The left plot of Fig.~\ref{fig:spect_1D} illustrates the absolute difference between one and the spectral radius of the matrix $A$,
	denoted $\rho(A)$,
	as a function of the mesh size for a CFL constant of $0.9$, 
	$c_H=1$ and various values of $m$.
{\color{black} For $m\leq4$, 
	we observe that the method is stable for a sufficiently small mesh size.}
In other words,
	the eigenvalues of $A$ are moving inside the unit circle as the mesh is refined.
This is expected since the terms in $\mathcal{H}_i^n$ impacting the stability
	scale as $\mathcal{O}(\ell_i^{2\,k+6})$ while the other terms in $J_i^n$ scale as $\mathcal{O}(\ell_i^{2\,k+5})$.
{\color{black} For $m=5$, 
	we do not observe a clear improvement as the mesh size diminishes for the considered CFL constant.}
This motivates us to diminish the CFL constant and the value of $c_H$ in order to improve the stability of the Hermite-Taylor correction function method.
	
The middle plot of Fig.~\ref{fig:spect_1D} illustrates the absolute difference between one and $\rho(A)$ 
	as a function of the CFL constant for $h=\tfrac{1}{80}$, 
	$c_H=1$ and various values of $m$.
{\color{black} We clearly have a stable method as the CFL constant diminishes.} 
	
The right plot of Fig.~\ref{fig:spect_1D} illustrates the absolute difference between one and the spectral radius of the matrix $A$ 
	as a function of $c_H$ for a CFL constant of $0.9$, 
	$h=\tfrac{1}{80}$ and various values of $m$.
For all $m$,
	a smaller value of the penalization coefficient $c_H$ helps to obtain a stable method. 
	
\begin{figure} 
\begin{adjustbox}{max width=1.10\textwidth,center}
 \centering
		\includegraphics[width=2.5in,trim={1.5cm 0cm 1.75cm 0cm},clip]{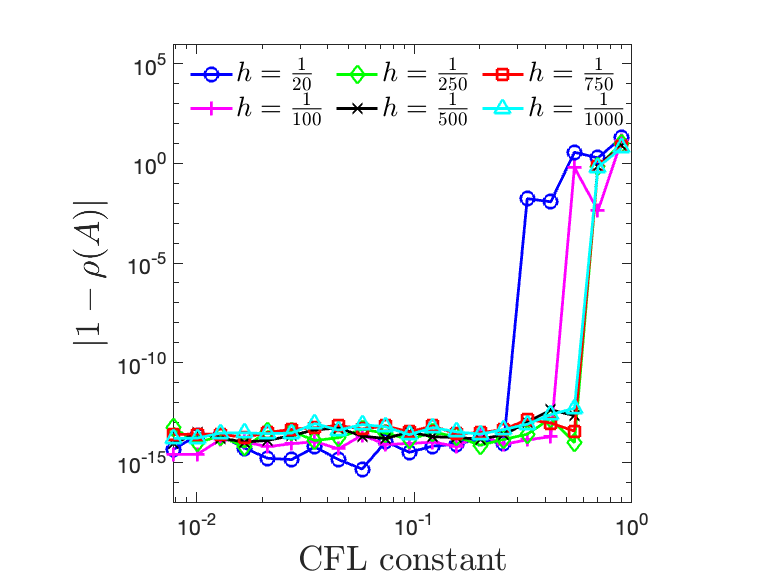}
		\includegraphics[width=2.5in,trim={1.5cm 0cm 1.75cm 0cm},clip]{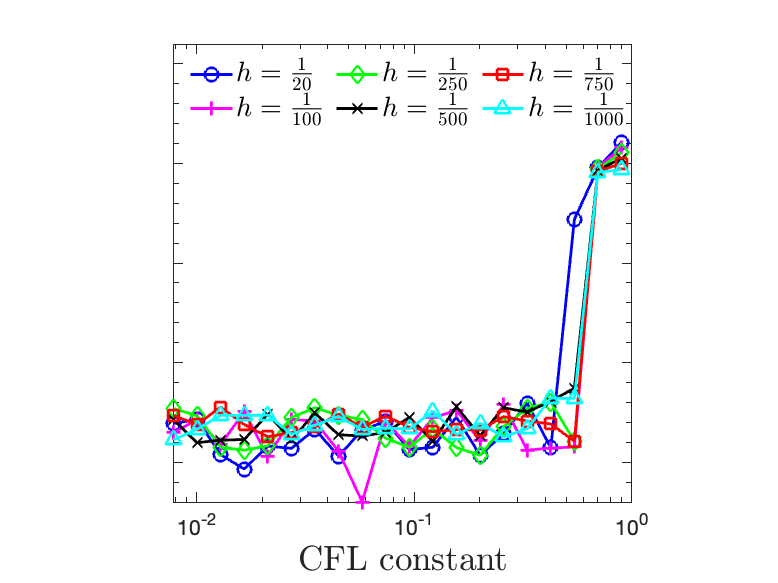}
	   	\includegraphics[width=2.5in,trim={1.5cm 0cm 1.75cm 0cm},clip]{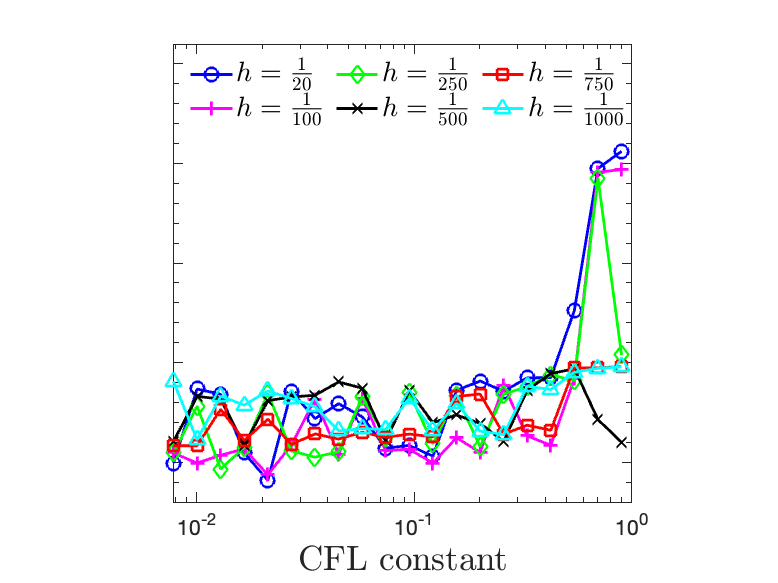}
\end{adjustbox}
  \caption{Absolute difference between one and the spectral radius of the matrix $A$ as a function of the CFL constant for {\color{black} $m=5$, 
  	and various mesh sizes and values of $c_H$. The left, 
	middle and right plots are respectively for $c_H=1$,
	$c_H=\frac{1}{10}$ and $c_H = \frac{1}{100}$.}} 
   \label{fig:spect_1D_CFL_m_3}
\end{figure}

To give further evidences of that, 
	Fig.~\ref{fig:spect_1D_CFL_m_3} illustrates the absolute difference between one and $\rho(A)$ as a function of the CFL constant for {\color{black} $m=5$},
	$h\in\big\{\tfrac{1}{20},\tfrac{1}{100},\tfrac{1}{250},\tfrac{1}{500},\tfrac{1}{750},\tfrac{1}{1000}\big\}$ and 
	$c_H \in \big\{1,\tfrac{1}{10},\tfrac{1}{100}\big\}$.
A smaller penalization coefficient $c_H$ does not improve the stability of the proposed method 
	for coarser meshes.
In these cases, 
	we therefore need to lower the CFL constant.
Based on these results, 
	the stability of the Hermite-Taylor correction function method improves by reducing the CFL constant and the value of 
	the penalization coefficient $c_H$.	
Moreover,
	the stability of this method improves as the mesh size diminishes, 
	which suggests that larger CFL constants could be used for finer meshes.

\subsubsection{Condition Number of CFM Matrices}
\begin{figure} 
\begin{adjustbox}{max width=1.0\textwidth,center}
 \centering
		\includegraphics[width=2.5in,trim={1.5cm 0cm 1.75cm 0cm},clip]{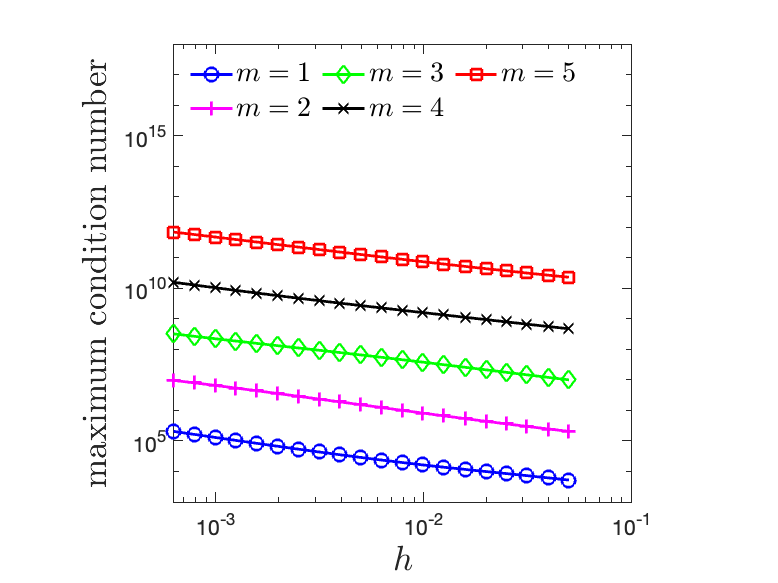}
		\includegraphics[width=2.5in,trim={1.5cm 0cm 1.75cm 0cm},clip]{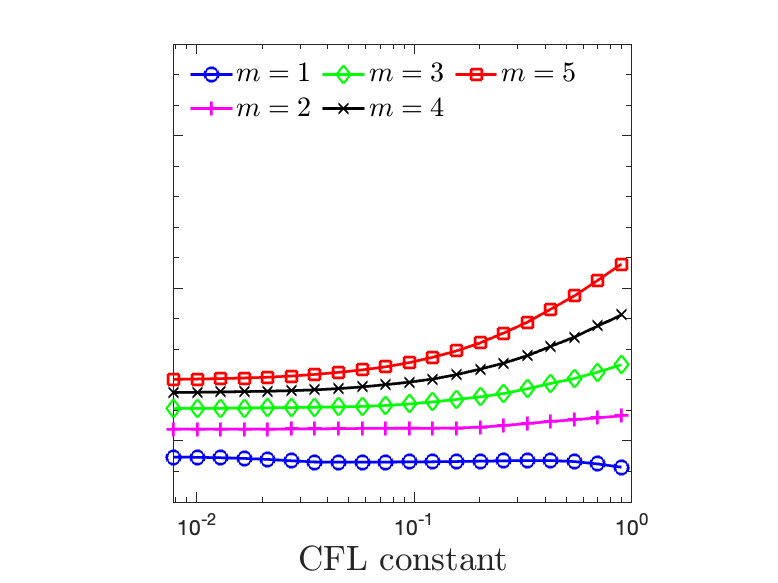}
	   	\includegraphics[width=2.5in,trim={1.5cm 0cm 1.75cm 0cm},clip]{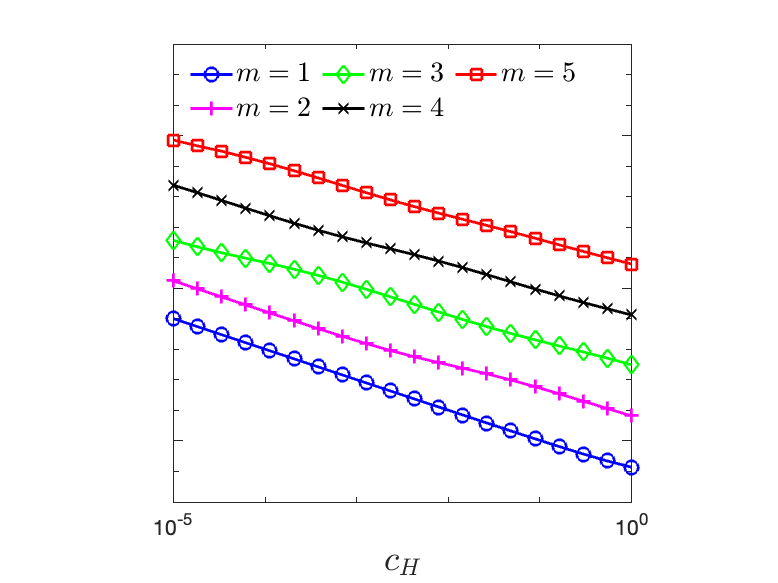}
\end{adjustbox}
  \caption{Maximum condition number of the matrices coming from the CFM as a function of the mesh size, the CFL constant and the penalization parameter $c_H$ for various values of $m$. For the left plot, the CFL constant is set to $0.9$ and $c_H = 1$. For the middle plot, the mesh size is $h=\tfrac{1}{80}$ and $c_H=1$. For the right plot, the CFL constant is set to $0.9$ and $h=\tfrac{1}{80}$.}
   \label{fig:cond_1D}
\end{figure}

Let us now investigate the impact of $h$, 
	$c_H$ and the CFL constant on the condition number of the matrices $M_i$ coming from the minimization procedure used in the CFM.
Fig.~\ref{fig:cond_1D} illustrates the maximum condition number of these matrices as a function of the mesh size, 
	the CFL constant and the penalization parameter $c_H$ for various values of $m$.
We observe that the condition number increases as the mesh size diminishes and, 
	more precisely, 
	scales as $\tfrac{1}{h}$ for all different settings. 
{\color{black} We also notice that the condition number first diminishes as the CFL constant decreases, 
	then appears to stabilize at a constant.
Finally, 
	the condition number increases as $c_H$ diminishes and scales as $\frac{1}{c_H}$.}
Hence,
	an arbitrary small value of $c_H$ cannot be taken to avoid poorly conditioned matrices coming from the CFM.
{\color{black} It is then preferable to diminish the CFL constant to obtain a stable method.}

\subsubsection{Accuracy}
{\color{black} In the following, 
	we use $c_H=1$ and $k=2\,m$ for all settings. 
We set the CFL constant at $0.9$ for $m=1$ and $m=2$, 
	$0.5$ for $m=3$ and $m=4$, 
	and $0.25$ for $m=5$.
The computed spectral radius is maximum one up to an error of $10^{-12}$ for all considered mesh sizes.}

\begin{figure}[htb]
 \centering
	\includegraphics[width=2.5in]{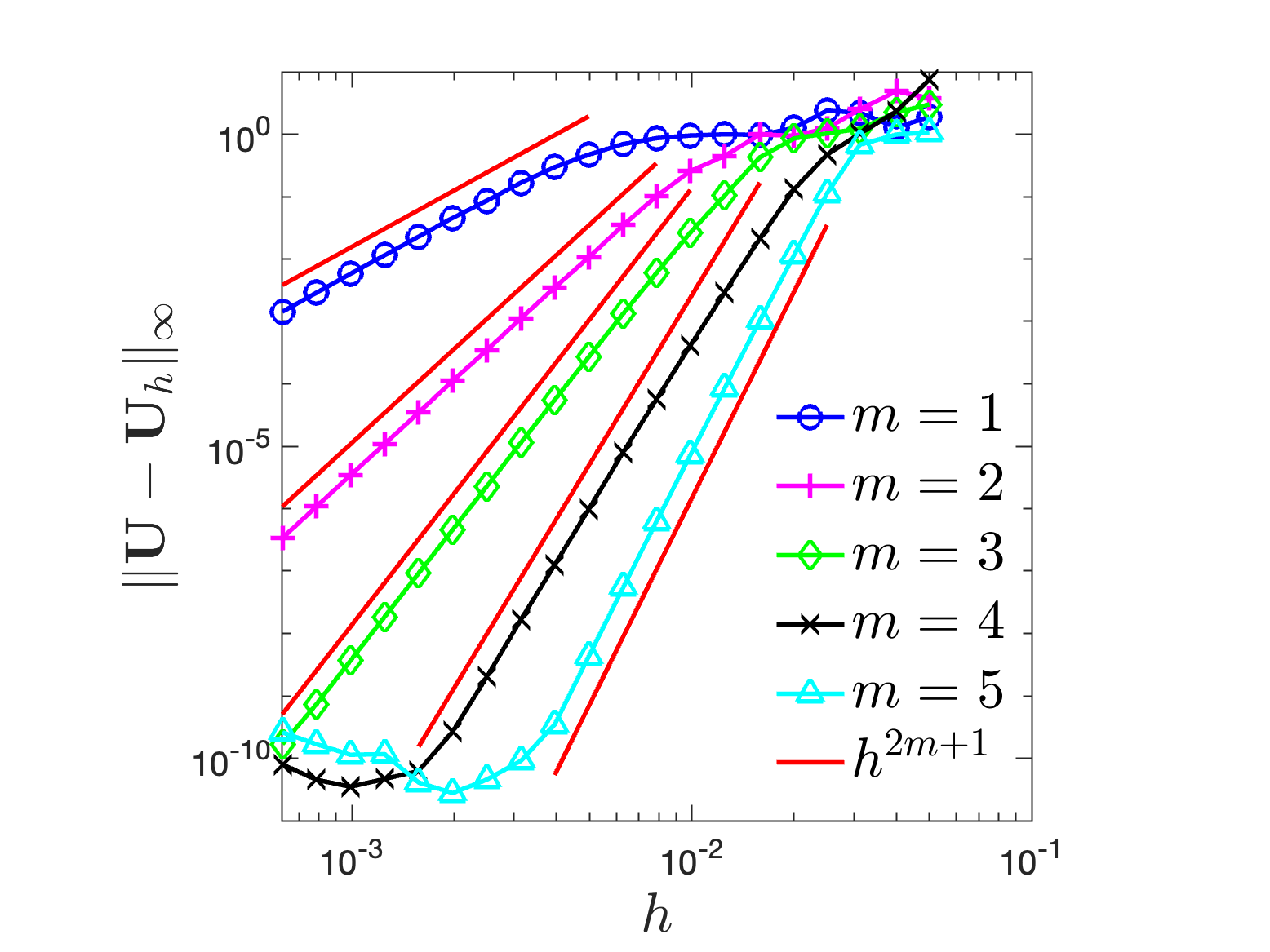}
  \caption{Convergence plots in the maximum norm for a standing mode problem using the Hermite-Taylor correction function method with different values of $m$ in 1-D. Here $\mathbold{U} = [H,E]^T$. }
   \label{fig:convPlots_1D}
\end{figure}

Let us now verify the convergence order of the proposed method. 
We consider a domain $\Omega = [\tfrac{1}{3},\tfrac{4}{3}]$, 
    a time interval $I = [0,1]$,
    $\mu=1$ and 
    $\epsilon=1$.
We set the initial and boundary data so find that the solution to the problem is 
\begin{equation} \label{eq:sol1D}
	\begin{aligned}
		H(x,t) =&\,\, \sin(250\,x)\,\sin(250\,t), \\
		E(x,t)=&\,\, \cos(250\,x)\,\cos(250\,t).
	\end{aligned}
\end{equation}
Fig.~\ref{fig:convPlots_1D} shows how the errors follow the expected $(2\,m+1)$ rates of convergence. 

\subsection{Examples in Two Dimensions}
\begin{figure}
 	\centering
	\includegraphics[width=1.8in]{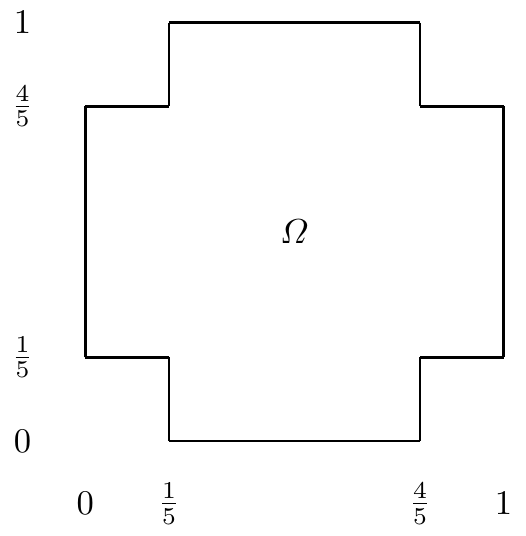}
  \caption{Geometry of a cross domain in 2-D.}
   \label{fig:cross_geo_2D}
\end{figure}

Let us consider the transverse magnetic (TM$_z$) mode.
We seek approximate solutions to Maxwell's equations 
\begin{equation} \label{eq:TMzSyst}
	\begin{aligned}
	\mu\,\partial_t H_x + \partial_y E_z =&\,\, 0,\\
	\mu\,\partial_t H_y - \partial_x E_z  =&\,\, 0,\\
	\epsilon\,\partial_t E_z - \partial_x H_y + \partial_y H_x =&\,\, 0,\\
	\partial_x H_x + \partial_y H_y =&\,\, 0,\\
	\end{aligned}
\end{equation}
	in the domain $\Omega\subset \mathbb{R}^2$ and the time interval $I$, 
	with initial conditions for $H_x$, 
	$H_y$ and $E_z$.
The boundary conditions are either 
\begin{equation} \label{eq:Ez_bnd_cdns_2D}
	E_z = g_E,
\end{equation}
\begin{equation} \label{eq:H_bnd_cdns_2D}
	n_x\,H_y  - n_y\,H_x =g_H
\end{equation}
{\color{black} or 
\begin{equation} \label{eq:impedance_bnd_cdns_2D}
	\begin{bmatrix}
	-n_y\,E_z + Z\,n_y\,(n_y\,H_x-n_x\,H_y)\\
	n_x\,E_z - Z\,n_x\,(n_y\,H_x-n_x\,H_y)
	\end{bmatrix}
	= \mathbold{g}.
\end{equation}}
We consider two geometries of the domain, 
	that is a square $\Omega =  [\tfrac{1}{3},\tfrac{4}{3}]\times[\tfrac{1}{6},\tfrac{7}{6}]$ and one with reentrant corners, 
	which is named cross domain and is illustrated in Fig.~\ref{fig:cross_geo_2D}.
{\color{black} We set $k=2\,m$ and $c_H=1$.
The CFL constant is $0.9$ for $m=1$,
	$0.5$ for $m=2$ and $0.25$ for $m=3$.}
In the following, 
	we numerically investigate the stability of the Hermite-Taylor correction function method and perform convergence studies for both geometries.

\subsubsection{Stability}
Since the total number of degrees of freedom on the primal mesh in 2-D,
	given by $3\,(N_x+1)\,(N_y+1)\,(m+1)^2$ where $N_x$ and $N_y$ are the number of cells 
	in respectively the $x$ and $y$ direction,
	is very large, 
	we cannot compute the spectral radius of the matrix $A$ for small mesh sizes, 
	as in 1-D.
To provide numerical evidences of the stability of the proposed method, 
	we therefore compute the maximum norm of the electromagnetic fields over 10000 time steps using 
	the trivial solution, 
	but with initial data, 
	that is the electromagnetic fields and their {derivatives through order $m$}, 
	to be random numbers in $]-10\,\epsilon_M,10\,\epsilon_M[$.
Here $\epsilon_M$ is the machine precision. 
We set $\mu=1$ and $\epsilon=1$.
{\color{black}
Fig.~\ref{fig:long_time_stability_2D_cross_domain} illustrates the evolution of the maximum norm of the electromagnetic fields 
	using different values of $m$ and boundary conditions for the cross domain and different mesh sizes. 
 \begin{figure} 
	\begin{adjustbox}{max width=1.0\textwidth,center}
	 \centering
		\includegraphics[width=2.5in,trim={1.5cm 0cm 1.75cm 0cm},clip]{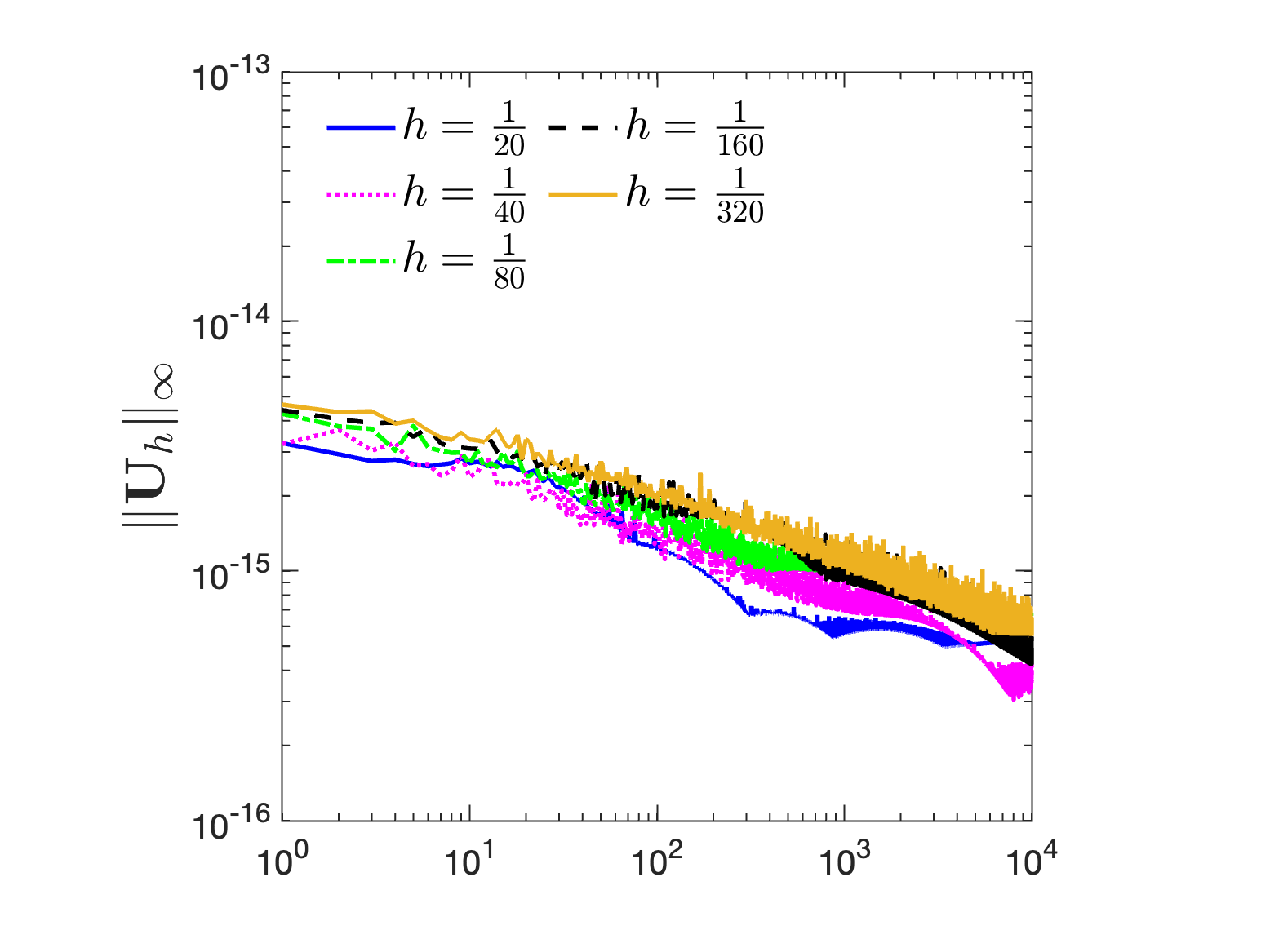}
		\includegraphics[width=2.5in,trim={1.5cm 0cm 1.75cm 0cm},clip]{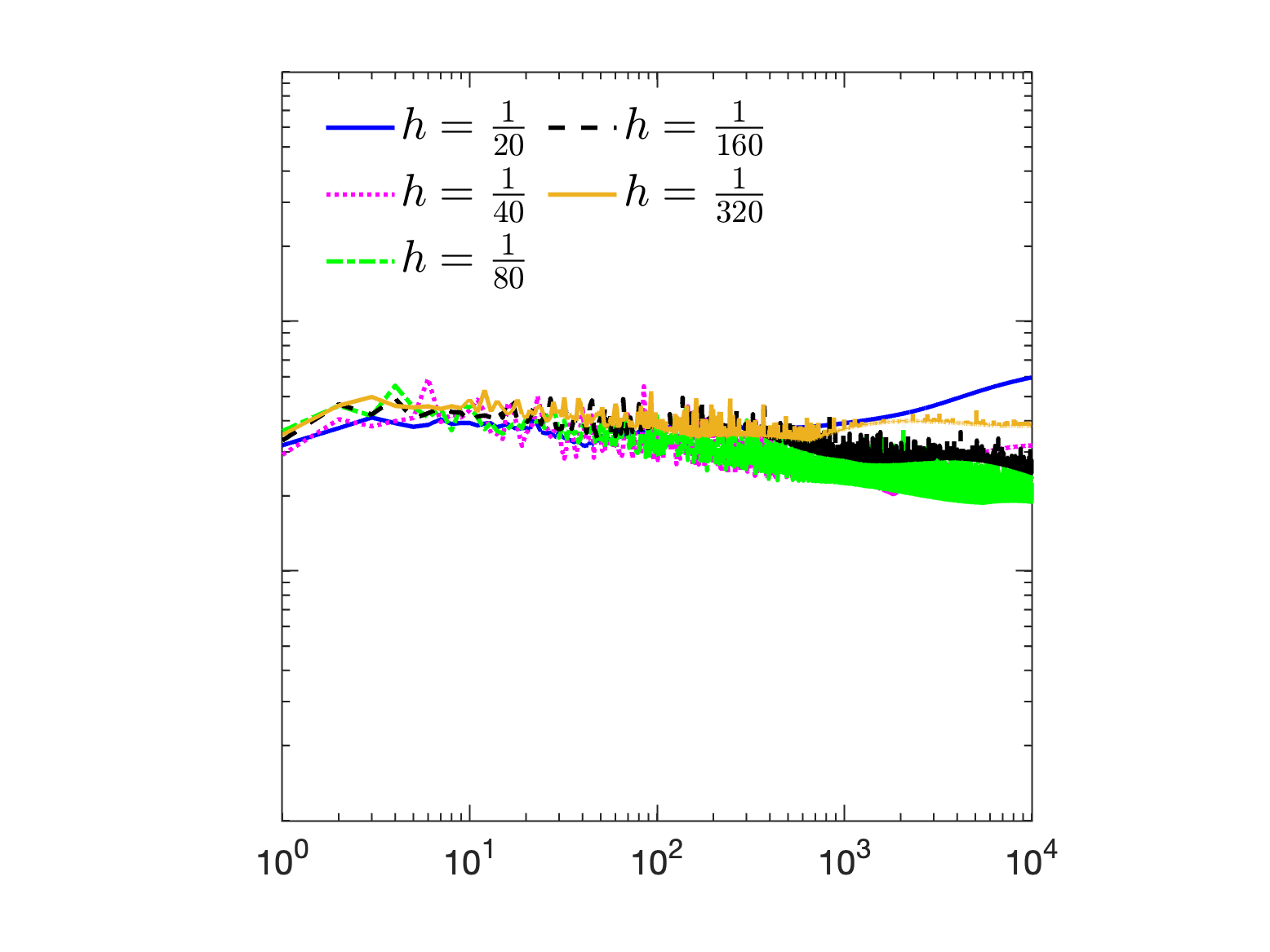}
	   	\includegraphics[width=2.5in,trim={1.5cm 0cm 1.75cm 0cm},clip]{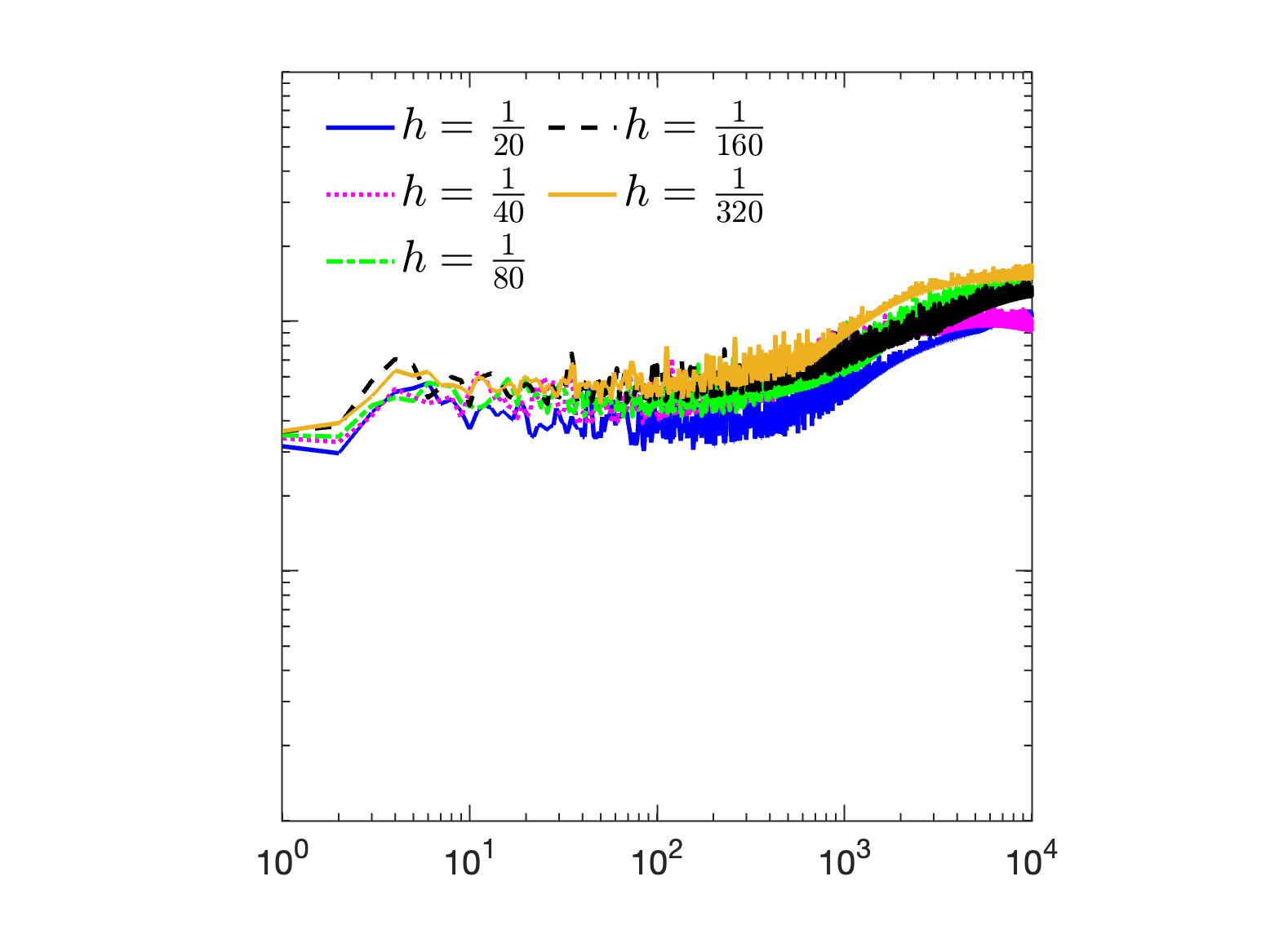} 
	\end{adjustbox}	
	\begin{adjustbox}{max width=1.0\textwidth,center}
	 \centering
		\includegraphics[width=2.5in,trim={1.5cm 0cm 1.75cm 0cm},clip]{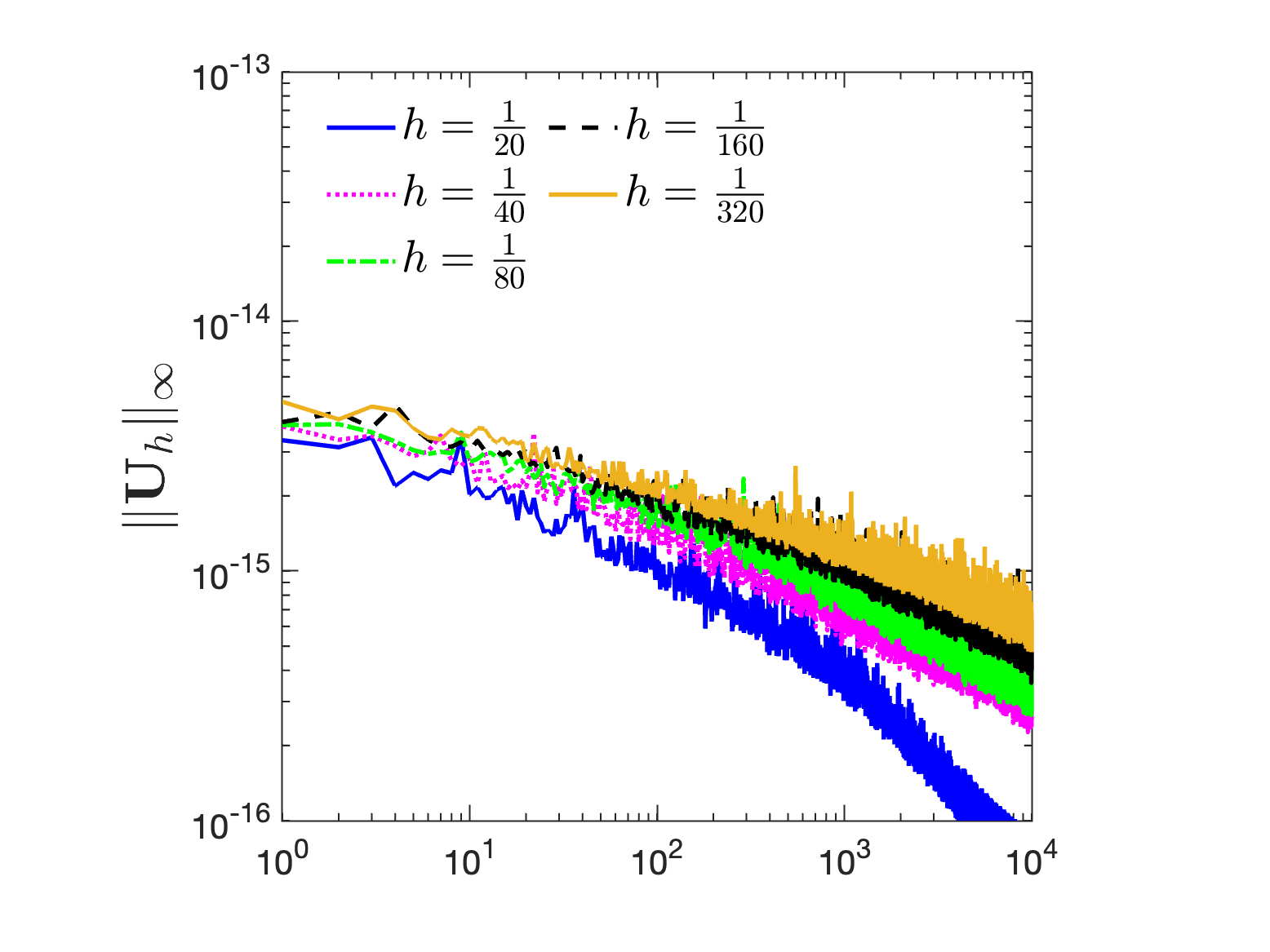}
		\includegraphics[width=2.5in,trim={1.5cm 0cm 1.75cm 0cm},clip]{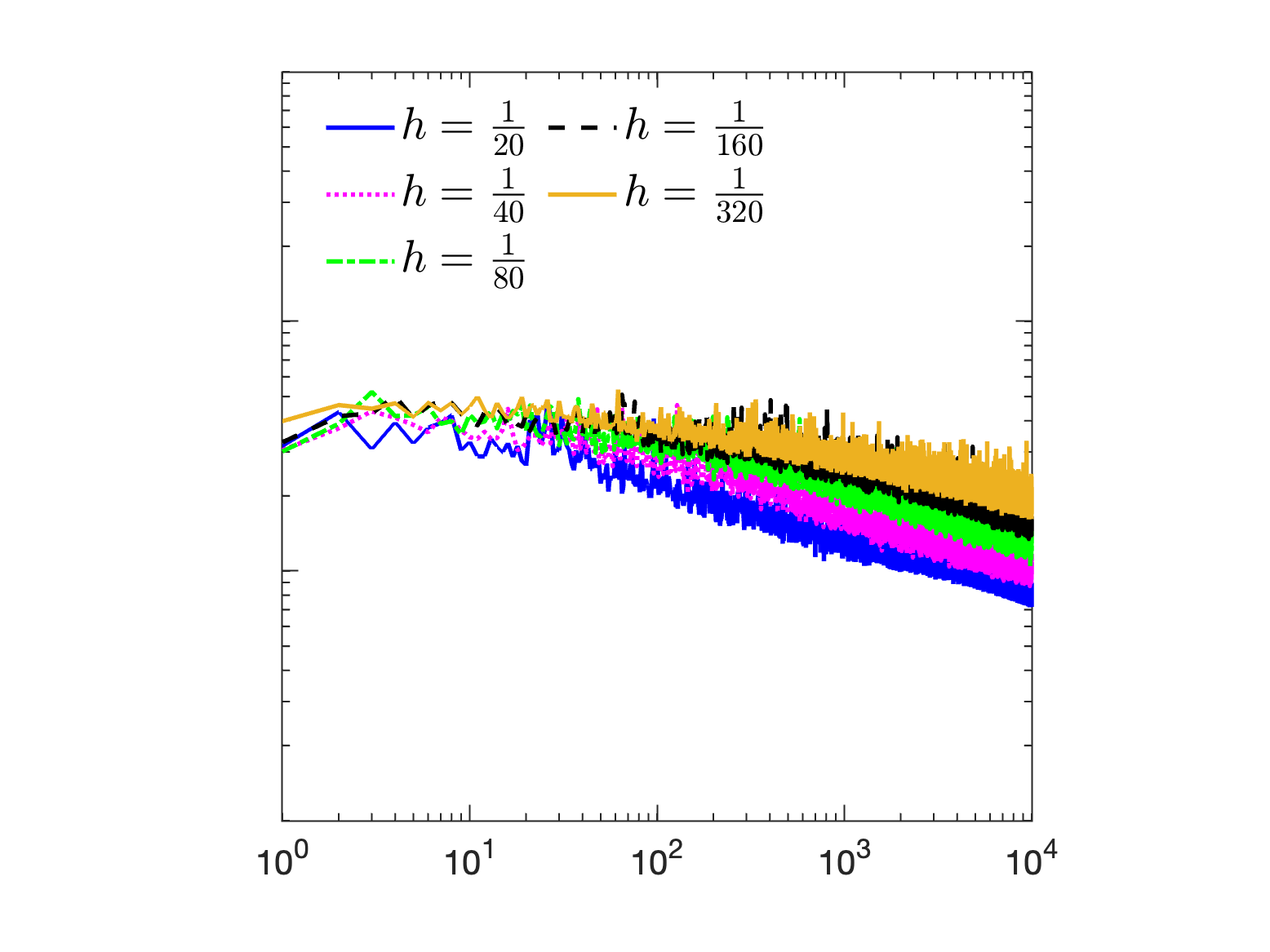}
	   	\includegraphics[width=2.5in,trim={1.5cm 0cm 1.75cm 0cm},clip]{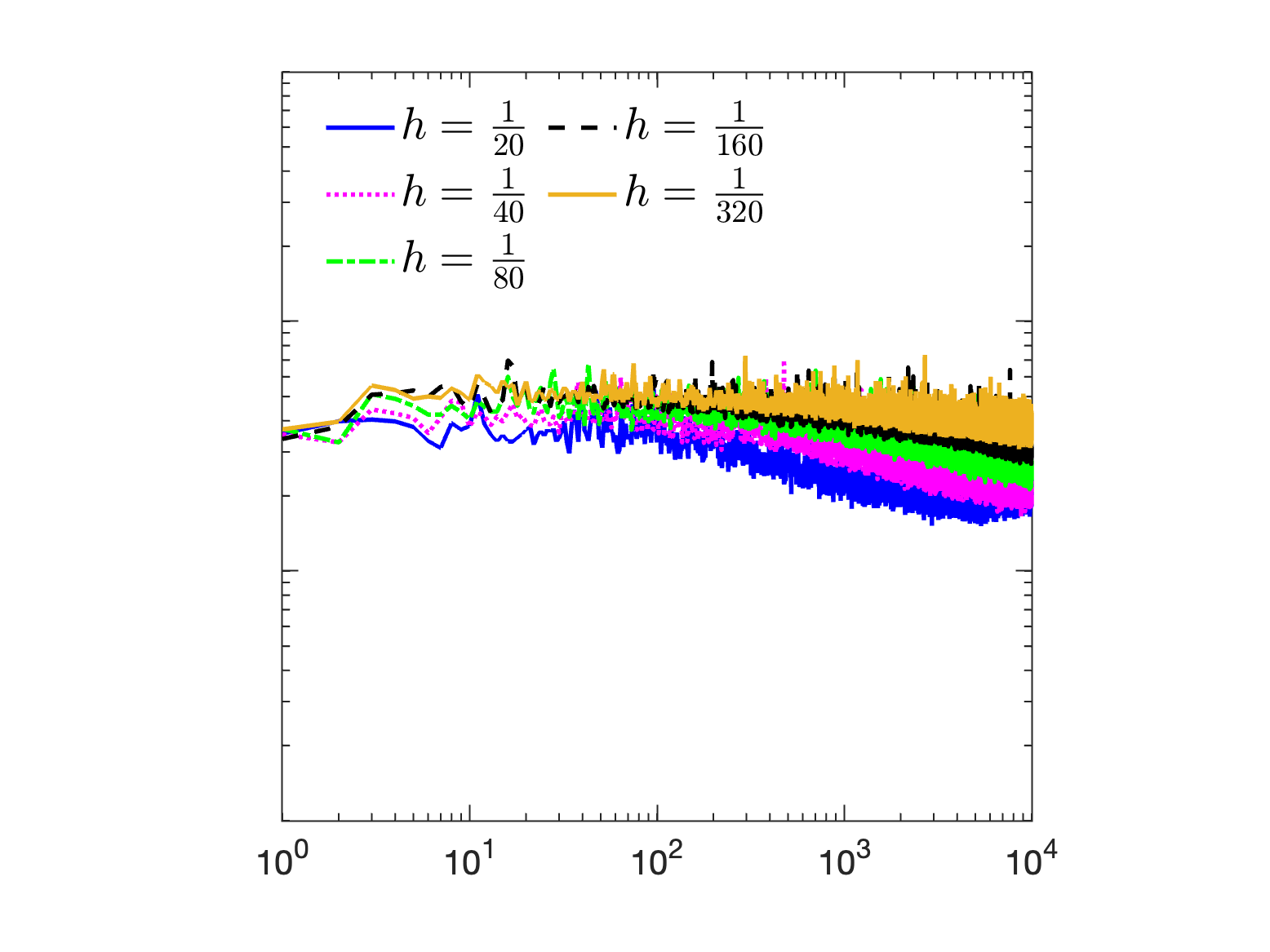} 
	\end{adjustbox}	
	\begin{adjustbox}{max width=1.0\textwidth,center}
	 \centering
		\includegraphics[width=2.5in,trim={1.5cm 0cm 1.75cm 0cm},clip]{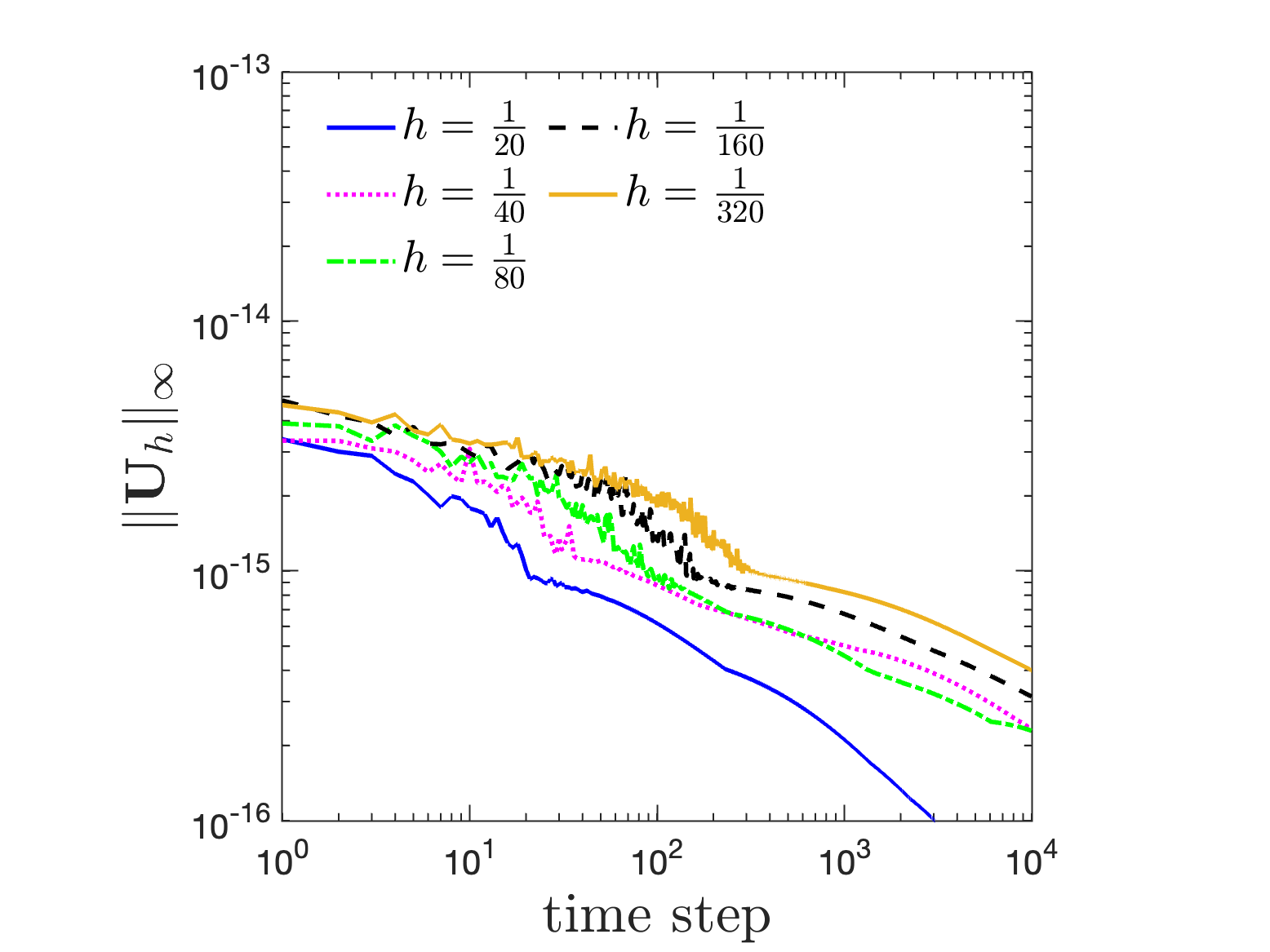}
		\includegraphics[width=2.5in,trim={1.5cm 0cm 1.75cm 0cm},clip]{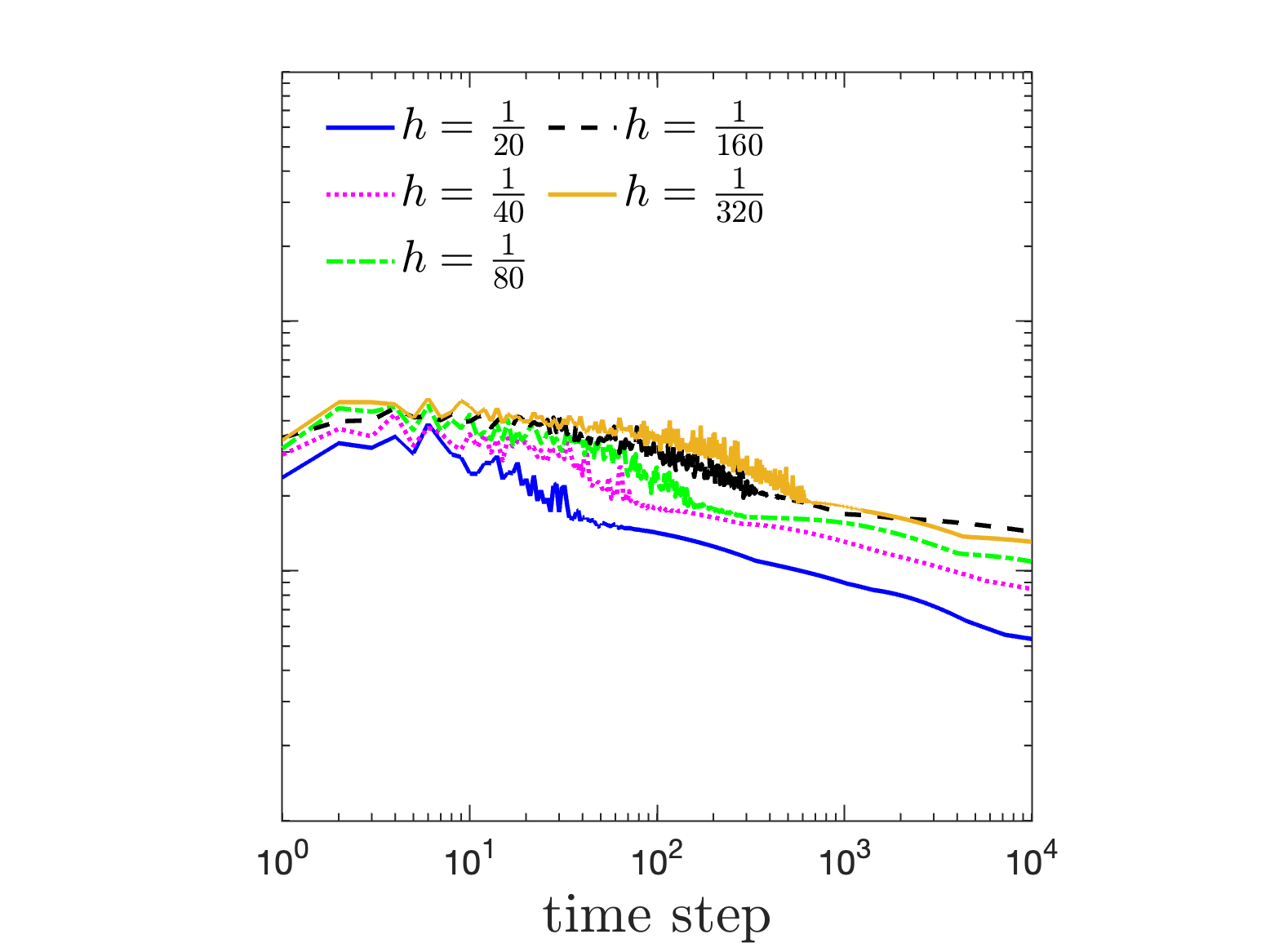}
	   	\includegraphics[width=2.5in,trim={1.5cm 0cm 1.75cm 0cm},clip]{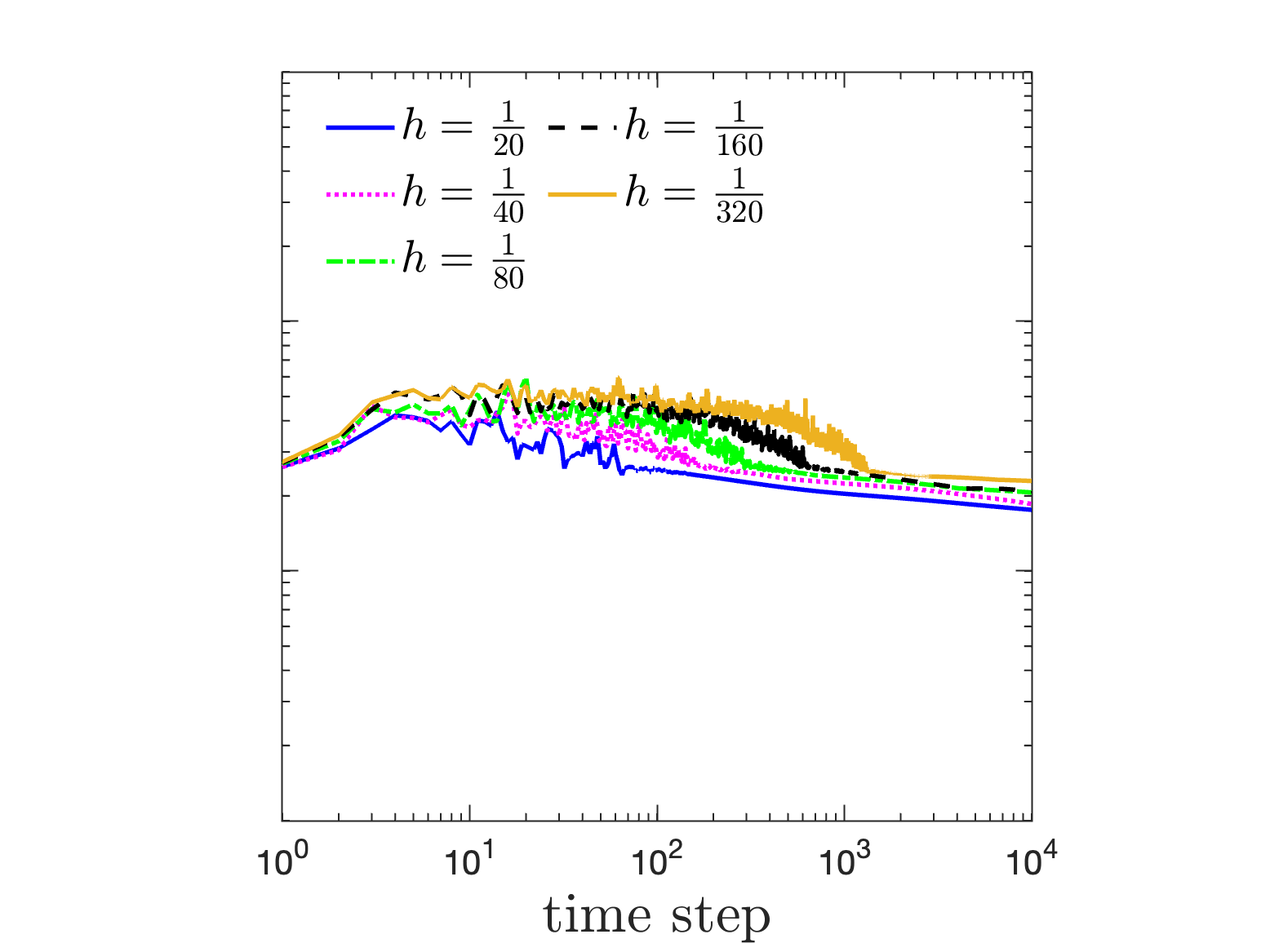}
	\end{adjustbox}	
  \caption{Evolution of the maximum norm of the numerical solution for different values of $m$ and boundary conditions using the cross domain in 2-D. 
  		The left, 
		middle and right columns are respectively for $m=1$, 
		$m=2$ and $m=3$. 
		The top, 
		middle and bottom rows are for the boundary conditions \eqref{eq:Ez_bnd_cdns_2D}, 
		\eqref{eq:H_bnd_cdns_2D} and \eqref{eq:impedance_bnd_cdns_2D}. 
		Here $\mathbold{U} = [H_x,H_y,E_z]^T$.}
   \label{fig:long_time_stability_2D_cross_domain}
\end{figure}
These results suggest that the method is stable.}

\subsubsection{Accuracy}
For the convergence studies, 
	we consider the time interval $I=[0,1]$,
	and set $\mu=1$ and $\epsilon = 1$.
The initial conditions and boundary conditions are chosen in such a way that the solution is given by
\begin{equation*}
	\begin{aligned}
	H_x =&\,\, -\frac{1}{\sqrt{2}}\sin(\omega\,\pi\,x)\,\cos(\omega\,\pi\,y)\,\sin(\sqrt{2}\,\omega\,\pi\,t), \\
	H_y =&\,\, \frac{1}{\sqrt{2}}\cos(\omega\,\pi\,x)\,\sin(\omega\,\pi\,y)\,\sin(\sqrt{2}\,\omega\,\pi\,t), \\ 
	E_z =&\,\, \sin(\omega\,\pi\,x)\,\sin(\omega\,\pi\,y)\,\cos(\sqrt{2}\,\omega\,\pi\,t),
	\end{aligned}
\end{equation*}
	with $\omega = 20$.
{\color{black}
Fig.~\ref{fig:convPlots_2D_EM} illustrates convergence plots for different values of $m$, 
  			boundary conditions and geometries in 2-D.
As expected, 
	we observe a $(2\,m+1)$ rate of convergence in the maximum norm for the 
	Hermite-Taylor correction function method.
 \begin{figure} 
 \centering
	\begin{adjustbox}{max width=1.0\textwidth,center}
	 \centering
		\includegraphics[width=2.5in,trim={1.5cm 0cm 1.75cm 0cm},clip]{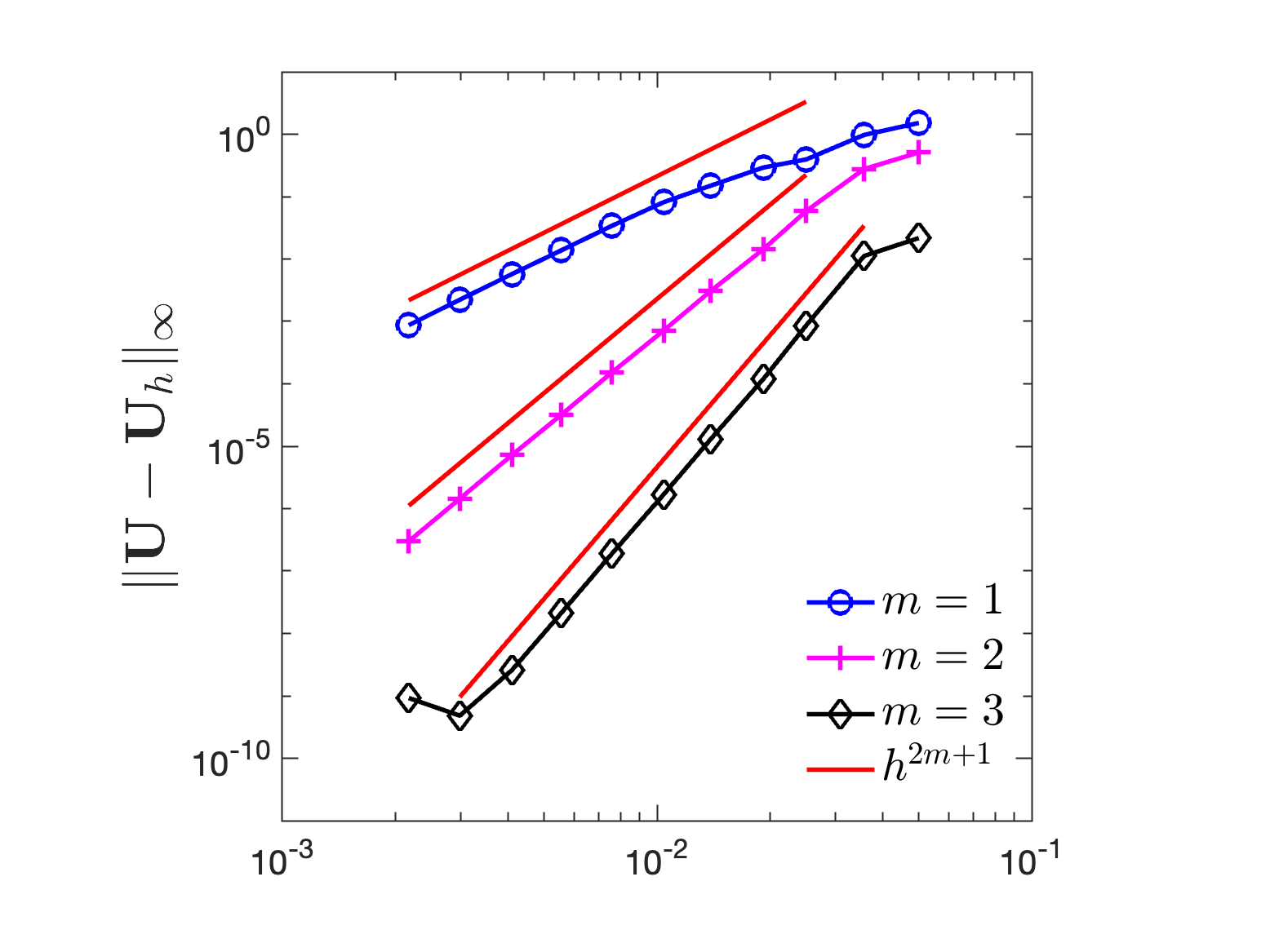}
		\includegraphics[width=2.5in,trim={1.5cm 0cm 1.75cm 0cm},clip]{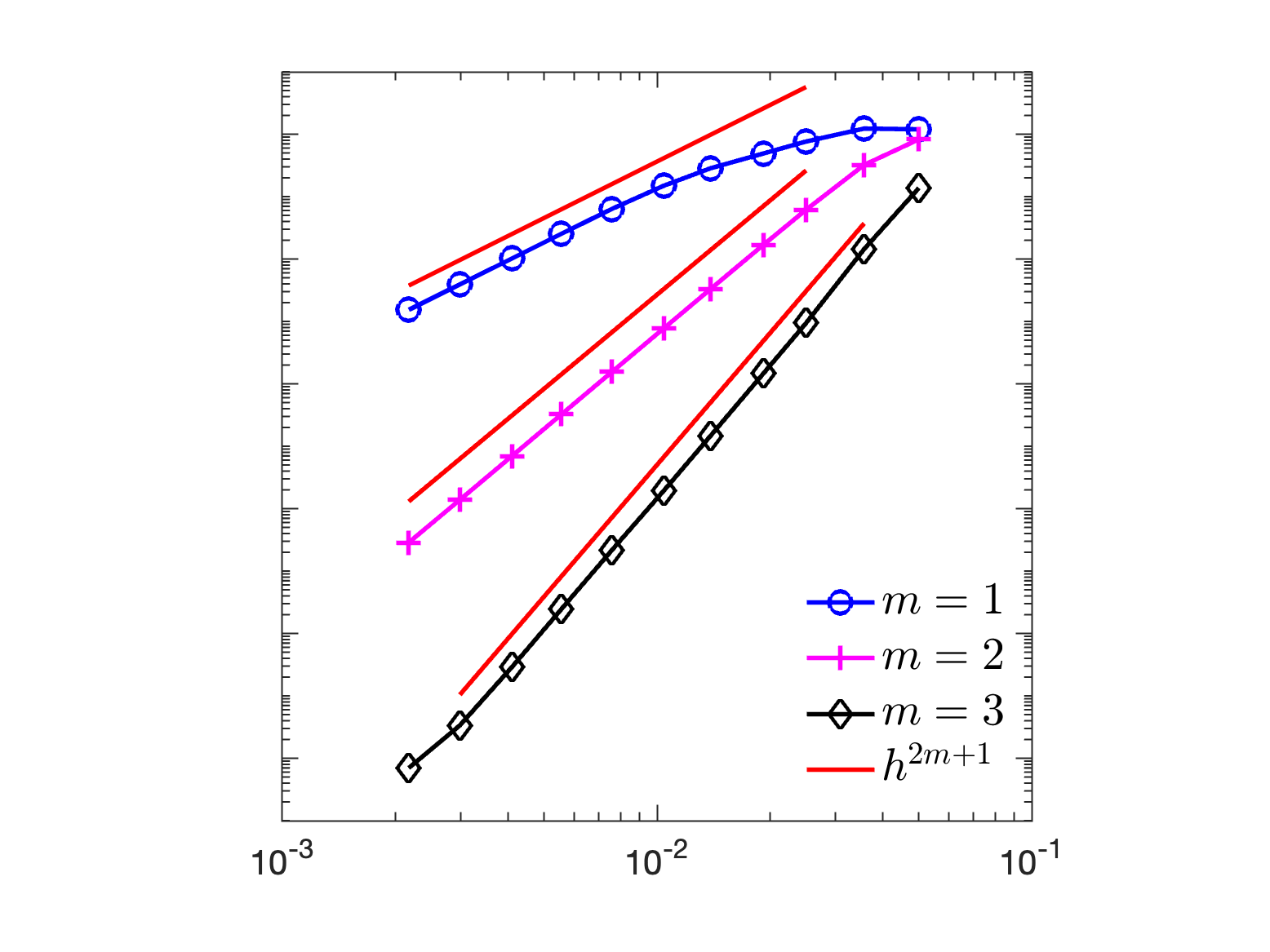}
	   	\includegraphics[width=2.5in,trim={1.5cm 0cm 1.75cm 0cm},clip]{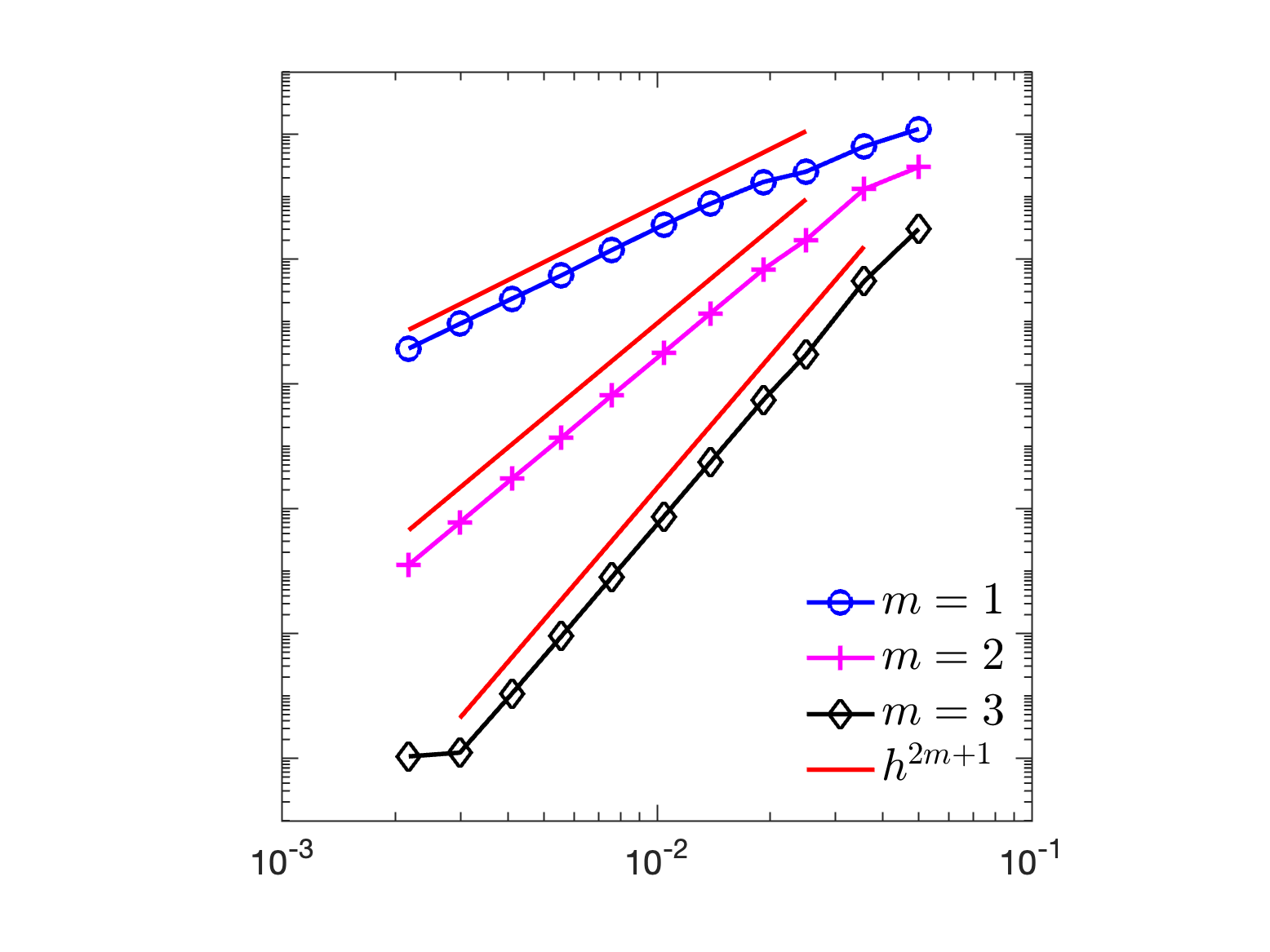} 
	\end{adjustbox}	
	\begin{adjustbox}{max width=1.0\textwidth,center}
	 \centering
		\includegraphics[width=2.5in,trim={1.5cm 0cm 1.75cm 0cm},clip]{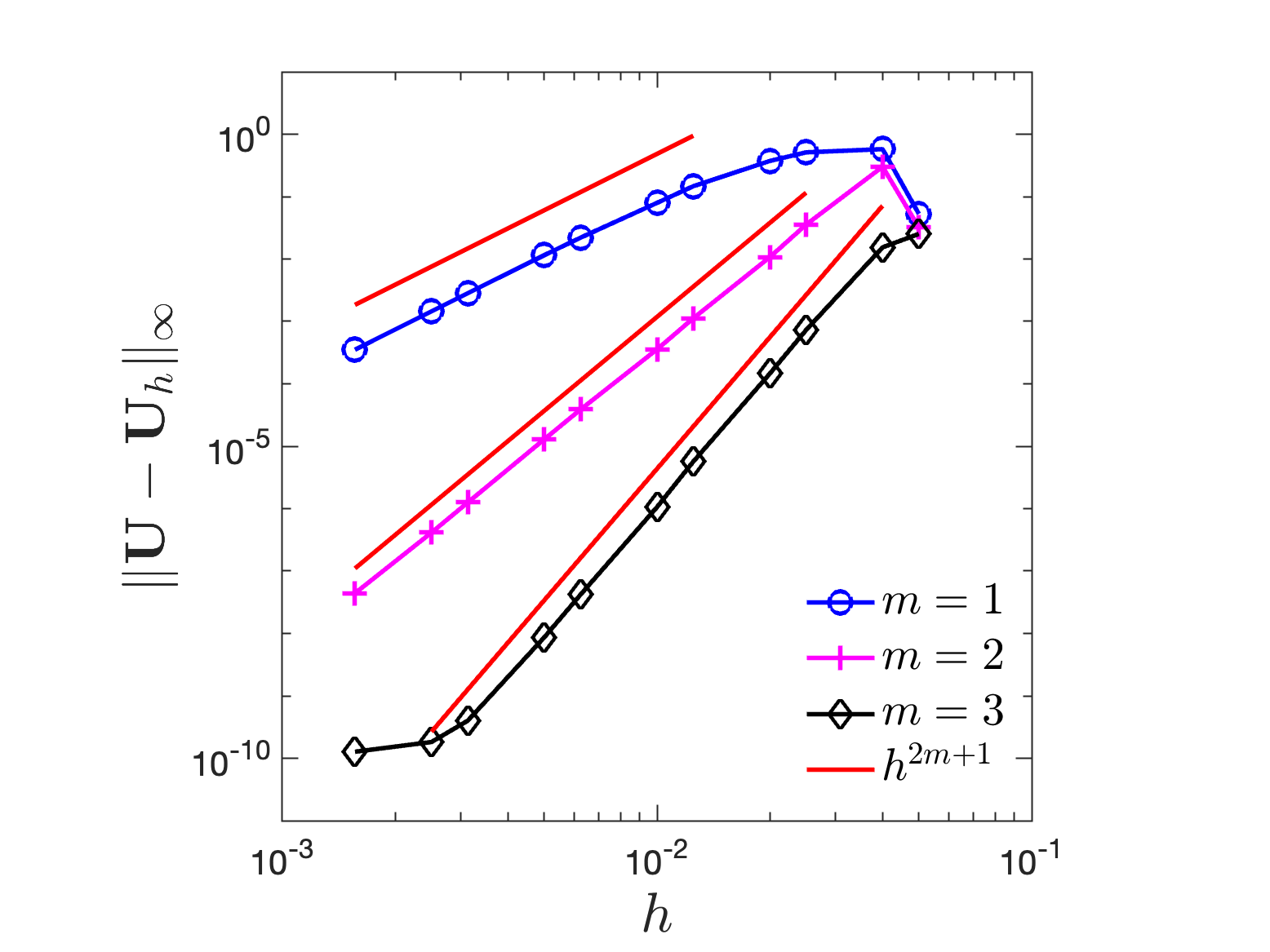}
		\includegraphics[width=2.5in,trim={1.5cm 0cm 1.75cm 0cm},clip]{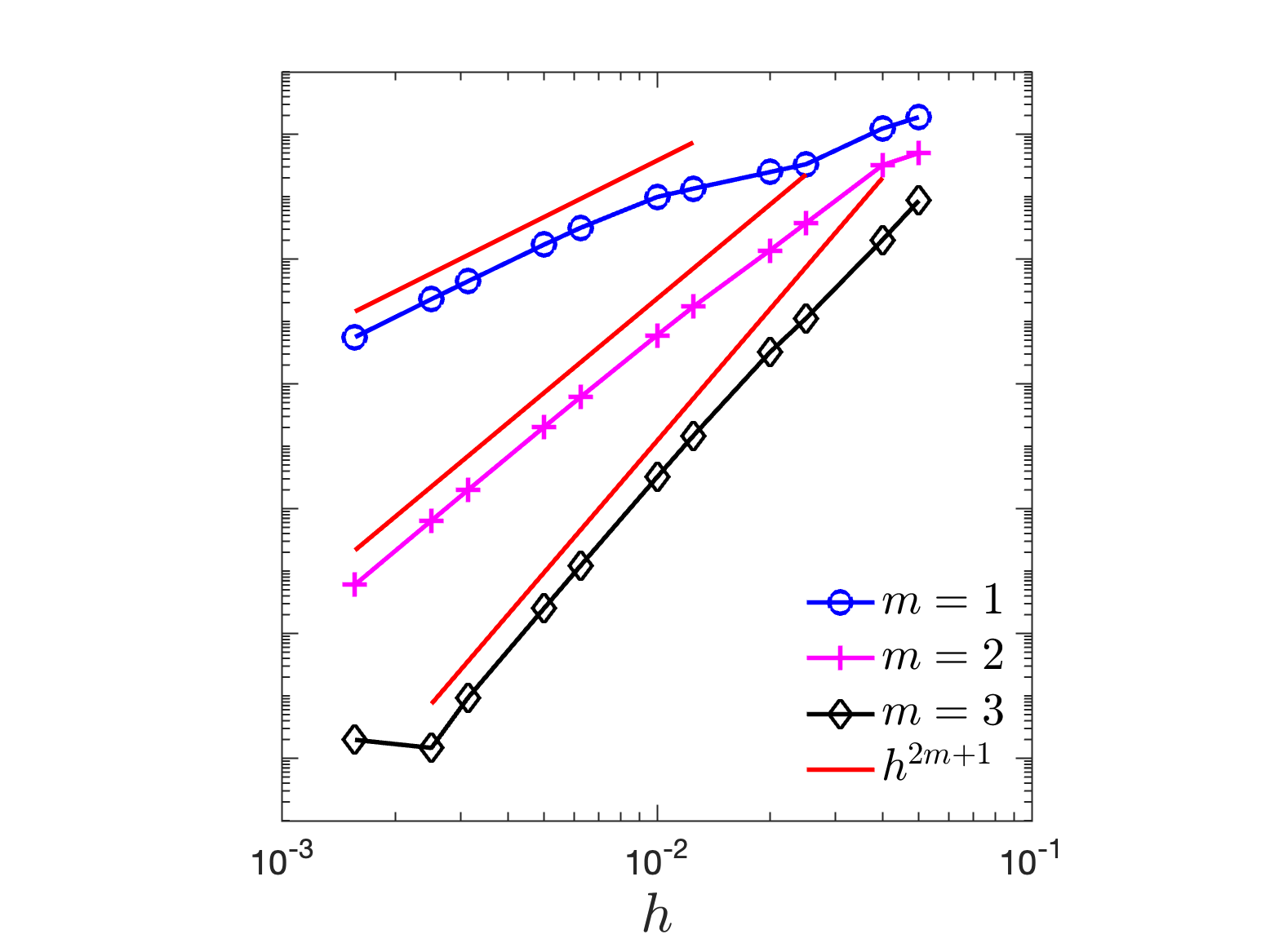}
	   	\includegraphics[width=2.5in,trim={1.5cm 0cm 1.75cm 0cm},clip]{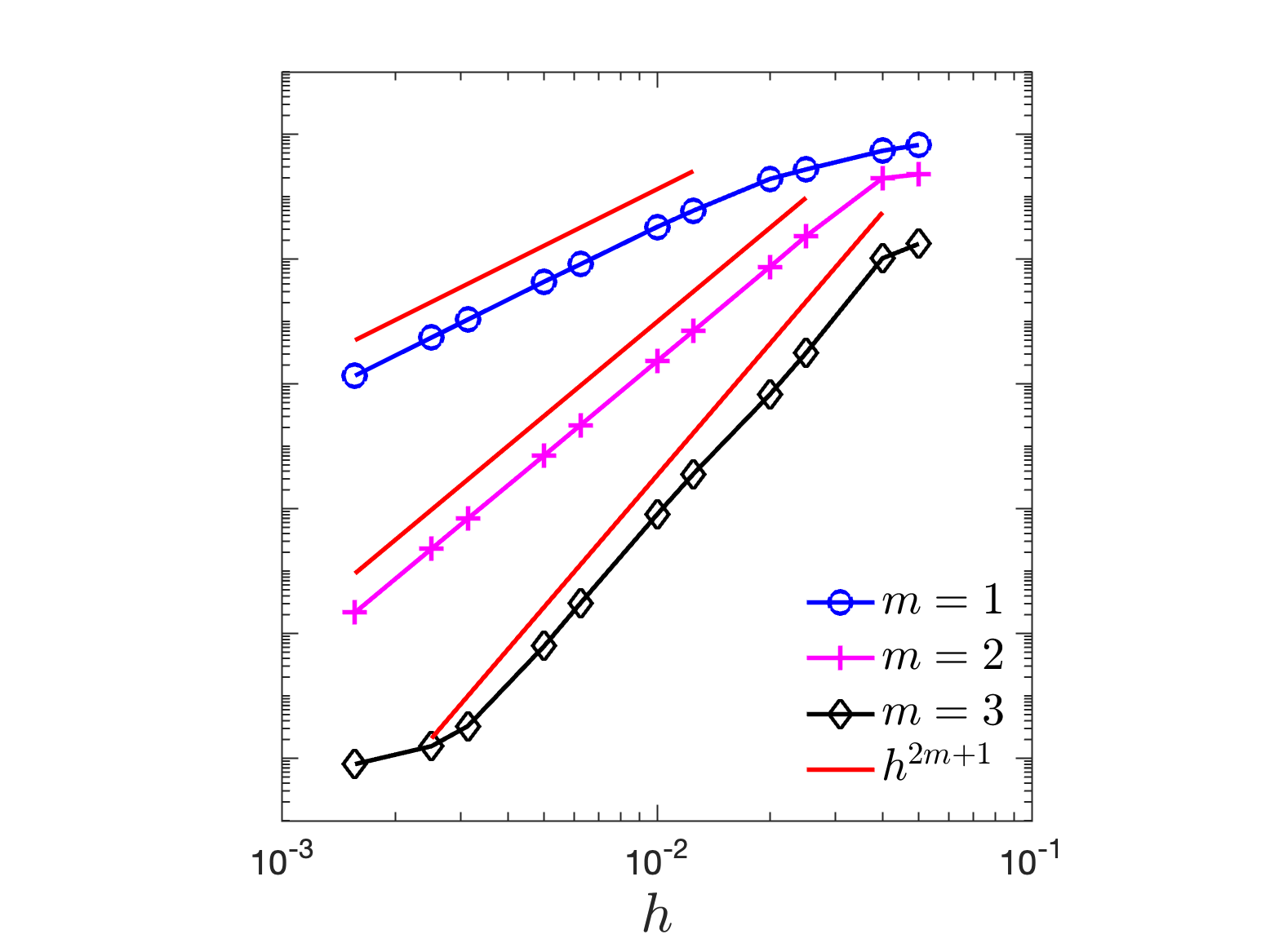} 
	\end{adjustbox}		
  \caption{Convergence plots in the maximum norm for a standing mode problem using for different values of $m$, 
  			boundary conditions and geometries in 2-D.
  		The left, 
			middle and right columns are respectively for the boundary conditions \eqref{eq:Ez_bnd_cdns_2D}, 
			\eqref{eq:H_bnd_cdns_2D} and \eqref{eq:impedance_bnd_cdns_2D}.
		The top and bottom rows are for the square and cross domains.  
		Here $\mathbold{U} = [H_x,H_y,E_z]^T$.}
   \label{fig:convPlots_2D_EM}
\end{figure}
Fig.~\ref{fig:convPlots_2D_divH} illustrates convergence plots for the divergence-free constraint on the magnetic field.
We observe a $2\,m$ rate of convergence as expected. 
 \begin{figure} 
 \centering
	\begin{adjustbox}{max width=1.0\textwidth,center}
	 \centering
		\includegraphics[width=2.5in,trim={1.5cm 0cm 1.75cm 0cm},clip]{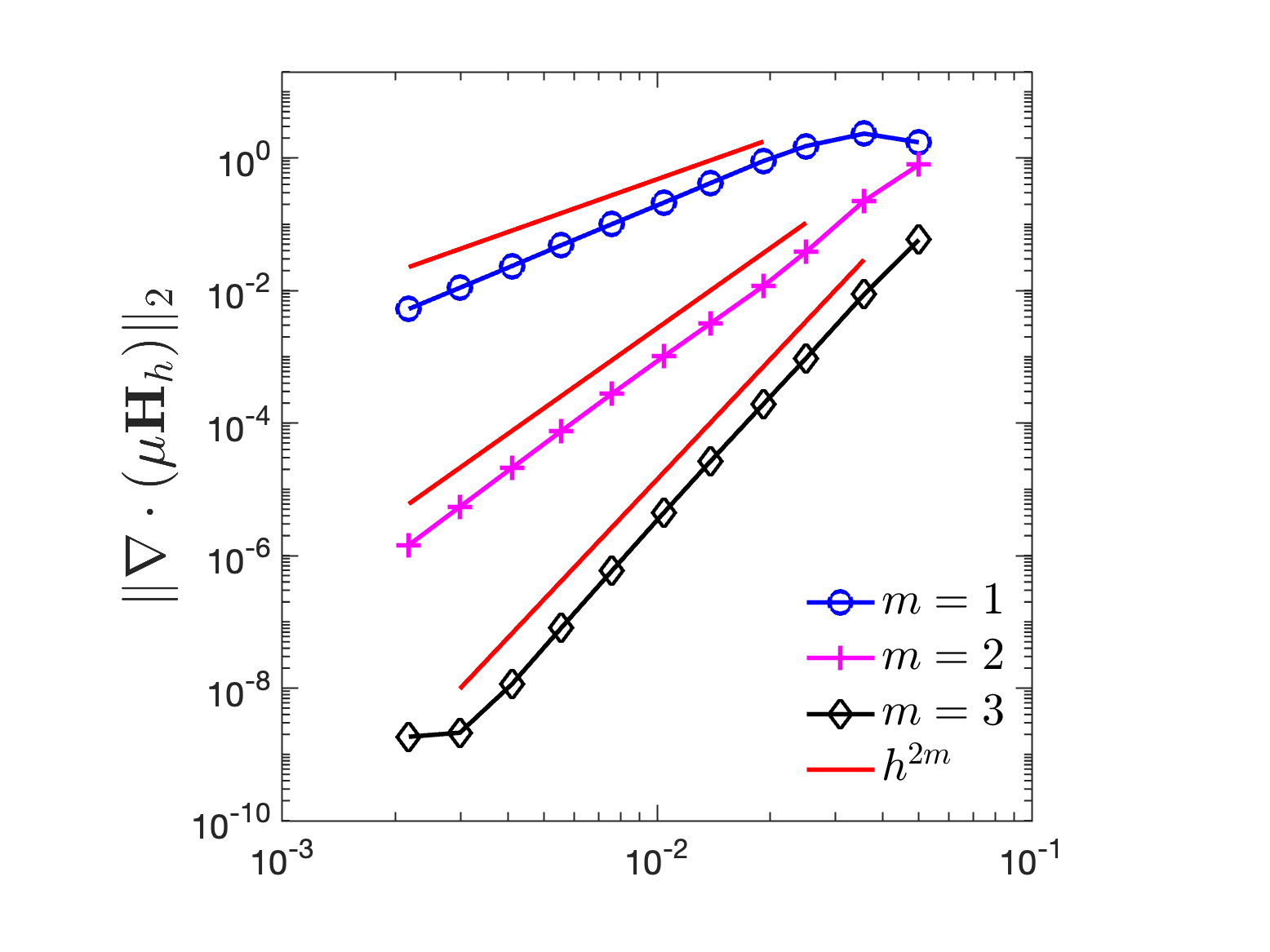}
		\includegraphics[width=2.5in,trim={1.5cm 0cm 1.75cm 0cm},clip]{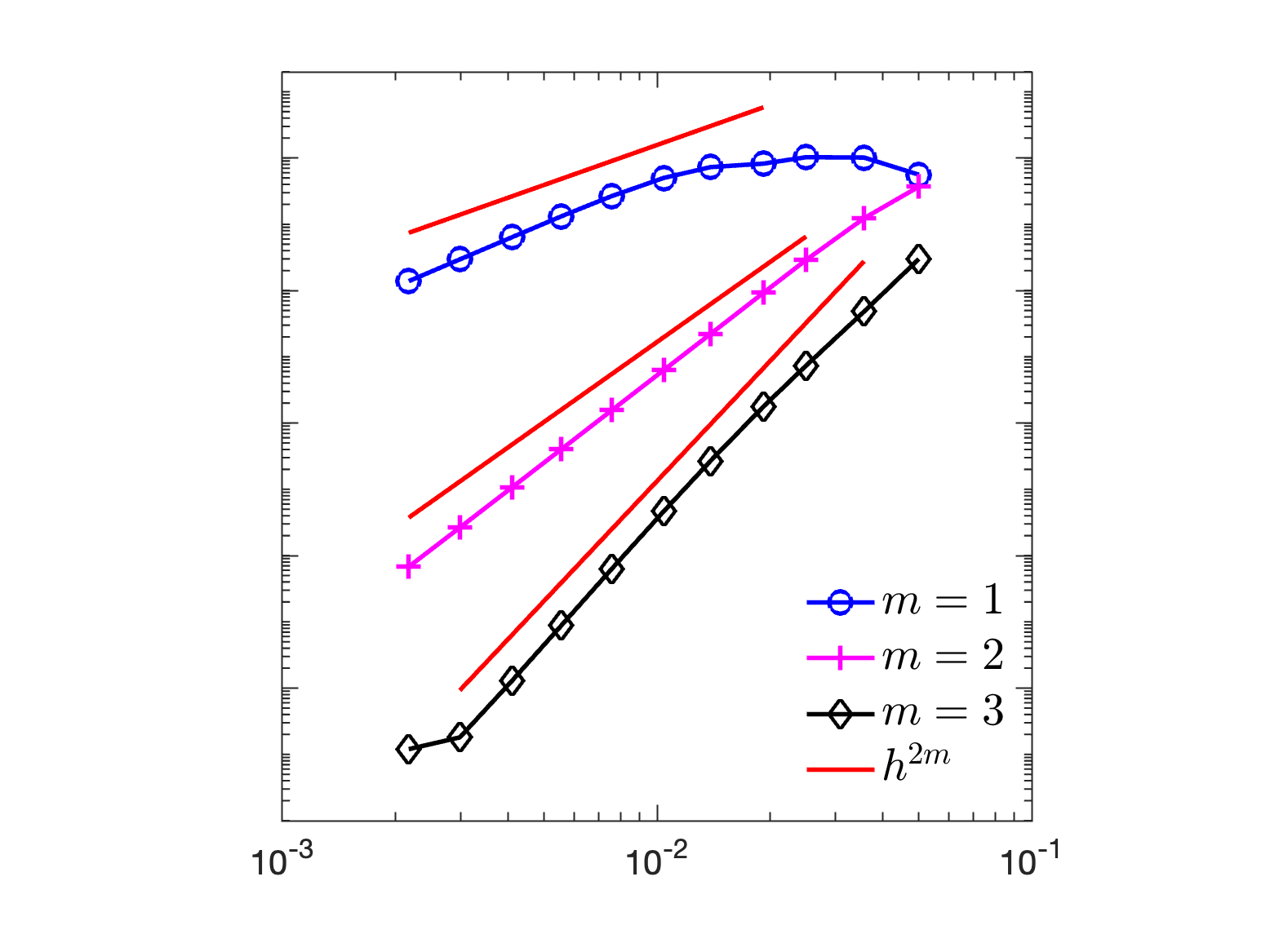}
	   	\includegraphics[width=2.5in,trim={1.5cm 0cm 1.75cm 0cm},clip]{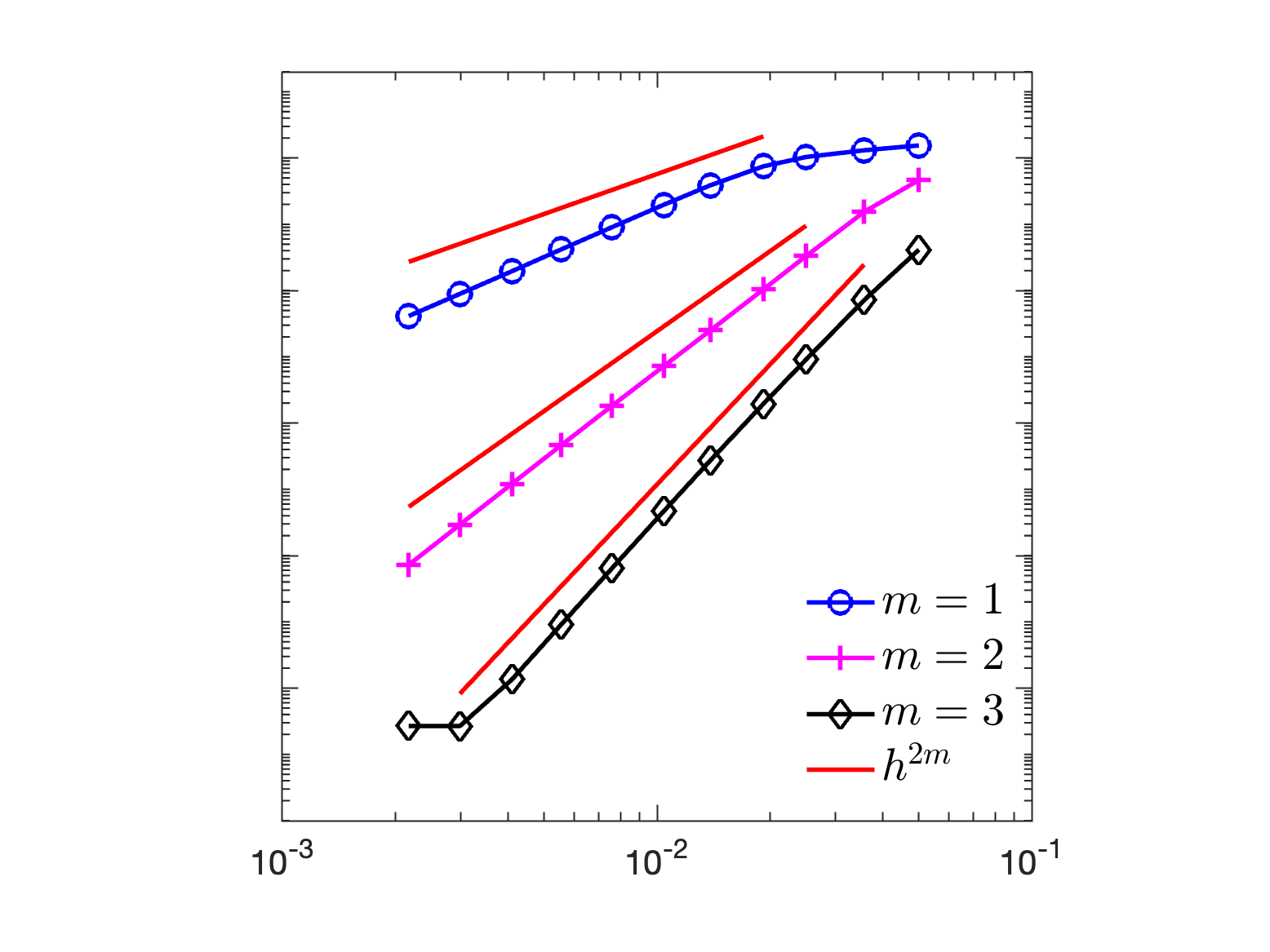} 
	\end{adjustbox}	
	\begin{adjustbox}{max width=1.0\textwidth,center}
	 \centering
		\includegraphics[width=2.5in,trim={1.5cm 0cm 1.75cm 0cm},clip]{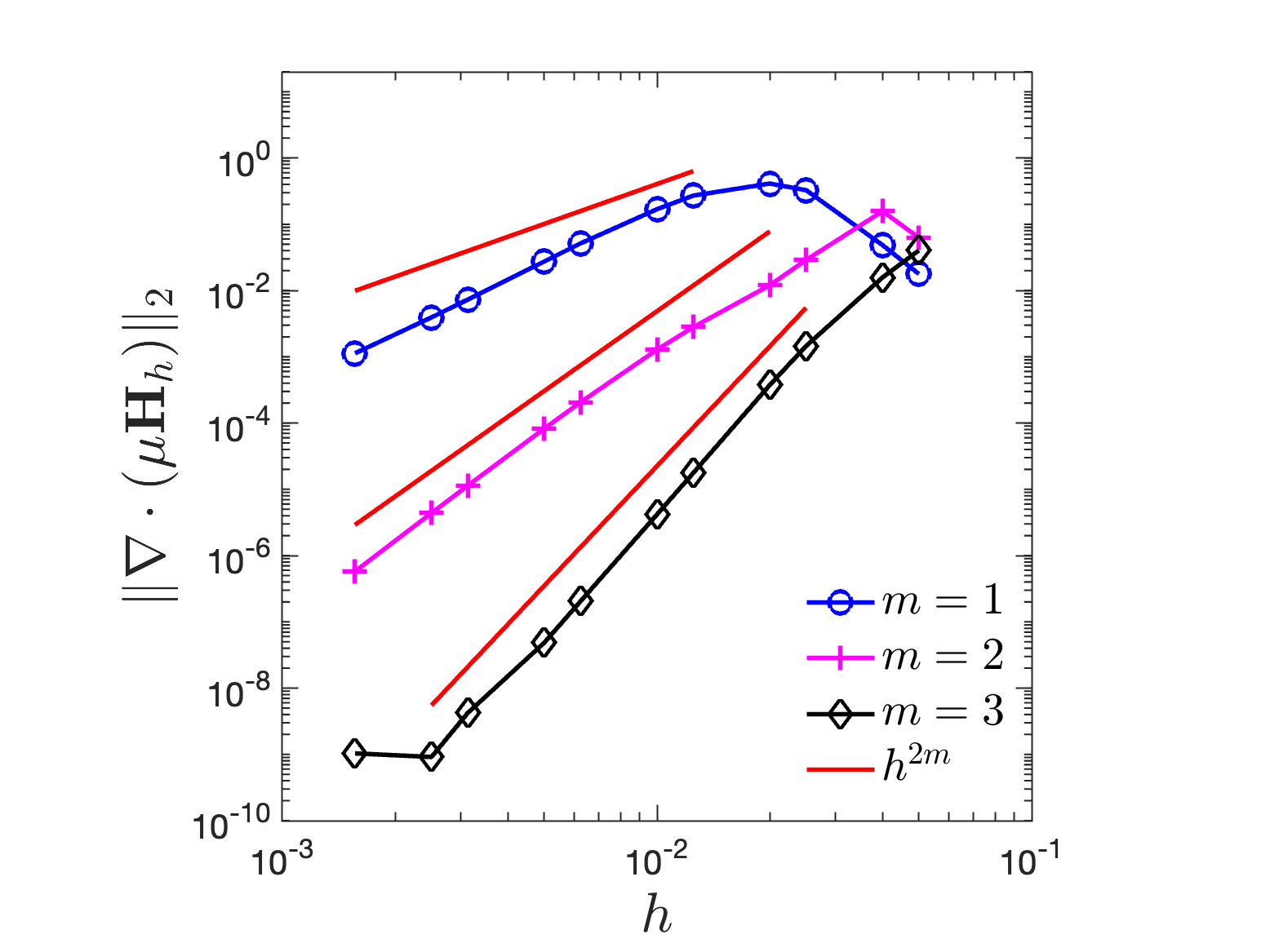}
		\includegraphics[width=2.5in,trim={1.5cm 0cm 1.75cm 0cm},clip]{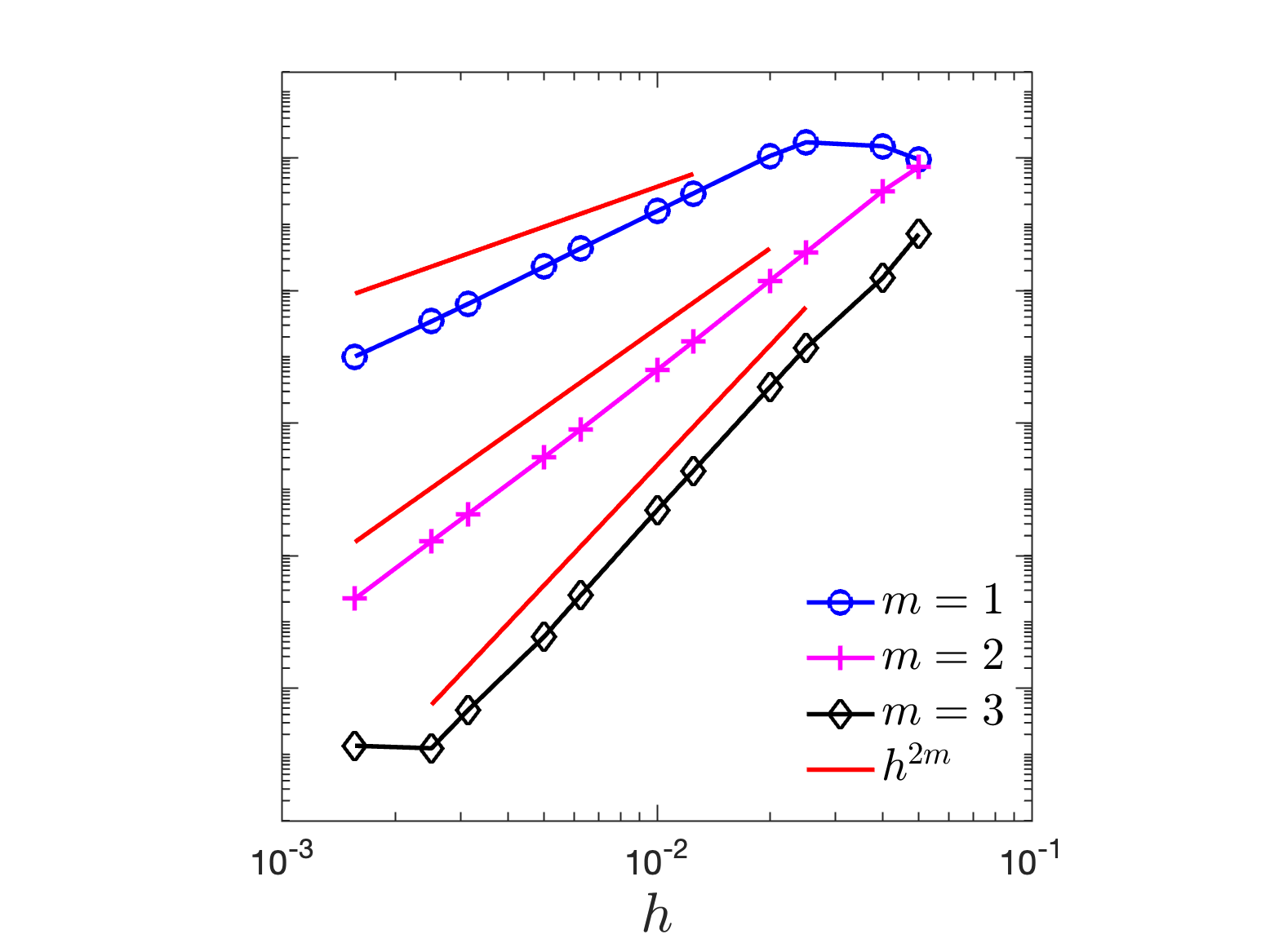}
	   	\includegraphics[width=2.5in,trim={1.5cm 0cm 1.75cm 0cm},clip]{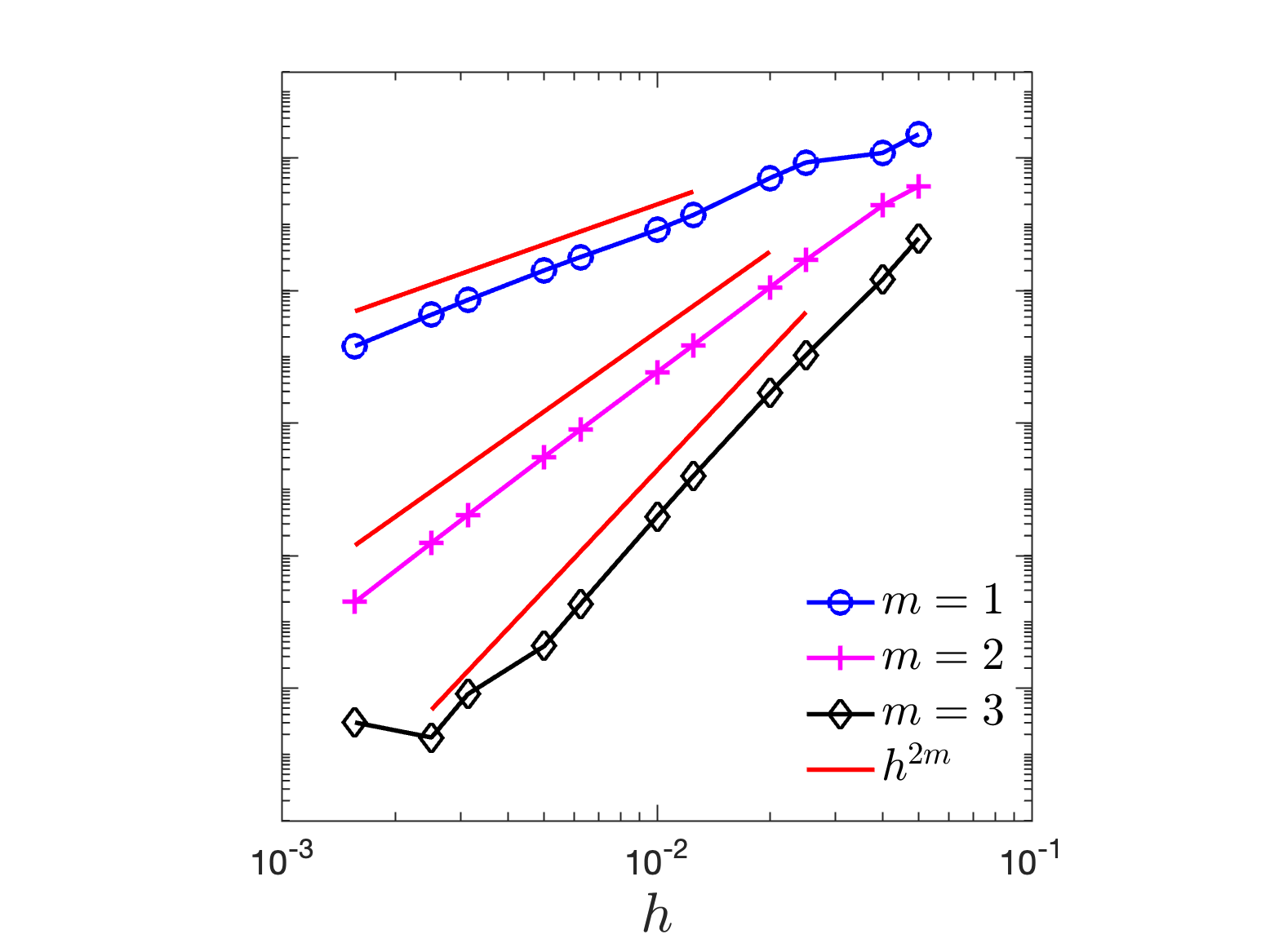} 
	\end{adjustbox}		
  \caption{Convergence plots of the divergence of the magnetic field in the $L^2$ norm for a standing mode problem using for different values of $m$, 
  			boundary conditions and geometries in 2-D.
  		The left, 
			middle and right columns are respectively for the boundary conditions \eqref{eq:Ez_bnd_cdns_2D}, 
			\eqref{eq:H_bnd_cdns_2D} and \eqref{eq:impedance_bnd_cdns_2D}.
		The top and bottom rows are for the square and cross domains.  
		Here $\mathbold{H} = [H_x,H_y]^T$.}
   \label{fig:convPlots_2D_divH}
\end{figure}
}

Let us now consider an initial Gaussian pulse on the electric field and PEC boundary conditions, 
	that is $E_z=0$, 
	on all the boundary of the domain. 
The square domain $\Omega = [0,1]\times[0,1]$ and the cross domain are considered with the time domain $I=[0,2]$.
The initial conditions are given by $H_x=H_y=0$ and 
\begin{equation*}
	E_z = e^{- \frac{r^2}{2\,\sigma^2}}.
\end{equation*}
Here $r^2 = (x-0.5)^2+(y-0.5)^2$ and $\sigma = 0.035$.
We set $\mu=1$ and $\epsilon=1$.

To our knowledge,
	there is no known analytic solution for this problem. 
Hence, 
	we perform self-convergence studies.
The reference solution $\mathbold{U}^* = [H_x^*,H_y^*,E_z^*]^T$
	is computed using the {\color{black} seventh-order} Hermite-Taylor correction function method with 
	$h=\frac{1}{800}$. 
We use meshes with $h=\big\{\frac{1}{25}, \frac{1}{50}, \frac{1}{100}, \frac{1}{200}, \frac{1}{400}\big\}$, 
	so all nodes used in the coarser meshes are also part of the reference solution mesh. 

\begin{figure} 
\begin{adjustbox}{max width=1.0\textwidth,center}
 \centering
	\includegraphics[width=3.0in]{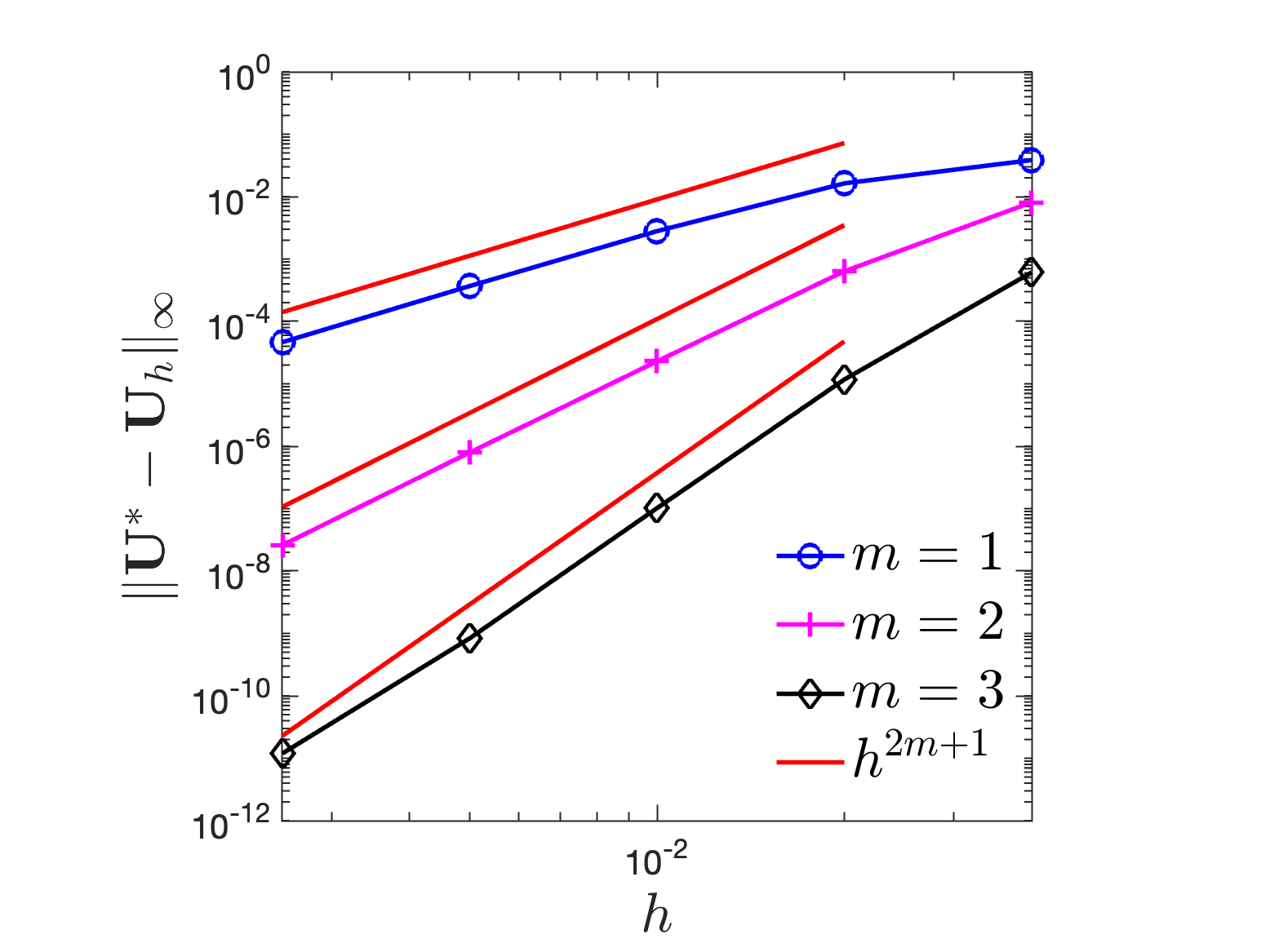}
	\includegraphics[width=3.0in]{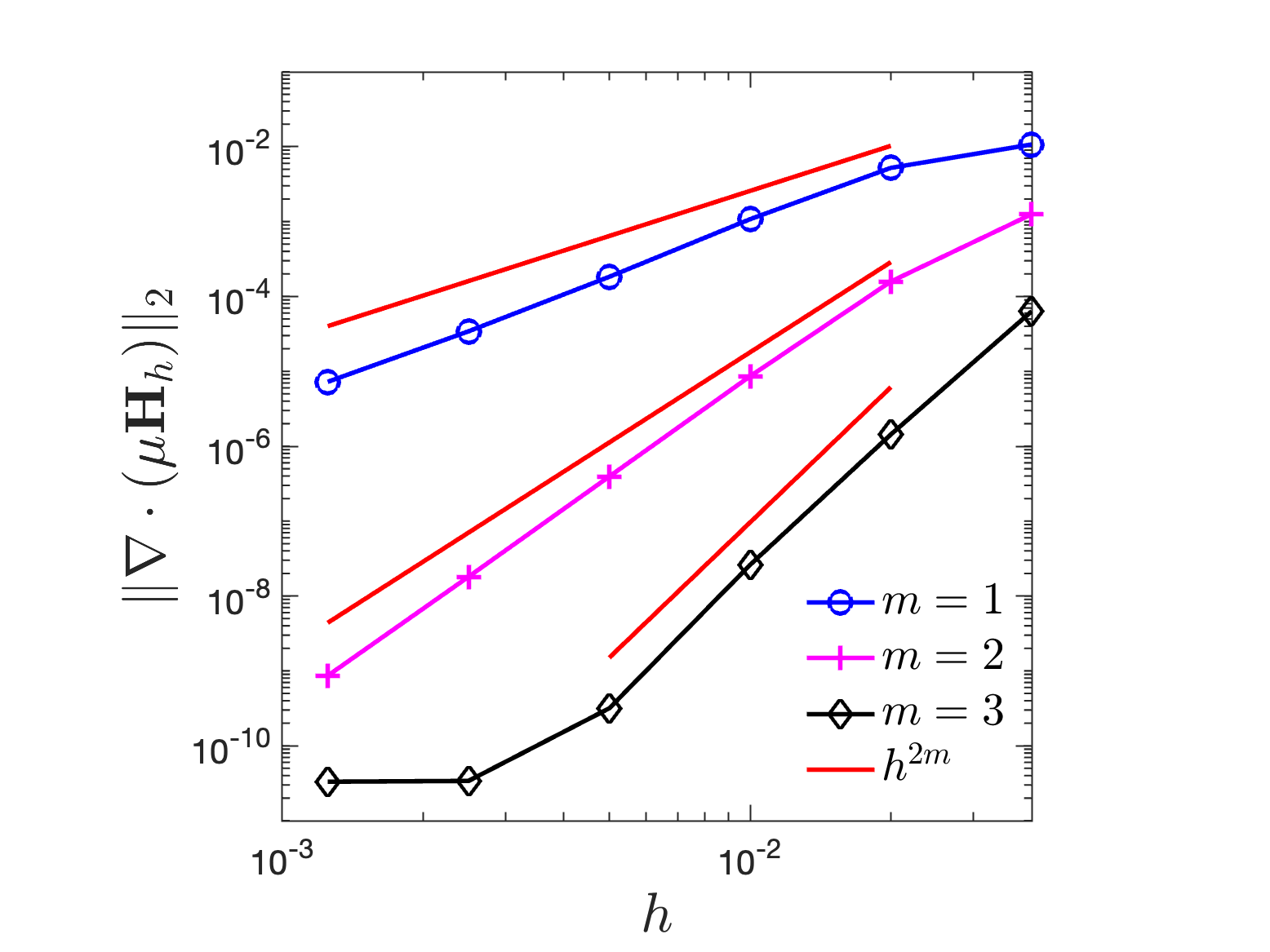}
\end{adjustbox}	
  \caption{Self-convergence plots for a Gaussian pulse problem using for different values of $m$ and the square domain in 2-D. 
  	The left and right plots are respectively for the error on the electromagnetic fields in the maximum norm 
			and the error on the divergence of the magnetic field in the $L^2$ norm.
		Here $\mathbold{U} = [H_x,H_y,E_z]^T$ and $\mathbold{H} = [H_x,H_y]^T$. }
   \label{fig:self_conv_plots_square}
\end{figure}

\begin{figure}   
	\centering
	\begin{adjustbox}{max width=1.0\textwidth,center}
		\includegraphics[width=2.5in,trim={2.5cm 0cm 1.75cm 0cm},clip]{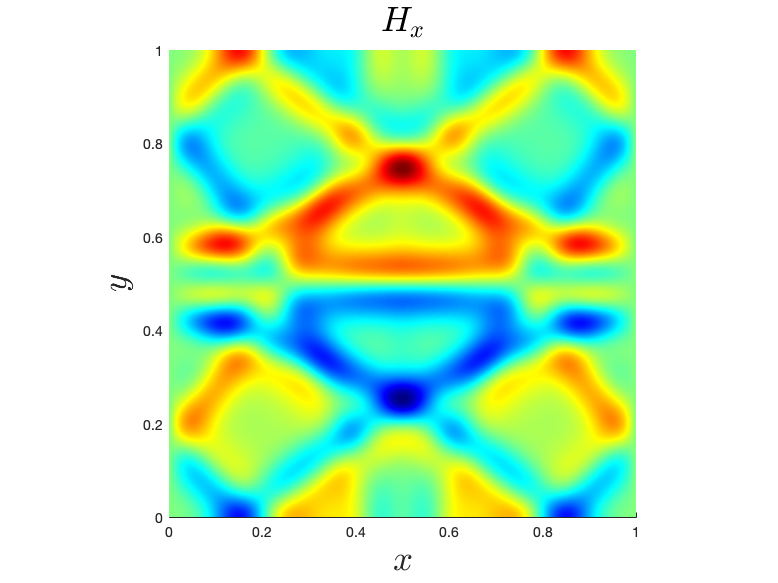}
		\includegraphics[width=2.5in,trim={2.5cm 0cm 1.75cm 0cm},clip]{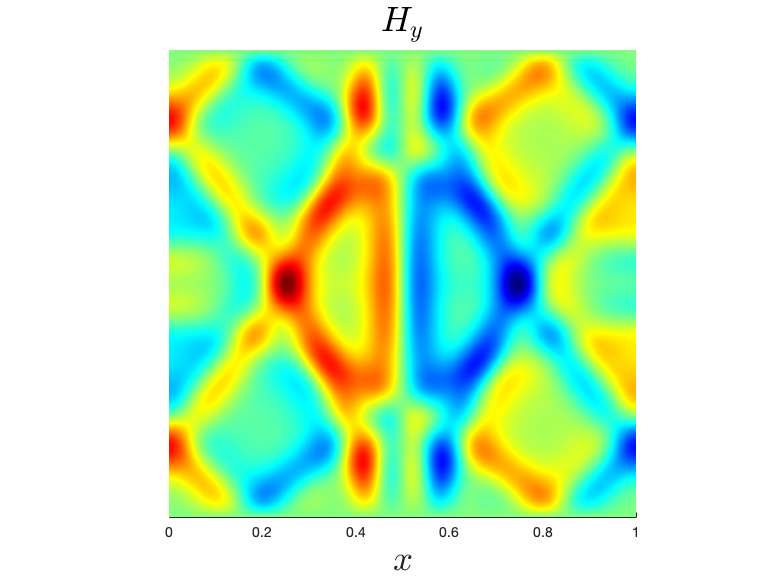}
	   	\includegraphics[width=2.5in,trim={2.5cm 0cm 1.75cm 0cm},clip]{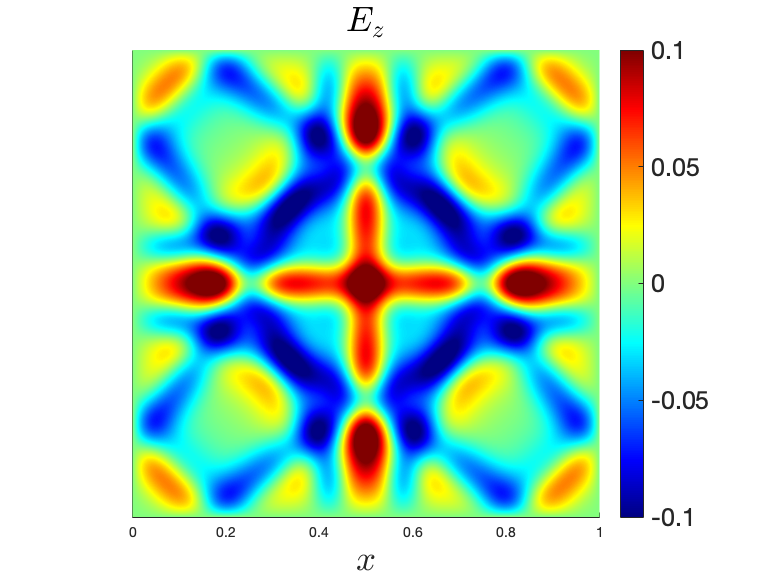}
	\end{adjustbox} 
       \caption{The components $H_x$, $H_y$ and $E_z$ for a Gaussian pulse problem using a square domain with the seventh-order Hermite-Taylor correction function method and $h=\frac{1}{800}$ at the final time.}
       \label{fig:square_gaussian_pulse_electromagnetic_fields}
\end{figure}

Fig.~\ref{fig:self_conv_plots_square} illustrates the self-convergence plots for the square domain. 
We obtain the expected $(2\,m+1)$ order of convergence for the electromagnetic fields, 
	while a $2\,m$ order of convergence is observed for the magnetic field divergence. 
The reference electromagnetic fields at the final time are shown in Fig.~\ref{fig:square_gaussian_pulse_electromagnetic_fields}.

\begin{figure} 
\begin{adjustbox}{max width=1.0\textwidth,center}
 \centering
	\includegraphics[width=3.0in]{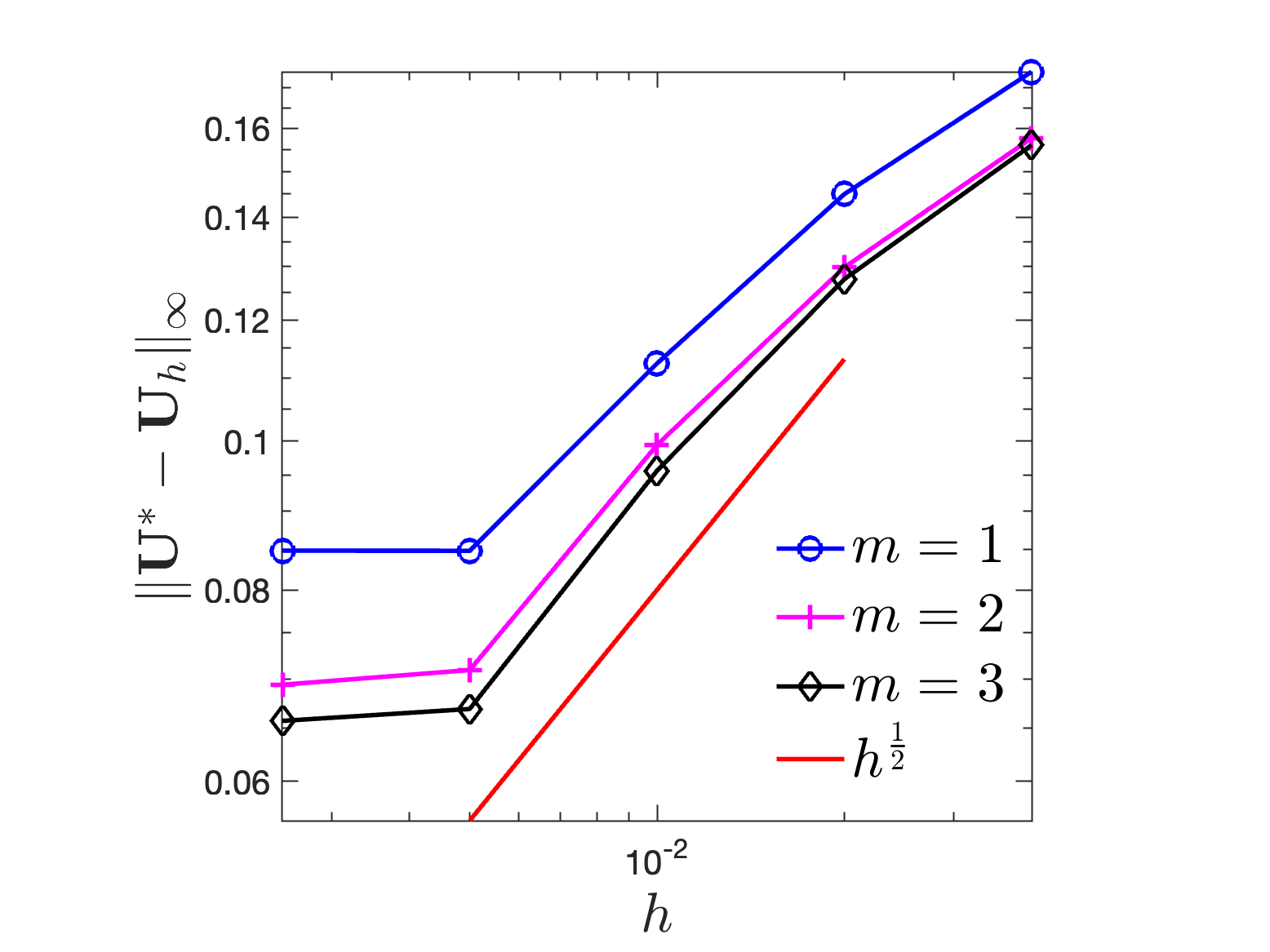}
	\includegraphics[width=3.0in]{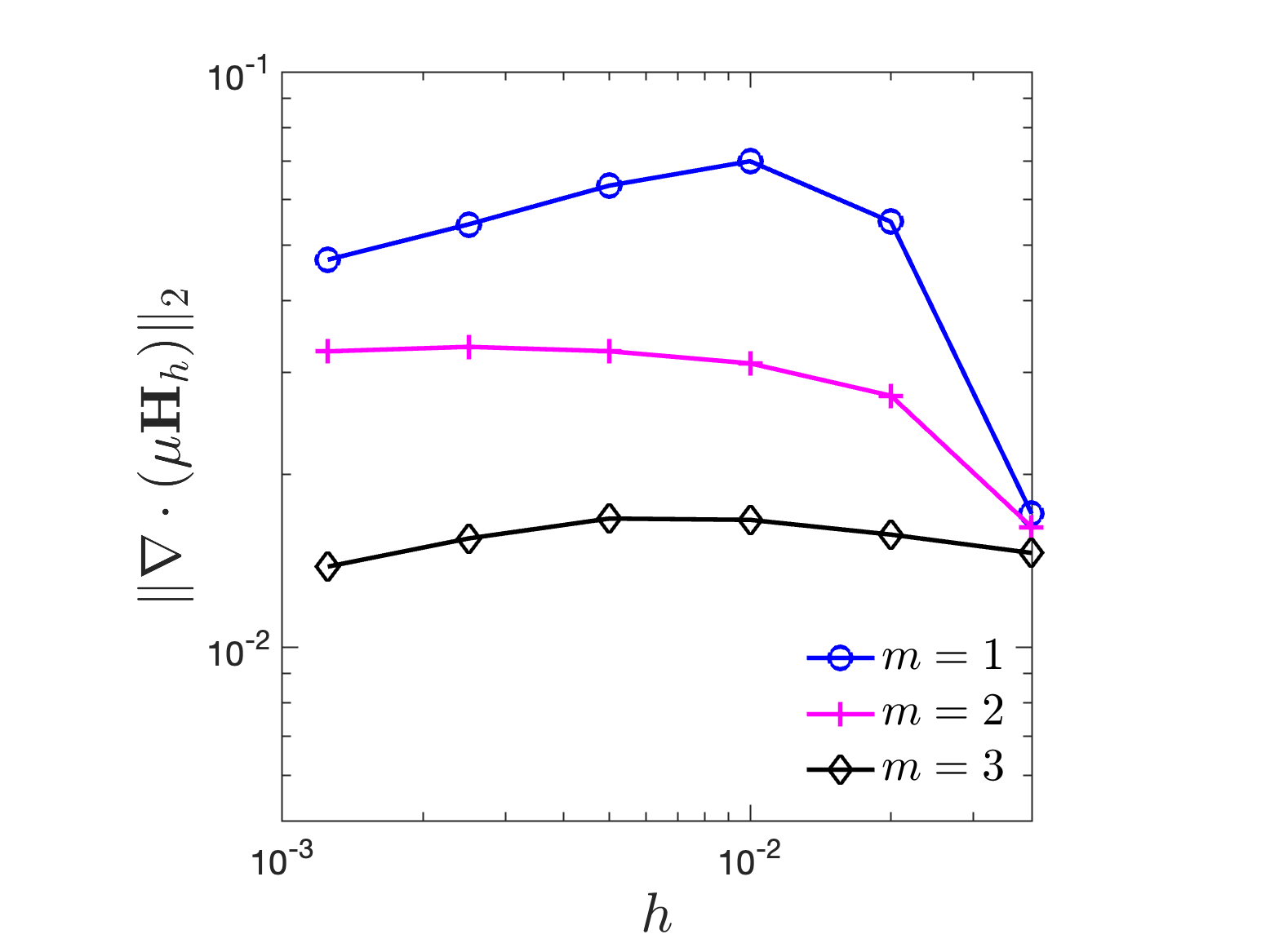}
\end{adjustbox}	
  \caption{Self-convergence plots for a Gaussian pulse problem using for different values of $m$ and the cross domain in 2-D. 
  	The left and right plots are respectively for the error on the electromagnetic fields in the maximum norm 
			and the error on the divergence of the magnetic field in the $L^2$ norm.
		Here $\mathbold{U} = [H_x,H_y,E_z]^T$ and $\mathbold{H} = [H_x,H_y]^T$. }
   \label{fig:self_conv_plots_cross}
\end{figure}

\begin{figure}   
	\centering
	\begin{adjustbox}{max width=1.0\textwidth,center}
		\includegraphics[width=2.5in,trim={2.5cm 0cm 1.75cm 0cm},clip]{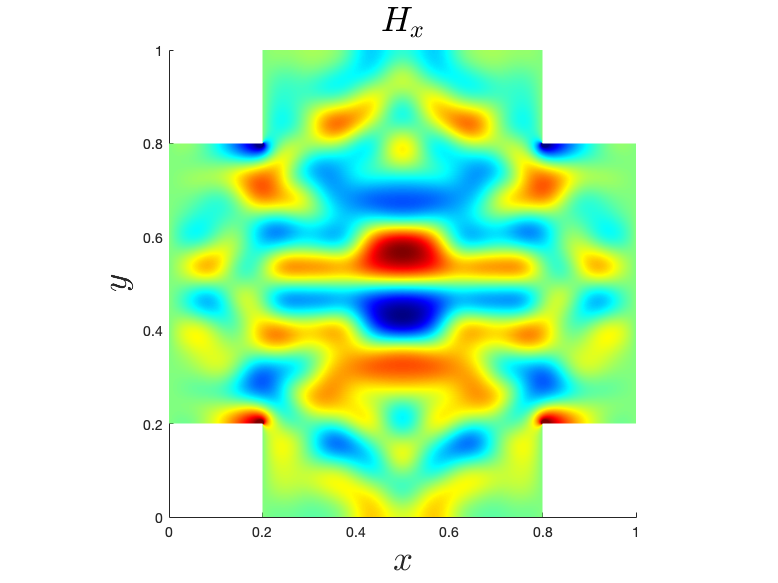}
		\includegraphics[width=2.5in,trim={2.5cm 0cm 1.75cm 0cm},clip]{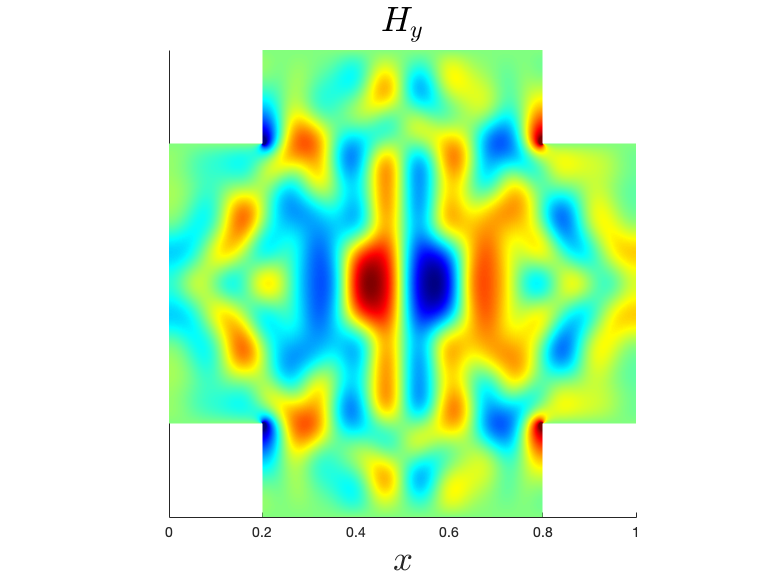}
	   	\includegraphics[width=2.5in,trim={2.5cm 0cm 1.75cm 0cm},clip]{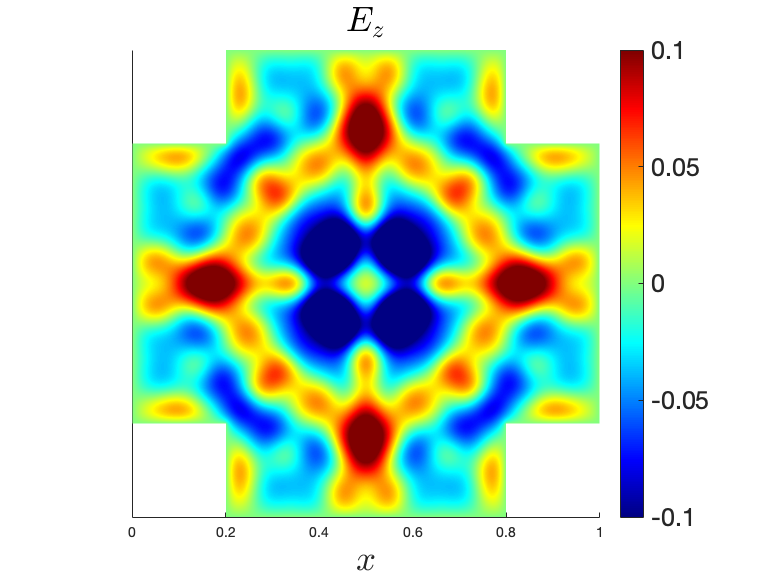}
	\end{adjustbox} 
       \caption{The components $H_x$, $H_y$ and $E_z$ for a Gaussian pulse problem using a cross domain with the seventh-order Hermite-Taylor correction function method and $h=\frac{1}{800}$ at the final time.}
       \label{fig:cross_gaussian_pulse_electromagnetic_fields}
\end{figure}

Fig.~\ref{fig:self_conv_plots_cross} illustrates the self-convergence plots for the cross domain. 
The numerical solution does not convergence in the maximum norm.
These results are explained by the reentrant corners in the cross domain. 
In this setting, 
	the solution has a singular part that hinders the performance of the numerical method \cite{Assous2000}.
The reference solution of the cross domain is shown in Fig.~\ref{fig:cross_gaussian_pulse_electromagnetic_fields},
	where strong variations in the magnetic field are observed at the reentrant corners.

{\color{black}
As a final numerical example, 
	we consider a variable coefficients problem.
In this situation,
	we use a manufactured solution given by
\begin{equation*}
	\begin{aligned}
	H_x =&\,\, -x\,e^{-x\,y}\,\sin(2\,\pi\,t), \\
	H_y =&\,\, y\,e^{-x\,y}\,\sin(2\,\pi\,t),\\
	E_z =&\,\, \sin(2\,\pi\,x\,y)\,\cos(2\,\pi\,t)
	\end{aligned}
\end{equation*}
	with $\mu(x,y) = \sin(5\,\pi\,x\,y)+2$ and $\epsilon(x,y) = 2\,e^{x\,y}$.
Note that source terms were considered in Maxwell's equations \eqref{eq:TMzSyst}. 
We consider the cross domain and the time interval $I=[0,1]$.
We enforce impedance boundary condition \eqref{eq:impedance_bnd_cdns_2D}.
Fig.~\ref{fig:convPlots_2D_var_coeff} illustrates the convergence plots for the electromagnetic fields and the divergence of 
	the magnetic field for different values of $m$.
 \begin{figure} 
 \centering
	\begin{adjustbox}{max width=1.0\textwidth,center}
	 \centering
		\includegraphics[width=2.5in,trim={1.5cm 0cm 1.75cm 0cm},clip]{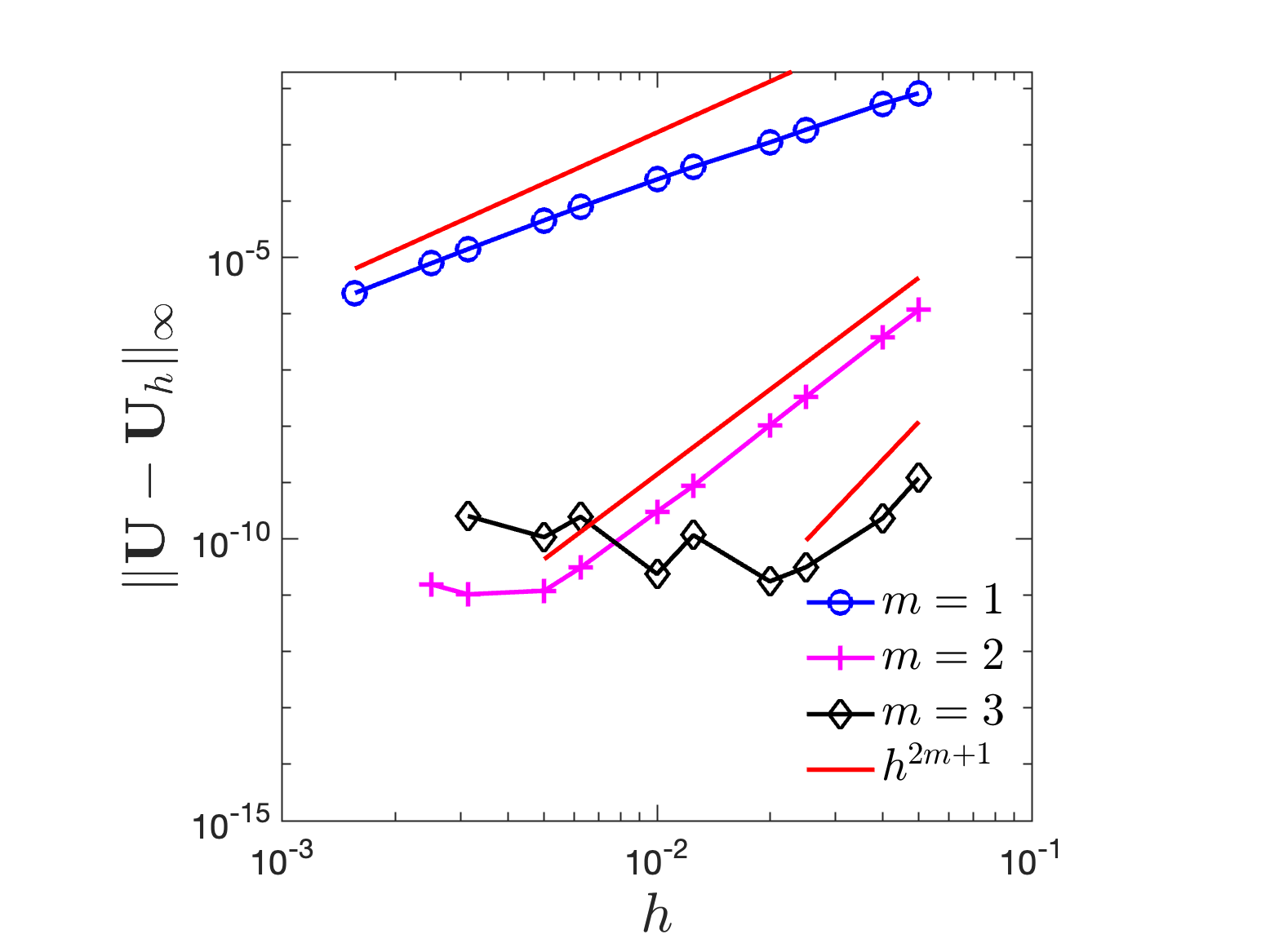}
	   	\includegraphics[width=2.5in,trim={1.5cm 0cm 1.75cm 0cm},clip]{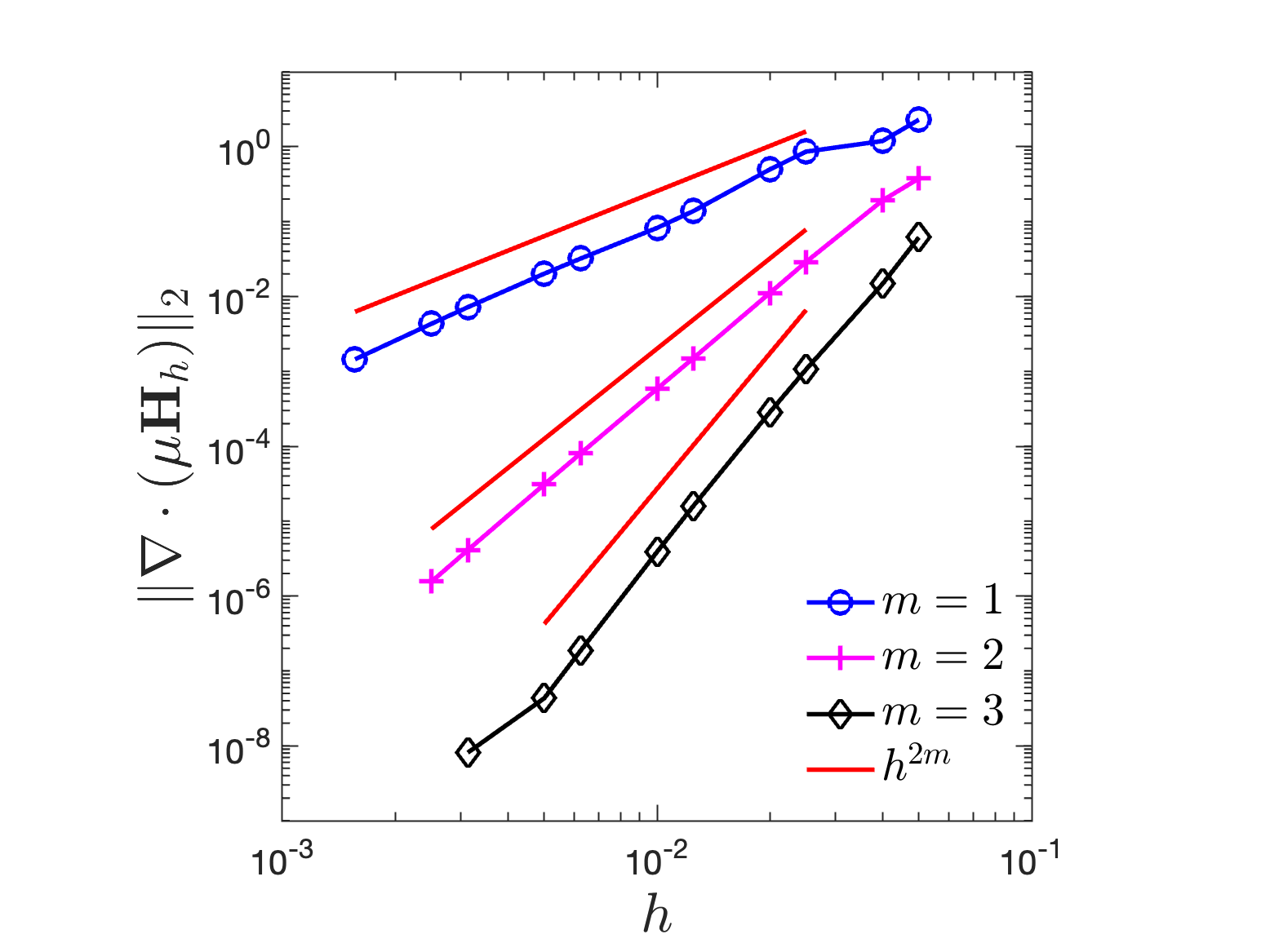} 
	\end{adjustbox}		
  \caption{Convergence plots for a manufactured solution problem with variable coefficients using for different values of $m$, 
  			the boundary condition \eqref{eq:impedance_bnd_cdns_2D} and the cross domain in 2-D.
  		The left and right plots are respectively for the error on the electromagnetic fields in the maximum norm 
			and the error on the divergence of the magnetic field in the $L^2$ norm.
		Here $\mathbold{U} = [H_x,H_y,E_z]^T$ and $\mathbold{H} = [H_x,H_y]^T$.}
   \label{fig:convPlots_2D_var_coeff}
\end{figure}
For $m=1$, 
	we observe a lower rate of convergence than expected. 
Based on previous numerical examples and the convergence of the magnetic field divergence, 
	finer meshes should exhibit the expected order of convergence.
As for $m=2$,
	we obtain the expected convergence order for the electromagnetic fields and the divergence-free constraint on the magnetic field.
Finally,
	the error on the electromagnetic fields with $m=3$ is already very low for coarser meshes, 
	making it hard to observe the seventh-order convergence of the method,
	while a clear sixth-order convergence is observed for the divergence of the magnetic field.
}

\section{Conclusion}

In this work, we have proposed a new method to handle boundary conditions for the Hermite-Taylor method for first order hyperbolic problems	based on the correction function method. Our method relies on a functional to be minimized that is a square measure of the residual of Maxwell's equations, the boundary conditions and the polynomial approximations of the electromagnetic fields coming from the Hermite-Taylor method. Once the minimization problems are solved, the information needed on the boundary, that is both electromagnetic fields and their {\color{black} space derivatives through order $m$}, are computed. {\color{black} Numerical examples suggest that the Hermite-Taylor correction function method is stable under a loose CFL constant and value of the penalization coefficient. Convergence rates of the Hermite-Taylor correction function method have been verified in 1-D and 2-D with different boundary conditions and geometries of the domain. Future work will focus on embedded boundary and interface problems.} 

\section*{Declarations}
\section*{Funding} 
This work was supported in part by Grant NSF- 2208164 and 2210286. Any opinions, findings, and conclusions or recommendations expressed in this material are those of the authors and do not necessarily reflect the views of the NSF.
\section*{Conflicts of interest/Competing interests}
On behalf of all authors, the corresponding author states that there is no conflict of interest.

%
%

\bibliographystyle{spmpsci}      
\bibliography{references}   


\end{document}